\documentclass[12pt,a4paper]{amsart}
\usepackage{amssymb}
\usepackage{mathrsfs}
\usepackage[cmtip,all]{xy}
\usepackage{hyperref}

\calclayout

\newcommand{\+}{\protect\nobreakdash-}
\renewcommand{\:}{\colon}

\newcommand{\rarrow}{\longrightarrow}
\newcommand{\larrow}{\longleftarrow}
\newcommand{\ot}{\otimes}
\newcommand{\ocn}{\odot}
\newcommand{\tim}{\rightthreetimes}

\newcommand{\lrarrow}{\mskip.5\thinmuskip\relbar\joinrel\relbar\joinrel
 \rightarrow\mskip.5\thinmuskip\relax}

\newcommand{\bu}{{\text{\smaller\smaller$\scriptstyle\bullet$}}}

\DeclareMathOperator{\Mor}{Mor}
\DeclareMathOperator{\Hom}{Hom}
\DeclareMathOperator{\Ext}{Ext}
\DeclareMathOperator{\Spec}{Spec}
\DeclareMathOperator{\Spf}{Spf}

\DeclareMathOperator{\Ind}{\mathsf{Ind}}

\newcommand{\Modl}{{\operatorname{\mathsf{--Mod}}}}
\newcommand{\Qcoh}{{\operatorname{\mathsf{--Qcoh}}}}
\newcommand{\Contra}{{\operatorname{\mathsf{--Contra}}}}
\newcommand{\Proqc}{{\operatorname{\mathsf{--Proqc}}}}

\newcommand{\fl}{\mathsf{fl}}
\newcommand{\vfl}{\mathsf{vfl}}
\newcommand{\cta}{\mathsf{cta}}
\renewcommand{\cot}{\mathsf{cot}}

\newcommand{\Sets}{\mathsf{Sets}}
\newcommand{\Sch}{\mathsf{Sch}}
\newcommand{\AffSch}{\mathsf{AffSch}}
\newcommand{\CSch}{\mathsf{CSch}}
\newcommand{\Com}{\mathsf{Com}}
\newcommand{\Ac}{\mathsf{Ac}}
\newcommand{\Fil}{\mathsf{Fil}}
\newcommand{\Sub}{\mathrm{Sub}}

\newcommand{\sop}{\mathsf{op}}
\newcommand{\red}{\mathrm{red}}

\newcommand{\ilim}
 {\mathop{\text{\normalfont``$\varinjlim$''\!\!}}\nolimits}

\newcommand{\rE}{\mathrm E}
\newcommand{\bP}{\mathbf P}
\newcommand{\bQ}{\mathbf Q}

\newcommand{\sD}{\mathsf D}
\newcommand{\sE}{\mathsf E}
\newcommand{\sF}{\mathsf F}
\newcommand{\sK}{\mathsf K}
\newcommand{\sL}{\mathsf L}
\newcommand{\sM}{\mathsf M}
\newcommand{\sP}{\mathsf P}
\newcommand{\sS}{\mathsf S}

\newcommand{\cO}{\mathcal O}
\newcommand{\cM}{\mathcal M}
\newcommand{\cI}{\mathcal I}
\newcommand{\cK}{\mathcal K}
\newcommand{\cF}{\mathcal F}
\newcommand{\cG}{\mathcal G}
\newcommand{\cH}{\mathcal H}
\newcommand{\cC}{\mathcal C}

\newcommand{\fX}{\mathfrak X}
\newcommand{\fY}{\mathfrak Y}
\newcommand{\fP}{\mathfrak P}
\newcommand{\fR}{\mathfrak R}
\newcommand{\fS}{\mathfrak S}
\newcommand{\fU}{\mathfrak U}
\newcommand{\fV}{\mathfrak V}
\newcommand{\fI}{\mathfrak I}
\newcommand{\fO}{\mathfrak O}
\newcommand{\fF}{\mathfrak F}

\newcommand{\ff}{\mathfrak f}
\newcommand{\fg}{\mathfrak g}
\newcommand{\fii}{\mathfrak i}
\newcommand{\fj}{\mathfrak j}

\newcommand{\rC}{\mathscr C}

\newcommand{\rF}{\mathscr F}
\newcommand{\rG}{\mathscr G}
\newcommand{\rH}{\mathscr H}
\newcommand{\rS}{\mathscr S}

\newcommand{\boP}{\mathbb P}
\newcommand{\boZ}{\mathbb Z}

\newcommand{\Section}[1]{\bigskip\section{#1}\medskip}

\theoremstyle{plain}
\newtheorem{thm}{Theorem}[section]
\newtheorem{prop}[thm]{Proposition}
\newtheorem{lem}[thm]{Lemma}
\newtheorem{cor}[thm]{Corollary}

\theoremstyle{definition}
\newtheorem{rem}[thm]{Remark}

\newtheorem{ex}[thm]{Example}
\newtheorem{exs}[thm]{Examples}

\begin{document}

\title{Flat and very flat pro-quasi-coherent \\
pro-sheaves on ind-schemes}

\author{Leonid Positselski}

\address{Institute of Mathematics, Czech Academy of Sciences \\
\v Zitn\'a~25, 115~67 Prague~1 \\ Czech Republic} 

\email{positselski@math.cas.cz}

\begin{abstract}
 This is a paper about homological algebra of flat
pro-quasi-coherent pro-sheaves.
 We recall and collect various techniques such as accessibility
and deconstructibility, and adopt them to ind-schemes under
various assumptions.
 We also discuss the foundations of the notion of an ind-scheme.
\end{abstract}

\maketitle

\tableofcontents

\section{Introduction}
\bigskip

 According to the discussion in~\cite[Section~7.11]{BD2}, there are
two kinds of quasi-coherent module objects over an ind-scheme.
 In the terminology of the present author's
book~\cite[Chapters~2\+-3]{Psemten}, they are called
the \emph{quasi-coherent torsion sheaves} and
the \emph{pro-quasi-coherent pro-sheaves}.
 In this paper, we study the pro-quasi-coherent pro-sheaves.

 The category of arbitrary pro-quasi-coherent pro-sheaves
on ind-scheme or a formal scheme is \emph{not} well-behaved from
the homological algebra standpoint~\cite[Section~7.11.3]{BD2},
\cite[Examples~3.1 and~3.8(1)]{Psemten}, \cite[Example~2.7(1)]{Pcta}.
 The first approximate explanation is that the functors of inverse
image of quasi-coherent sheaves with respect to closed immersions of
schemes are \emph{not} exact (and the functor of projective limit of
countable sequences of abelian groups is \emph{not} exact, either).
 In this paper, we restrict ourselves to the category of \emph{flat
pro-quasi-coherent pro-sheaves}, which is a well-behaved exact
category where one can do homological algebra.

 The identification of flat pro-quasi-coherent pro-sheaves on
a strict ind-affine $\aleph_0$\+ind-scheme $\fU=\Spf\fR$ with flat
contramodules over the related topological commutative ring $\fR$ was
pointed out in~\cite[Example~3.8(2)]{Psemten}.
 It became a geometric motivation for the study of flat contramodules
over complete, separated topological rings with a countable base of
neighborhoods of zero consisting of open two-sided ideals in
the paper~\cite[Sections~7--13]{Pflcc} (see~\cite[Example~8.9]{Pflcc}).
 In the present paper, go back to the geometric language and
generalize some contramodule-theoretic results of~\cite{Pflcc}
(as well as~\cite[Appendix~E]{Pcosh}, \cite{PR}, etc.)\ to flat
pro-quasi-coherent pro-sheaves over non-ind-affine ind-schemes.

 \emph{Ind-schemes} (and morphisms of ind-schemes) are represented
by filtered inductive systems of schemes, with the transition
morphisms usually presumed to be closed immersions.
 Various classes of ind-schemes and morphisms of ind-schemes are
conventionally defined by imposing conditions on the filtered
inductive systems representing such ind-schemes and morphisms.

 A subtle issue in the foundations of the theory arises immediately.
 Suppose that we are given an ind-scheme $\fX$, represented by
a directed diagram of schemes
$(Y_{\delta'}\to Y_{\delta''})_{\delta'\le\delta''\in\Delta}$,
and this directed diagram has some good property~$\bP$.
 Suppose further that the same ind-scheme $\fX$ is represented by
another directed diagram of schemes
$(Z_{\epsilon'}\to Z_{\epsilon''})_{\epsilon'\le\epsilon''\in\rE}$,
and that second directed diagram representing $\fX$ has another good
property~$\bQ$.
 Can one find a directed diagram
$(X_{\gamma'}\to X_{\gamma''})_{\gamma'\le\gamma''\in\Gamma}$
representing $\fX$ that has \emph{both the good properties\/ $\bP$
and\/~$\bQ$}\,?

 The same question applies to directed diagrams of morphisms of
schemes representing morphisms of ind-schemes, and their properties.
 Without a positive answer to this kind of questions, requiring
an ind-scheme or a morphism of ind-schemes to be ``both ind\+$\bP$ and
ind\+$\bQ$'' becomes ambiguous at best, or useless at worst.

 In Section~\ref{preliminaries-on-ind-schemes-secn} of the present
paper, which builds upon~\cite[Chapter~1]{Psemten}, we discuss
the foundations of the notion of an ind-scheme, concentrating in
particular on the class of questions described above.
 The exact category $\fX\Proqc_\fl$ of flat pro-quasi-coherent
pro-sheaves on an ind-scheme $\fX$ and its full subcategory
$\fX\Proqc_\vfl\subset\fX\Proqc_\fl$ of \emph{very flat}
pro-quasi-coherent pro-sheaves are defined in the next
Section~\ref{flat-and-vfl-proqcoh-prosheaves-secn}.

 In the subsequent Sections~\ref{accessibility-secn}\+-%
\ref{orthogonality-secn}, we recall various techniques from
(set-theoretical and other) homological algebra and apply them
to flat pro-quasi-coherent pro-sheaves.
 In particular, we prove that the category $\fX\Proqc_\fl$ is
accessible ($\aleph_1$\+accessible for an ind-quasi-compact
ind-quasi-separated $\aleph_0$\+ind-scheme~$\fX$), and that the classes
of flat and very flat pro-quasi-coherent pro-sheaves $\fX\Proqc_\fl$
and $\fX\Proqc_\vfl$ are deconstructible.
 We also spell out an argument proving $\Ext^1$\+orthogonality
of flat pro-quasi-coherent pro-sheaves with certain levelwise
conditions along the directed diagram of schemes representing $\fX$,
which was invented in the recent book
manuscript~\cite[Section~8.5]{Pform}.
 All of these results are also extended to some classes of
complexes of flat and very flat pro-quasi-coherent pro-sheaves.

\subsection*{Acknowledgement}
 I~am grateful to Jan \v St\!'ov\'\i\v cek and Michal Hrbek
for helpful discussions.
 The author is supported by the GA\v CR project 26-22734S
and the Institute of Mathematics, Czech Academy of Sciences
(research plan RVO:~67985840).

\Section{Preliminaries on Ind-Schemes}
\label{preliminaries-on-ind-schemes-secn}

\subsection{The basics of ind-schemes}
\label{basics-of-ind-schemes-subsecn}
 We can recommend the notes~\cite{Rich}, the unfinished
manuscript~\cite[Section~7.11]{BD2}, and the book~\cite{Psemten}
as the reference sources on ind-schemes.

 In this section, we consider diagrams indexed by directed posets.
 The terms ``filtered inductive system'', ``direct system'', and
``directed diagram'' are used as synonymous.
 A substantial discussion of the difference between and the essential
equivalence of the directed and the filtered colimits can be found in
the book~\cite[Theorem~1.5]{AR} and the paper~\cite{AN}.

 Given a category $\sK$, the category $\Ind \sK$ of \emph{ind-objects}
in $\sK$ is defined as the following full subcategory of the category
$\Sets^{\sK^\sop}$ of contravariant functors from $\sK$ to the category
of sets $\Sets$.
 To every object $K\in\sK$, the representable functor
$\Mor_\sK({-},K)\:\sK^\sop\rarrow\Sets$ is assigned (where
$\Mor_\sK({-},{-})$ denotes the set of morphisms in~$\sK$).
 The full subcategory $\Ind\sK\subset\Sets^{\sK^\sop}$ consists of
all the filtered inductive limits of representable functors.

 Let $\Gamma$ be a directed poset and
$(K_\gamma\to K_\delta)_{\gamma<\delta\in\Gamma}$ be
a $\Gamma$\+indexed diagram in~$\sK$.
 Then the corresponding ind-object, denoted by
$\ilim_{\gamma\in\Gamma}K_\gamma\in\Ind\sK$, represents the functor
$\sK^\sop\rarrow\Sets$ given by the rule
$L\longmapsto\varinjlim_{\gamma\in\Gamma}\Mor_\sK(L,K_\gamma)$.
 Accordingly, the sets of morphisms in $\Ind\sK$ are computed
by the formula
\begin{equation} \label{morphisms-of-ind-objects}
 \Mor_{\Ind\sK}(\ilim_{\delta\in\Delta}L_\delta,\>
 \ilim_{\gamma\in\Gamma}K_\gamma)=
 \varprojlim\nolimits_{\delta\in\Delta}
 \varinjlim\nolimits_{\gamma\in\Gamma}\Mor_\sK(L_\delta,K_\gamma)
\end{equation}
for any two directed posets $\Delta$ and $\Gamma$ and any two
diagrams $(L_\delta\in\sK)_{\delta\in\Delta}$ and
$(K_\gamma\in\sK)_{\gamma\in\Gamma}$ in the category~$\sK$.

 The standard references on ind-objects in categories
are~\cite[Section~I.8]{SGA4} and~\cite[Section~6]{KS}.
 See also~\cite[Section~1.1]{Psemten}, \cite[Section~9]{Pextop},
and/or~\cite{Hen}.

\begin{rem} \label{ind-objects-in-schemes-problematic-remark}
 Somewhat confusingly, ind-schemes are \emph{not} defined as
ind-objects in the category of schemes~$\Sch$.
 This can be explained as follows.

 Let $\Xi$ be an infinite set and $(U_\xi)_{\xi\in\Xi}$ be a family
of nonempty affine schemes indexed by~$\Xi$.
 Then the coproduct $X=\coprod_{\xi\in\Xi}U_\xi$ exists in
the category of schemes; this is simply the disjoint/disconnected
union of the schemes~$U_\xi$.
 So $X$ is not an affine scheme (and in fact, not a quasi-compact
scheme); but $X$ is still a scheme.

 On the other hand, denote by $\Gamma$ the directed poset of all
finite subsets of $\Xi$, ordered by inclusion.
 Then $\fX=\ilim_{\gamma\in\Gamma}\coprod_{\xi\in\gamma} U_\xi$ is
an ind-object in the category of schemes.
 In fact, $\fX$ is ind-object in the category of affine schemes, as it
is represented by the filtered inductive system of affine schemes
$\coprod_{\xi\in\gamma} U_\xi$.

 Now there is a natural morphism $\fX\rarrow X$ of ind-objects
in the category of schemes.
 But this morphism is \emph{never} an isomorphism, and in fact,
there is \emph{no} natural morphism from $X$ to $\fX$ in
the category of ind-objects in schemes.
 Indeed, for any scheme $Y$, any morphism $Y\rarrow\fX$ in
the category of ind-objects in schemes $\Ind(\Sch)$ would factorize
as $Y\rarrow\coprod_{\xi\in\gamma} U_\xi\rarrow\fX$ for some
$\gamma\in\Gamma$.

 So, say, the countable disjoint union of points $U_\xi=\Spec k$,
where $k$~is any field, gets represented by two quite different
objects $X$ and $\fX$ in the category $\Ind(\Sch)$ of ind-objects in
schemes.
 This is obviously an instance of undesirable behavior.
 The proper definition of an ind-scheme (as opposed to an ind-object
in schemes) is designed, in particular, to avoid this phenomenon.
\end{rem}

 Instead of thinking of ind-schemes as ind-objects in the category of
schemes, which means a certain kind of functors $\Sch^\sop\rarrow
\Sets$, the proper approach to defining ind-schemes is to view them as
functors $\AffSch^\sop\rarrow\Sets$ from the category of \emph{affine}
schemes $\AffSch$ to the category of sets.
 So, to any scheme $X$, the functor $\AffSch^\sop\rarrow\Sets$ given
by the rule $U\longmapsto\Mor_\Sch(U,X)$ is assigned.

 For a directed poset $\Gamma$ and a $\Gamma$\+indexed diagram of
schemes $(X_\gamma\to X_\delta)_{\gamma<\delta\in\Gamma}$, one
considers the functor $\AffSch^\sop\rarrow\Sets$ taking any affine
scheme $U$ to the set
$\varinjlim_{\gamma\in\Gamma}\Mor_\Sch(U,X_\gamma)$.
 An \emph{ind-scheme} is a functor $\AffSch^\sop\rarrow\Sets$
represented by some directed diagram of schemes
$(X_\gamma)_{\gamma\in\Gamma}$ in this way.

 An alternative approach uses functors on the category of concentrated
(i.~e., quasi-compact quasi-separated) schemes $\CSch$ instead of
the affine ones.
 For a directed poset $\Gamma$ and a $\Gamma$\+indexed diagram of
schemes $(X_\gamma\to X_\delta)_{\gamma<\delta\in\Gamma}$, consider
the functor $\CSch^\sop\rarrow\Sets$ given by the same rule
$Z\longmapsto\varinjlim_{\gamma\in\Gamma}\Mor_\Sch(Z,X_\gamma)$.
 Then the category of ind-schemes is equivalent to the category
of functors $\CSch^\sop\rarrow\Sets$ represented by directed diagrams
of schemes.

 The point is that, for any scheme $X$, the functor
$\CSch^\sop\rarrow\Sets$ represented by $X$ is a \emph{sheaf}
in the Zariski topology on the category of concentrated schemes.
 Moreover, the filtered inductive limit functors in the categories
of presheaves of sets and sheaves of sets on $\CSch$ agree, due to
the finite nature of the Zariski topology on concentrated schemes.
 Finally, the category of sheaves of sets on $\CSch$ is equivalent
to the category of sheaves of sets on $\AffSch$.
 This argument can be found in~\cite[Section~1.2]{Psemten}.

 It follows from the discussion in the preceding two paragraphs
that the full subcategory of ind-schemes representable by
filtered inductive systems of concentrated schemes is equivalent
to the category $\Ind(\CSch)$ of ind-objects in $\CSch$.
 Ind-schemes representable by filtered inductive systems of
concentrated schemes are said to be \emph{ind-concentrated}
(or \emph{ind-quasi-compact ind-quasi-separated}).
 In this paper, we are only interested in ind-concentrated ind-schemes.

 An ind-scheme is said to be \emph{strict} if it is representable
by a filtered inductive system
$(X_\gamma\to X_\delta)_{\gamma<\delta\in\Gamma}$ whose transition
morphisms $i_{\gamma\delta}\:X_\gamma\rarrow X_\delta$ are
\emph{closed immersions} of schemes.
 An \emph{$\aleph_0$\+ind-scheme} is an ind-scheme representable by
a direct system $(X_\gamma)_{\gamma\in\Gamma}$ indexed by
a \emph{countable} directed poset~$\Gamma$ (or equivalently, by
the directed poset of positive integers $\Gamma=\boZ_{\ge1}$).
 In this paper, we are only interested in strict
$\aleph_0$\+ind-schemes.

 The previous two paragraphs raise an immediate question.
 Is imposing such and other similar conditions on a filtered inductive
system $(X_\gamma\to X_\delta)_{\gamma<\delta\in\Gamma}$
\emph{separately} equivalent to imposing them \emph{simultaneously}?
 Given an ind-scheme $\fX$, suppose that it is represented by
a filtered inductive system of schemes $(Y_\delta)_{\gamma\in\Delta}$
having some good properties, and the same ind-scheme $\fX$ is also
represented by a filtered inductive system of schemes
$(Z_\lambda)_{\gamma\in\Lambda}$ having some other good properties.
 Does it follow that one can represent $\fX$ by a filtered inductive
system $(X_\gamma)_{\gamma\in\Gamma}$ \emph{combining both sets of 
good properties} at the same time?

 We continue the discussion of this issue started
in~\cite[Section~1.2]{Psemten}.
 It turns out that some of such questions can be answered positively
with the na\"\i ve definition of an ind-scheme as an object of
$\Ind(\Sch)$, but get negative answers with the more advanced
definition of an ind-scheme as an object of $\Sets^{\AffSch^\sop}$
or $\Sets^{\CSch^\sop}$.

 Let us start with restating a particular lemma from The Stacks project,
which is highly relevant to the basics of ind-schemes.

\begin{lem} \label{loc-closed-immersion-cancellation-lemma}
 Let $X\rarrow Y\rarrow Z$ be morphisms of schemes. \par
\textup{(a)} Assume that the composition $X\rarrow Z$ is a locally
closed immersion.
 Then the morphism $X\rarrow Y$ is a locally closed immersion. \par
\textup{(b)} Assume that the composition $X\rarrow Z$ is
a quasi-compact locally closed immersion, while the morphism
$Y\rarrow Z$ is quasi-separated.
 Then $X\rarrow Y$ is a quasi-compact locally closed immersion. \par
\textup{(c)} Assume that the composition $X\rarrow Z$ is a closed
immersion, while the morphism $Y\rarrow Z$ is separated.
 Then $X\rarrow Y$ is a closed immersion. \par
\textup{(d)} In particular, if the composition $X\rarrow Z$ is
a closed immersion and $Y$ is a quasi-separated scheme, then
the morphism $X\rarrow Y$ is a quasi-compact locally closed
immersion. \par
\textup{(e)} In particular, if the composition $X\rarrow Z$ is
a closed immersion and $Y$ is a separated scheme, then
the morphism $X\rarrow Y$ is a closed immersion.
\end{lem}

\begin{proof}
 Parts~(a\+-c) are~\cite[Definition Tag~01IO and Lemma Tag~07RK]{SP}.
 Parts~(a) and~(c) were also restated
in~\cite[Lemma~1.2(a\+-b)]{Psemten}.
 Part~(d) follows from~(b), since any closed immersion of schemes
is an affine (hence quasi-compact) morphism, and any morphism from
a quasi-separated scheme is quasi-separated~\cite[Lemma Tag~01KV]{SP}.
 Part~(e) similarly follows~(c); see~\cite[Lemma~1.2(c)]{Psemten}.
\end{proof}

\begin{cor} \label{Ind-of-Sch-strict-qs-qc-ss-lemma}
\textup{(a)} Suppose that\/ $\ilim_{\delta\in\Delta}Y_\delta\simeq
\ilim_{\lambda\in\Lambda}Z_\lambda$ in\/ $\Ind(\Sch)$, where
$(Y_\delta)_{\delta\in\Delta}$ is a direct system of quasi-separated
schemes, while $(Z_\lambda)_{\lambda\in\Lambda}$ is a direct system
of closed immersions of schemes.
 Then all the schemes $Z_\lambda$ are quasi-separated. \par
\textup{(b)} Suppose that\/ $\ilim_{\delta\in\Delta}Y_\delta\simeq
\ilim_{\lambda\in\Lambda}Z_\lambda$ in\/ $\Ind(\Sch)$, where
$(Y_\delta)_{\delta\in\Delta}$ is a direct system of semi-separated
schemes, while $(Z_\lambda)_{\lambda\in\Lambda}$ is a direct system
of closed immersions of schemes.
 Then all the schemes $Z_\lambda$ are semi-separated. \par
\textup{(c)} Suppose that\/ $\ilim_{\delta\in\Delta}Y_\delta\simeq
\ilim_{\lambda\in\Lambda}Z_\lambda$ in\/ $\Ind(\Sch)$, where
$(Y_\delta)_{\delta\in\Delta}$ is a direct system of Noetherian
schemes, while $(Z_\lambda)_{\lambda\in\Lambda}$ is a direct system
of closed immersions of schemes.
 Then all the schemes $Z_\lambda$ are Noetherian. \par
\textup{(d)} Suppose that\/ $\ilim_{\delta\in\Delta}Y_\delta\simeq
\ilim_{\lambda\in\Lambda}Z_\lambda$ in\/ $\Ind(\Sch)$, where
$(Y_\delta)_{\delta\in\Delta}$ is a direct system of quasi-compact
quasi-separated schemes, while $(Z_\lambda)_{\lambda\in\Lambda}$ is
a direct system of closed immersions of schemes.
 Then all the schemes $Z_\lambda$ are quasi-compact and
quasi-separated. \par
\textup{(e)} Suppose that\/ $\ilim_{\delta\in\Delta}Y_\delta\simeq
\ilim_{\lambda\in\Lambda}Z_\lambda$ in\/ $\Ind(\Sch)$, where
$(Y_\delta)_{\delta\in\Delta}$ is a direct system of affine schemes,
while $(Z_\lambda)_{\lambda\in\Lambda}$ is a direct system of closed
immersions of schemes.
 Then all the schemes $Z_\lambda$ are affine. \par
\end{cor}

\begin{proof}
 Put $\fX=\ilim_{\delta\in\Delta}Y_\delta=
\ilim_{\lambda\in\Lambda}Z_\lambda\in\Ind(\Sch)$.
 Let $\lambda'\in\Lambda$ be an index.
 By the definition of the category $\Ind(\Sch)$, the natural morphism
$Z_{\lambda'}\rarrow\fX$ factorizes through the natural morphism
$Y_\delta\rarrow\fX$ for some $\delta\in\Delta$, while the morphism
$Y_\delta\rarrow\fX$ factorizes through the natural morphism
$Z_{\lambda''}\rarrow\fX$ for some $\lambda''\in\Lambda$.
 Furthermore, all closed immersions of schemes are monomorphisms in
$\Sch$; so $(Z_\lambda)_{\lambda\in\Lambda}$ is a direct system of
monomorphisms in $\Sch$, and it follows that $Z_{\lambda''}\rarrow\fX$
is a monomorphism in $\Ind(\Sch)$.
 Hence the composition $Z_{\lambda'}\rarrow Y_\delta\rarrow
Z_{\lambda''}$ is equal to the transition morphism
$i_{\lambda'\lambda''}\:Z_{\lambda'}\rarrow Z_{\lambda''}$ in
the direct system $(Z_\lambda)_{\lambda\in\Lambda}$.

 By Lemma~\ref{loc-closed-immersion-cancellation-lemma}(a), it
follows that $Z_{\lambda'}\rarrow Y_\delta$ is a locally closed
immersion of schemes.
 Hence, if $Y_\delta$ is a quasi-separated scheme, then so
is~$Z_{\lambda'}$.
 For the same reason, if $Y_\delta$ is a semi-separated scheme,
then so is~$Z_{\lambda'}$.
 Similarly, if $Y_\delta$ is a Noetherian scheme, then so
is~$Z_{\lambda'}$.
 Furthermore, if $Y_\delta$ is a quasi-separated scheme, then
$Z_{\lambda'}\rarrow Y_\delta$ is a quasi-compact morphism by
Lemma~\ref{loc-closed-immersion-cancellation-lemma}(d).
 Hence, if $Y_\delta$ is a quasi-compact quasi-separated scheme,
then $Z_{\lambda'}$ is a quasi-compact quasi-separated scheme, too.
 Finally, if $Y_\delta$ is an affine scheme, then
$Z_{\lambda'}\rarrow Y_\delta$ is a closed immersion by
Lemma~\ref{loc-closed-immersion-cancellation-lemma}(e).
 Hence $Z_{\lambda'}$ is an affine scheme, too.
\end{proof}

 Parts~(b), (c), and~(e) of
Corollary~\ref{Ind-of-Sch-strict-qs-qc-ss-lemma} are also applicable
to direct systems of concentrated schemes and isomorphisms in
$\Ind(\CSch)$.
 Notice that, for any full subcategory $\sL\subset\sK$, the category
of ind-objects $\Ind\sL$ is a full subcategory in $\Ind\sK$.
 In particular, $\Ind(\CSch)$ is a full subcategory in $\Ind(\Sch)$;
so the assertions specialize to the case of concentrated schemes
automatically.

\begin{lem} \label{strict-aleph0-ind-schemes-lemma}
\textup{(a)} Suppose that\/ $\ilim_{\delta\in\Delta}Y_\delta\simeq
\ilim_{\lambda\in\Lambda}Z_\lambda$ in\/ $\Ind(\Sch)$, where
$(Y_\delta)_{\delta\in\Delta}$ is a direct system of schemes indexed
by a countable directed poset\/ $\Delta$,  while
$(Z_\lambda)_{\lambda\in\Lambda}$ is a direct system
of closed immersions of schemes.
 Then there exists a countable directed subposet\/ $\Lambda'\subset
\Lambda$ such that the induced morphism of ind-objects\/
$\ilim_{\lambda'\in\Lambda'}Z_{\lambda'}\rarrow
\ilim_{\lambda\in\Lambda}Z_\lambda$ is an isomorphism
in\/ $\Ind(\Sch)$. \par
\textup{(b)} Suppose that\/ $\ilim_{\delta\in\Delta}Y_\delta\simeq
\ilim_{\lambda\in\Lambda}Z_\lambda$ in\/ $\Ind(\CSch)$, where
$(Y_\delta)_{\delta\in\Delta}$ is a direct system of concentrated
schemes indexed by a countable directed poset\/ $\Delta$, while
$(Z_\lambda)_{\lambda\in\Lambda}$ is a direct system
of closed immersions of concentrated schemes.
 Then there exists a countable directed subposet\/ $\Lambda'\subset
\Lambda$ such that the induced morphism of ind-objects\/
$\ilim_{\lambda'\in\Lambda'}Z_{\lambda'}\rarrow
\ilim_{\lambda\in\Lambda}Z_\lambda$ is an isomorphism
in\/ $\Ind(\CSch)$.
\end{lem}

\begin{proof}
 Quite generally, suppose that $(Y_\delta)_{\delta\in\Delta}$ is
a direct system of objects indexed by a countable directed poset
$\Delta$ in a category $\sK$, while $(Z_\lambda)_{\lambda\in\Lambda}$
is a direct system of monomorphisms in $\sK$, and there is
an isomorphism $\ilim_{\delta\in\Delta}Y_\delta\simeq\fX=
\ilim_{\lambda\in\Lambda}Z_\lambda$ in\/ $\Ind(\sK)$.
 Then there exists a countable directed subposet\/ $\Lambda'\subset
\Lambda$ such that the induced morphism
$\ilim_{\lambda'\in\Lambda'}Z_{\lambda'}\rarrow
\ilim_{\lambda\in\Lambda}Z_\lambda$ is an isomorphism
in $\Ind(\sK)$.
 Specifically, given an index $\delta\in\Delta$, choose an index
$\lambda(\delta)\in\Lambda$ such that the natural morphism
$Y_\delta\rarrow\fX$ factorizes through the natural morphism
$Z_{\lambda(\delta)}\rarrow\fX$.
 Then any countable directed subposet $\Lambda'\subset\Lambda$
such that for every $\delta\in\Delta$ there exists
$\lambda'\in\Lambda'$ with $\lambda(\delta)\le\lambda'$ in $\Lambda$
has the desired property.
 Indeed, for any object $W\in\sK$, any morphism $W\rarrow\fX$ in
$\Ind(\sK)$ factorizes through $Y_\delta\rarrow\fX$ for some
$\delta\in\Delta$, hence also through $Z_{\lambda'}\rarrow\fX$
for some $\lambda'\in\Lambda'$.
 So, setting $\fX'=\ilim_{\lambda'\in\Lambda'}Z_{\lambda'}
\in\Ind(\sK)$, the map of sets $\Mor(W,\fX')\rarrow\Mor(W,\fX)$
is surjective.
 The injectivity of this map follows from the assumption that
$(Z_\lambda)_{\lambda\in\Lambda}$ is a direct system of monomorphisms.
\end{proof}

\begin{exs} \label{strict-ind-concentrated-counterexamples}
 The following counterexamples show that the notion of a strict
ind-quasi-compact ind-quasi-separated ind-scheme, viewed as an object
of $\Sets^{\AffSch^\sop}$ or $\Sets^{\CSch^\sop}$, is not as
well-behaved as one might hope for.

\smallskip
 (1)~Let $X=\bigcup_{\gamma\in\Gamma}U_\gamma$ be a quasi-separated
scheme representable as a union of its quasi-compact open subschemes
$U_\gamma$ indexed by a directed poset $\Gamma$ such that
$U_{\gamma'}\varsubsetneq U_{\gamma''}$ whenever
$\gamma'<\gamma''\in\Gamma$.
 Assume further that, for every $\gamma\in\Gamma$, the open
subscheme $U_\gamma$ is dense in~$X$.
 An example of a separated locally Noetherian scheme $X$ of this kind
(of Krull dimension~$2$ and with $\Gamma=\boZ_{\ge1}$) can be found
in~\cite[end of Section~II.7]{HarRD}.

 Denote by $\fX\in\Sets^{\CSch^\sop}$ the functor represented by~$X$.
 Then any morphism into $\fX$ from a quasi-compact quasi-separated
scheme $W$ factorizes uniquely as $W\rarrow U_\gamma\rarrow\fX$ for
large enough $\gamma\in\Gamma$; hence we have
$\fX=\ilim_{\gamma\in\Gamma}U_\gamma$ in $\Sets^{\CSch^\sop}$.
 On the other hand, one obviously has $\fX=\ilim_\Delta X$, where
$\Delta=\{*\}$ denotes the singleton poset.
 So we have a representation of the functor $\fX\in\Sets^{\CSch^\sop}$
by a directed diagram of quasi-compact quasi-separated
schemes~$U_\gamma$ (indexed by~$\Gamma$), and another representation
of $\fX$ by a directed diagram of closed immersions of quasi-separated
schemes (indexed by~$\Delta$).

 However, it is intuitively pretty clear that $\fX$ cannot be
represented by a directed diagram of closed immersions of
quasi-compact quasi-separated schemes.
 Indeed, assume for the sake of contradiction that
$\fX=\ilim_{\lambda\in\Lambda}Z_\lambda$ is such a representation.
 Pick an index $\gamma_0\in\Gamma$.
 Then the natural morphism $U_{\gamma_0}\rarrow\fX$ factorizes
through the natural morphism $Z_{\lambda_0}\rarrow\fX$ for some
$\lambda_0\in\Lambda$.
 Pick an index $\lambda'\ge\lambda_0$ in~$\Lambda$.
 Then the natural morphism $Z_{\lambda'}\rarrow X$
factorizes through the natural morphism $U_{\gamma'}\rarrow X$
for some $\gamma'\in\Gamma$, \ $\gamma'\ge\gamma_0$.
 Now we have a commutative square diagram of morphisms of schemes
$$
 \xymatrix{
  U_{\gamma_0} \ar@{>->}[r] \ar@{>->}[d] & U_{\gamma'} \\
  Z_{\lambda_0} \ar@{>->}[r]  & Z_{\lambda'} \ar@{>->}[u]
 }
$$
with an open immersion $U_{\gamma_0}\rarrow U_{\gamma'}$
and a closed immersion $Z_{\lambda_0}\rarrow Z_{\lambda'}$.
 The commutativity of the square diagram follows from the fact
that $(U_\gamma)_{\gamma\in\Gamma}$ is a direct system of
monomorphisms in $\CSch$.

 Arguing similarly to the proof of 
Corollary~\ref{Ind-of-Sch-strict-qs-qc-ss-lemma} and using
Lemma~\ref{loc-closed-immersion-cancellation-lemma}(a), we see
that $U_{\gamma_0}\rarrow Z_{\lambda_0}$ and
$Z_{\lambda'}\rarrow U_{\gamma'}$ are locally closed immersions
of schemes.
 We recall that, by assumption, $U_{\gamma_0}$ is dense in
$U_{\gamma'}$.
 Hence $Z_{\lambda_0}$ is dense in $Z_{\lambda'}$; so the closed
immersion of schemes $Z_{\lambda_0}\rarrow Z_{\lambda'}$
is bijective (in fact, a homeomorphism of the underlying
topological spaces).
 We have proved that this holds for all $\lambda'\ge\lambda_0$.

 Put $\mu_0=\lambda_0$, and pick an index $\mu'\in\Lambda$,
$\mu'\ge\mu_0$.
 Then the natural morphism $Z_{\mu_0}\rarrow X$
factorizes through the natural morphism $U_{\delta_0}\rarrow X$
for some $\delta_0\in\Gamma$, and the natural morphism
$Z_{\mu'}\rarrow X$ factorizes through the natural morphism
$U_{\delta'}\rarrow X$ for some $\delta'\in\Gamma$,
\,$\delta'\ge\delta_0$.
 Now we have another commutative square diagram of morphisms of
schemes
$$
 \xymatrix{
  U_{\delta_0} \ar@{>->}[r] & U_{\delta'} \\
  Z_{\mu_0} \ar@{>->}[r] \ar@{>->}[u]
  & Z_{\mu'} \ar@{>->}[u] \ar@{-->}[lu]
 }
$$
with an open immersion $U_{\delta_0}\rarrow U_{\delta'}$
and a bijective closed immersion $Z_{\mu_0}\rarrow Z_{\mu'}$.
 Once again, both $Z_{\mu_0}\rarrow U_{\delta_0}$ and
$Z_{\mu'}\rarrow U_{\delta'}$ are locally closed immersions
of schemes (but we do not need to use this fact).

 It follows that the morphism of schemes $Z_{\mu'}\rarrow
U_{\delta'}$ factorizes uniquely through the open immersion
$U_{\delta_0}\rarrow U_{\delta'}$.
 As this holds for all $\mu'\ge\mu_0$ and
$\fX=\ilim_{\mu'\ge\mu_0}Z_{\mu'}$, we arrive to the conclusion that
the morphism $U_{\delta_0}\rarrow\fX$ is an isomorphism; an obvious
contradiction.

\smallskip
 (2)~Let $U$ be an affine scheme with an open subscheme $V\subset U$
that is \emph{not} quasi-compact.
 Then $V$ can be represented as a union of quasi-compact open
subschemes $W_\gamma\subset V$ indexed by a directed poset $\Gamma$
such that $W_{\gamma'}\subset W_{\gamma''}$ whenever
$\gamma'\le\gamma''\in\Gamma$.
 Denote by $X$ the scheme obtained by gluing two copies of $U$
along $V$, and by $Y_\gamma$ the scheme obtained by gluing two copies
of $U$ along~$W_\gamma$.
 So the schemes $Y_\gamma$ are quasi-separated while the scheme $X$
is \emph{not} quasi-separated.
 Denote by $\Sets^{\CSch^\sop}$ the functor represented by~$X$.
 The schemes $Y_\gamma$ form a directed diagram indexed by $\Gamma$,
with a natural morphism $\ilim_{\gamma\in\Gamma}Y_\gamma\rarrow\fX$
in $\Sets^{\CSch^\sop}$.

 We claim that the morphism of functors
$\ilim_{\gamma\in\Gamma}Y_\gamma\rarrow\fX$ is an isomorphism
in $\Sets^{\CSch^\sop}$.
 Indeed, let $Z$ be a quasi-compact quasi-separated scheme, and
let $f\:Z\rarrow X$ be a morphism of schemes.
 Denote by $U'$, $U''\subset X$ the two copies of $U$ embedded
naturally into~$X$; so $U'\cap U''=V\subset X$.
 Then the open subschemes $f^{-1}(U')$ and $f^{-1}(U'')\subset Z$
are usually \emph{not} quasi-compact.
 Pick quasi-compact open coverings $f^{-1}(U')=\bigcup_\alpha
T'_\alpha$ and $f^{-1}(U'')=\bigcup_\beta T'_\beta$, and put
$T_{\alpha,\beta}=T'_\alpha\cup T''_\beta\subset Z$.
 Then $Z=\bigcup_{\alpha,\beta}T_{\alpha,\beta}$ is an open covering
of~$Z$.
 As such, it has a finite subcovering; so we have
$Z=\bigcup_{i=1}^mT_{\alpha_i,\beta_i}$ for some integer $m\ge0$
and some indices $\alpha_i$,~$\beta_i$.

 Now the intersections $T'_\alpha\cap T''_\beta\subset Z$ are
quasi-compact (since the scheme $Z$ is quasi-separated).
 Hence there exists an index $\gamma\in\Gamma$ such that
for all $1\le i\le m$ one has $f(T'_{\alpha_i}\cap T''_{\beta_i})
\subset W_\gamma\subset V=U'\cap U''$.

 Let $T$ be a scheme represented as the union of two open subschemes,
$T=T'\cup T''$.
 To specify a morphism of schemes $g\:T\rarrow X$ such that
$g(T')\subset U'$ and $g(T'')\subset U''$, it suffices to pinpoint two
morphisms $g'\:T'\rarrow U'$ and $g''\:T''\rarrow U''$ that both take
the intersection $T'\cap T''$ into $V=U'\cap U''\subset X$ and agree
on that intersection.
 Morphisms of schemes $h\:T'\cup T''\rarrow Y_\gamma$ are described
similarly, except that the intersection is now
$U'_\gamma\cap U''_\gamma=W_\gamma$, with the obvious notation
$U'_\gamma$ and $U''_\gamma$ for the two copies of $U$ embedded
naturally into~$Y_\gamma$.

 It follows easily that the morphism $Z\rarrow X$ factorizes
uniquely as $Z\rarrow Y_\gamma\rarrow X$, where the index
$\gamma$~satisfies the condition above, and the morphism
$h\:Z\rarrow Y_\gamma$ satisfies the conditions that
$h(T'_{\alpha_i})\subset U'_\gamma$ and
$h(T''_{\beta_i})\subset U''_\gamma$ for all $1\le i\le m$.
 
 This shows that the morphism of set-valued functors
$\ilim_{\gamma\in\Gamma}Y_\gamma\rarrow\fX$ is surjective.
 To prove the injectivity, suppose that we are given an index
$\delta\in\Gamma$ and two morphisms of schemes
$h_1$, $h_2\:Z\rarrow Y_\delta$ such that the two compositions
$Z\rightrightarrows Y_\delta\rarrow X$ agree.
 Denote these two equal compositions by $f\:Z\rarrow X$.
 For every index $\epsilon\ge\delta$, denote the compositions
$Z\rightrightarrows Y_\delta\rarrow Y_\epsilon$ by
$h_{1,\epsilon}$, $h_{2,\epsilon}\:Z\rarrow Y_\epsilon$.
 Then we have $f^{-1}(U')=\bigcup_{\epsilon\ge\delta}
h_{1,\epsilon}^{-1}(U'_\epsilon)=\bigcup_{\epsilon\ge\delta}
h_{2,\epsilon}^{-1}(U'_\epsilon)$ and $f^{-1}(U'')=
\bigcup_{\epsilon\ge\delta}h_{1,\epsilon}^{-1}(U''_\epsilon)=
\bigcup_{\epsilon\ge\delta}h_{2,\epsilon}^{-1}(U''_\epsilon)$.

 In the notation above, we have $f^{-1}(U')=\bigcup_\alpha
T'_\alpha$ and $f^{-1}(U'')=\bigcup_\beta T'_\beta$, and
then $Z=\bigcup_{i=1}^mT_{\alpha_i,\beta_i}$ for
$T_{\alpha,\beta}=T'_\alpha\cup T''_\beta$ and
some indices $\alpha_i$,~$\beta_i$, and finally
$f(T'_{\alpha_i}\cap T''_{\beta_i}) \subset W_\gamma\subset
V=U'\cap U''$ for all $1\le i\le m$ and some $\gamma\in\Gamma$.
 Let $\epsilon\in\Gamma$ be an index such that $\epsilon\ge\delta$
and $\epsilon\ge\gamma$, and also $T'_{\alpha_i}\subset
h_{1,\epsilon}^{-1}(U'_\epsilon)$, \
$T'_{\alpha_i}\subset h_{2,\epsilon}^{-1}(U'_\epsilon)$ and
$T''_{\beta_i}\subset h_{1,\epsilon}^{-1}(U''_\epsilon)$, \
$T''_{\beta_i}\subset h_{2,\epsilon}^{-1}(U''_\epsilon)$ for all
$1\le i\le m$.

 We have a natural morphism $Y_\epsilon\rarrow U$, which factorizes
naturally as $Y_\epsilon\rarrow X\rarrow U$.
 Now let $T$ be a scheme and $g_1$, $g_2\:T\rightrightarrows Y_\epsilon$
be two morphisms of schemes.
 To check that $g_1=g_2$, it suffices to cover $T$ by two open
subschemes, $T=T'\cup T''$, such that $g_j(T')\subset U'_\epsilon$
and $g_j(T'')\subset U''_\epsilon$ for $j=1$, $2$, and the two
compositions $T\rightrightarrows Y_\epsilon\rarrow U$ are equal
to each other.

 It follows easily, by applying this criterion, that
$h_{1,\epsilon}=h_{2,\epsilon}$ for the index $\epsilon\in\Gamma$
chosen as above.
 This proves that the morphism of set-valued functors
$\ilim_{\gamma\in\Gamma}Y_\gamma\rarrow\fX$ is also injective,
as desired.

 On the other hand, one obviously has $\fX=\ilim_\Delta X$,
where $\Delta=\{*\}$ is the singleton poset.
 Thus we have a representation of the functor
$\fX\in\Sets^{\CSch^\sop}$ by a directed diagram of quasi-compact
quasi-separated schemes~$Y_\gamma$ (indexed by~$\Gamma$), and
another representation of $\fX$ by a directed diagram of closed
immersions of quasi-compact schemes (indexed by~$\Delta$).
 We leave it as an exercise for the reader to check that $\fX$
cannot be represented by a directed diagram of closed immersions
of quasi-compact quasi-separated schemes.
\end{exs}

 Based on the discussion above in this section, for the rest of
this paper we will mean by ``ind-scheme'' an ind-object in
the category of concentrated schemes $\CSch$.
 A ``strict ind-scheme'' is an ind-scheme representable by
a directed diagram of closed immersions of concentrated schemes.
 An ``$\aleph_0$\+ind-scheme'' is an ind-scheme representable by
a countable directed diagram of concentrated schemes.
 By Lemma~\ref{strict-aleph0-ind-schemes-lemma}(b), any strict
$\aleph_0$\+ind-scheme is representable by a countable directed
diagram of closed immersions of concentrated schemes.

 We refer to~\cite[Section~1.2]{Psemten} for a discussion of
the issue of uniqueness of a directed diagram of closed immersions
of schemes representing a given strict ind-scheme.
 To give a precise formulation, the arguments
in~\cite[Section~1.2]{Psemten} prove that any isomorphism in
$\Ind(\CSch)$ or $\Ind(\Sch)$ between two ind-schemes represented
by directed diagrams of closed immersions of schemes is also
an isomorphism of ind-objects in the category of (concentrated or
arbitrary) schemes \emph{with closed immersions as morphisms}.
 No other morphisms of schemes than closed immersions can appear as
components of such an isomorphism, the ``components'' being understood
in the sense of the formula~\eqref{morphisms-of-ind-objects}.

 An ind-scheme is said to be \emph{ind-semi-separated} if it is
representable by a directed diagram of semi-separated schemes,
\emph{ind-affine} if it is representable by a directed diagram
of affine schemes, and \emph{ind-Noetherian} if it is representable
by a directed diagram of Noetherian schemes.
 By Corollary~\ref{Ind-of-Sch-strict-qs-qc-ss-lemma}(b,c,e) and
the subsequent paragraph, the notions of a strict ind-semi-separated
ind-scheme, strict ind-affine ind-scheme, and strict ind-Noetherian
ind-scheme are well-behaved, in the sense of the discussion above
in this section.
 So are the notions of a strict ind-semi-separated, strict
ind-affine, or strict ind-Noetherian $\aleph_0$\+ind-scheme,
by Lemma~\ref{strict-aleph0-ind-schemes-lemma}(b).

\subsection{Finite diagrams and fibered products of ind-schemes}
\label{fibered-products-of-ind-schemes-subsecn}
 Let $D$ be a finite poset.
 Given a category $\sK$, we denote by $\sK^D$ the category of
$D$\+indexed diagrams in~$\sK$.
 We are particularly interested in two special cases: for the ordered
set of two elements $D=(\bullet\to\bullet)$, the category $\sK^D$ is
the category of morphisms in~$\sK$; while for the poset
$D=(\bullet\rightarrow\bullet\leftarrow\bullet)$, the category
$\sK^D$ is the category of \emph{cospans} in~$\sK$.

 The limits of cospans are called the \emph{pullbacks} or
the \emph{fibered products}.
 We denote the fibered product of a cospan $X\rarrow S\larrow Y$
by $X\times_SY$.

\begin{prop} \label{finite-diagrams-in-ind-objects}
 Let\/ $\sK$ be a category and $D$ be a finite poset.
 Then the natural functor\/ $\Ind(\sK^D)\rarrow\Ind(\sK)^D$
is a category equivalence.
\end{prop}

\begin{proof}
 This is~\cite[Expos\'e~I, Proposition~8.8.5]{SGA4},
\cite[page~55]{Mey}, or a special case of~\cite[Theorem~1.3]{Hen}.
\end{proof}

 In particular, by Proposition~\ref{finite-diagrams-in-ind-objects},
any morphism $\ff\:\fY\rarrow\fX$ in $\Ind\sK$ can be represented by
a directed diagram of morphisms $(f_\xi\: Y_\xi\to X_\xi)_{\xi\in\Xi}$
in $\sK$, indexed by some directed poset~$\Xi$.
 So we have $\ff=\ilim_{\xi\in\Xi}f_\xi$.
 Given two directed diagrams $(X_\gamma)_{\gamma\in\Gamma}$ and
$(Y_\delta)_{\delta\in\Delta}$ in $\sK$ representing the ind-objects
$\fX$ and $\fY\in\Ind\sK$, an explicit construction of a directed poset
$\Xi$ and a $\Xi$\+indexed diagram of morphisms $(f_\xi\:Y_\xi\to
X_\xi)$ in $\sK$ representing the morphism~$\ff$ can be found, e.~g.,
in~\cite[Section~9]{Pextop}.

 The nature of the explicit construction mentioned in the previous
paragraph is such that all the objects in the diagrams
$(X_\xi)_{\xi\in\Xi}$ and $(Y_\xi)_{\xi\in\Xi}$ appear already in
the diagrams $(X_\gamma)_{\gamma\in\Gamma}$ and
$(Y_\delta)_{\delta\in\Delta}$, and all the transition morphisms in
the diagrams $(X_\xi)_{\xi\in\Xi}$ and $(Y_\xi)_{\xi\in\Xi}$ appear
already in the diagrams $(X_\gamma)_{\gamma\in\Gamma}$ and
$(Y_\delta)_{\delta\in\Delta}$, respectively.
 If $(X_\gamma)_{\gamma\in\Gamma}$ and $(Y_\delta)_{\delta\in\Delta}$
are \emph{countable diagrams of monomorphisms} in $\sK$, then so
are the diagrams $(X_\xi)_{\xi\in\Xi}$ and $(Y_\xi)_{\xi\in\Xi}$.

 Furthermore, let $\fX\rarrow\fS\larrow\fY$ be a cospan in $\Ind\sK$.
 By Proposition~\ref{finite-diagrams-in-ind-objects}, the cospan
$\fX\rarrow\fS\larrow\fY$ can be represented by a directed diagram of
cospans $(X_\xi\rightarrow S_\xi\leftarrow Y_\xi)$ in $\sK$,
indexed by some directed poset~$\Xi$.
 Given three directed diagrams $(X_\gamma)_{\gamma\in\Gamma}$,
$(Y_\delta)_{\delta\in\Delta}$, and $(S_\epsilon)_{\epsilon\in\rE}$
in $\sK$ representing the ind-objects $\fX$, $\fY$, and $\fS$,
an explicit construction of a directed poset $\Xi$ and a $\Xi$\+indexed
diagram of cospans $(X_\xi\rightarrow S_\xi\leftarrow Y_\xi)$ in $\sK$
representing the cospan $\fX\rarrow\fS\larrow\fY$ in $\Ind\sK$ can be
found, e.~g., in~\cite[Section~1.1]{Psemten}.

 Once again, all the objects in the diagrams $(X_\xi)_{\xi\in\Xi}$,
$(Y_\xi)_{\xi\in\Xi}$, and $(S_\xi)_{\xi\in\Xi}$ appear already in
the diagrams $(X_\gamma)_{\gamma\in\Gamma}$,
$(Y_\delta)_{\delta\in\Delta}$, and $(S_\epsilon)_{\epsilon\in\rE}$,
and all the transition morphisms in the diagrams $(X_\xi)_{\xi\in\Xi}$,
$(Y_\xi)_{\xi\in\Xi}$, and $(S_\xi)_{\xi\in\Xi}$ appear already
in the diagrams $(X_\gamma)_{\gamma\in\Gamma}$,
$(Y_\delta)_{\delta\in\Delta}$, and $(S_\epsilon)_{\epsilon\in\rE}$,
respectively.
 If $(X_\gamma)_{\gamma\in\Gamma}$, $(Y_\delta)_{\delta\in\Delta}$,
and $(S_\epsilon)_{\epsilon\in\rE}$ are countable diagrams of
monomorphisms in $\sK$, then so are the diagrams $(X_\xi)_{\xi\in\Xi}$, $(Y_\xi)_{\xi\in\Xi}$, and $(S_\xi)_{\xi\in\Xi}$.

 Now let $\sK$ be a category with fibered coproducts.
 Let $\fX\rarrow\fS\larrow\fY$ be a cospan in $\Ind\sK$, represented
by a directed diagram of cospans $(X_\xi\rightarrow S_\xi\leftarrow
Y_\xi)$ in $\sK$, indexed by a directed poset~$\Xi$.
 Then one can easily see that the directed diagram
$(X_\xi\times_{S_\xi}Y_\xi)_{\xi\in\Xi}$ in $\sK$ represents
the fibered coproduct $\fX\times_\fS\fY$ in $\Ind\sK$.
 It suffices to compute morphisms in $\Ind\sK$ by
the formula~\eqref{morphisms-of-ind-objects} and use the fact that,
in the category of sets, the finite limits commute with both
the directed limits and the directed colimits.
 So $\Ind\sK$ is a category with fibered coproducts as well.

 Recall that the full subcategory of concentrated schemes $\CSch$ is
closed under fibered products in~$\Sch$ \,\cite[Remark Tag~0816]{SP}.
 Hence the category of concentrated schemes $\sK=\CSch$ has fibered
products.
 It follows that the category of (ind-quasi-compact ind-quasi-separated)
ind-schemes $\Ind(\CSch)$ is a category with fibered products, which
can be computed termwise on the cospans of directed diagrams in~$\CSch$.
 Furthermore, the full subcategories of strict ind-schemes and strict
$\aleph_0$\+ind-schemes are closed under fibered products in
$\Ind(\CSch)$.

\subsection{Ind-affine ind-schemes}
 Given a scheme $X$, we denote, as usual, by $\cO_X$ the structure
sheaf of~$X$.
 In particular, for an affine scheme $U$, the notation $R=\cO(U)=
\cO_U(U)$ stands for the commutative ring such that $U=\Spec R$.

 Following Corollary~\ref{Ind-of-Sch-strict-qs-qc-ss-lemma}(e) and
the subsequent paragraph, any strict ind-affine ind-scheme can be
represented by a directed diagram of closed immersions of affine
schemes.
 Furthermore, by Lemma~\ref{strict-aleph0-ind-schemes-lemma}(b), any
strict ind-affine $\aleph_0$\+ind-scheme can be represented by
a countable directed diagram of closed immersions of affine schemes.
 This was already mentioned at the end of
Section~\ref{basics-of-ind-schemes-subsecn}.

 What are the strict ind-affine ind-schemes?
 Let $\fU=\ilim_{\gamma\in\Gamma}U_\gamma$ be an ind-affine ind-scheme
represented by a directed diagram of closed immersions of affine
schemes.
 So, for every pair of indices $\gamma'\le\gamma''\in\Gamma$,
the transition map $\cO(U_{\gamma''})\rarrow\cO(U_{\gamma'})$ is
a surjective homomorphism of commutative rings.

 Let us say that the ind-affine ind-scheme $\fU$ is
\emph{limit-epimorphic} if the projection map
$\varprojlim_{\gamma\in\Gamma}\cO(U_\gamma)\rarrow\cO(U_\delta)$
is surjective for every index $\delta\in\Gamma$.
 One can easily see that this property does not depend on the choice
of a directed diagram of closed immersions of affine schemes
$(U_\gamma)_{\gamma\in\Gamma}$ representing the given strict
ind-affine ind-scheme~$\fU$.
 By the standard property of countable projective limits, any strict
ind-affine $\aleph_0$\+ind-scheme is limit-epimorphic.

\begin{lem} \label{limit-epimorphic-ind-affine-ind-schemes}
 The full subcategory of limit-epimorphic strict ind-affine ind-schemes
in\/ $\Ind(\CSch)$ is naturally anti-equivalent to the category
of complete, separated topological commutative rings where open ideals
form a base of neighborhoods of zero (with continuous ring
homomorphisms as the morphisms).
 In particular, the full subcategory of strict ind-affine\/
$\aleph_0$\+ind-schemes in\/ $\Ind(\CSch)$ is naturally
anti-equivalent to the category of complete, separated topological
rings having a countable base of neighborhoods of zero consisting of
open ideals.
\end{lem}

\begin{proof}
 This is~\cite[Section~7.11.2(i)]{BD2} or~\cite[Examples~1.6
and the subsequent paragraph]{Psemten}.
 Another relevant reference is~\cite[Theorem~9.4(a)]{Pextop}.

 To spell out some details, the desired equivalence of categories
assigns to a limit-epimorphic strict ind-affine ind-scheme $\fU=
\ilim_{\gamma\in\Gamma}U_\gamma$ as above the topological commutative
ring $\fR=\varprojlim_{\gamma\in\Gamma}\cO(U_\gamma)$ with
the topology of projective limit of discrete rings $\cO(U_\gamma)$.
 The inverse functor assigns to a complete, separated topological
commutative ring $\fR$ with a topology base of open ideals
the ind-affine ind-scheme $\fU=\ilim_{\fI\subset\fR}\Spec\fR/\fI$.
 Here $\fI$ ranges the poset of all open ideals $\fI\subset\fR$
(ordered by reverse inclusion), and $\Spec\fR/\fI$ is the spectrum
of the discrete commutative ring~$\fR/\fI$.
\end{proof}

 We denote the equivalence of categories from
Lemma~\ref{limit-epimorphic-ind-affine-ind-schemes} by
$\fR=\fO(\fU)$ and $\fU=\Spf\fR$ (where $\fU$ and $\fR$ are
the ind-affine ind-scheme and the topological ring corresponding
to each other under this category equivalence).

\subsection{Classes of morphisms of ind-schemes}
\label{classes-of-morphisms-of-ind-schemes-subsecn}
 Let $\bP$ be a property of morphisms in $\CSch$ that is preserved by
the pullbacks in~$\CSch$.
 We will say that a morphism $\ff\:\fY\rarrow\fX$ in $\Ind(\CSch)$ is
\emph{lax-ind\+$\bP$} if $\ff$~can be represented by a directed diagram of morphisms $(f_\gamma\:Y_\gamma\to X_\gamma)_{\gamma\in\Gamma}$
in $\CSch$, i.~e., $\ff=\ilim_{\gamma\in\Gamma}f_\gamma$, in such
a way that all the morphisms~$f_\gamma$ have property~$\bP$
\,\cite[Definition~1.19]{Rich}.

\begin{lem} \label{morphisms-with-property-lax-ind-P-and-strictness}
 Let\/ $\bP$ be a property of morphisms in\/ $\CSch$ that is
stable under pullbacks, and let\/ $\fX\in\Ind(\CSch)$ be a strict
ind-scheme.
 Then a morphism\/ $\ff\:\fY\rarrow\fX$ in\/ $\Ind(\CSch)$ is
lax-ind\+$\bP$ if and only if\/ $\ff$~can be represented by a directed
diagram of morphisms $(f_\xi\:Y_\xi\to X_\xi)_{\xi\in\Xi}$ in $\CSch$,
i.~e., $\ff=\ilim_{\xi\in\Xi}f_\xi$, in such a way that all
the morphisms~$f_\xi$ have property\/ $\bP$ \emph{and} all
the transition morphisms in the directed diagram $(X_\xi)_{\xi\in\Xi}$
are closed immersions.
\end{lem}

\begin{proof}
 Let $\ff=\ilim_{\gamma\in\Gamma} g_\gamma$ be a representation
of~$\ff$ by a directed diagram of morphisms
$(g_\gamma\:Y_\gamma\to X_\gamma)_{\gamma\in\Gamma}$ in $\CSch$ such
that all the morphisms~$g_\gamma$ have property $\bP$, and let
$\fX\simeq\ilim_{\delta\in\Delta}X'_\delta$ be a representation of $\fX$
by a directed diagram of closed immersions in~$\CSch$.
 Put $\fX'=\ilim_{\delta\in\Delta}X'_\delta$; so we have
an (iso)morphism of ind-schemes $\fX'\rarrow\fX$.
 Following the discussion in
Section~\ref{fibered-products-of-ind-schemes-subsecn} and
the explicit constructions of representing direct systems of diagrams
in~\cite[Section~9]{Pextop} and~\cite[Section~1.1]{Psemten}, one
can produce a directed poset $\Xi$ and a $\Xi$\+indexed direct system
of cospans $(Y_\xi\rightarrow X_\xi\leftarrow X'_\xi)$ in $\CSch$
such that the morphisms $g_\xi\:Y_\xi\rarrow X_\xi$ have property $\bP$,
the transition morphisms in the direct system $(X'_\xi)_{\xi\in\Xi}$
are closed immersions, $\ff=\ilim_{\xi\in\Xi}g_\xi$, and the induced
morphism of ind-schemes $\fX'=\ilim_{\xi\in\Xi}X'_\xi\rarrow
\ilim_{\xi\in\Xi}X_\xi=\fX$ is our given (iso)morphism.
 Now we have $\fY\simeq\fX'\times_\fX\fY=\ilim_{\xi\in\Xi}
(X'_\xi\times_{X_\xi}\nobreak Y_\xi)$ by the same discussion
in Section~\ref{fibered-products-of-ind-schemes-subsecn}.
 It remains to point out that the induced morphism
$f_\xi\:X'_\xi\times_{X_\xi}Y_\xi\rarrow X'_\xi$ has property $\bP$
as a pullback of the morphism~$g_\xi$.
\end{proof}

 More generally, it is clear from the discussion of the fibered
products in Section~\ref{fibered-products-of-ind-schemes-subsecn}
that the class of all lax-ind\+$\bP$ morphism of ind-schemes is
closed under pullbacks.
 In other words, if $\fY\rarrow\fX$ is a lax-ind\+$\bP$ morphism and
$\fX'\rarrow\fX$ is a morphism in $\Ind(\CSch)$, then the pullback
$\fX'\times_\fX\fY\rarrow\fX'$ is also a lax-ind\+$\bP$ morphism.
 An argument similar to the proof of
Lemma~\ref{morphisms-with-property-lax-ind-P-and-strictness}
proves as much.

\begin{rem} \label{ind-P-morphisms-of-strict-ind-schemes-remark}
 Lemma~\ref{morphisms-with-property-lax-ind-P-and-strictness} only
solves a half of the natural question about the definition of
an ind\+$\bP$ morphism of strict ind-schemes.
 Let $\fX$ and $\fY$ be strict ind-schemes, and $\ff\:\fY\rarrow\fX$
be a lax-ind\+$\bP$ morphism.
 Does it follow that that $\ff$~can be represented by a directed diagram
of morphisms $(f_\gamma\:Y_\gamma\to X_\gamma)_{\gamma\in\Gamma}$ in
$\CSch$ such that all the morphisms~$f_\gamma$ have property $\bP$ and
the transition morphisms in \emph{both} the directed diagrams
$(X_\gamma)_{\gamma\in\Gamma}$ and $(Y_\gamma)_{\gamma\in\Gamma}$ are
closed immersions of schemes?
 We do not know.

 For this reason, from now on and for the rest of this paper
\emph{we will generally presume all our ind-schemes to be strict}, and
all our directed diagrams in $\CSch$ representing ind-schemes to have
closed immersions as transition morphisms.

 In particular, modifying the definition above, we will say that
a morphism~$\ff$ of strict ind-concentrated ind-schemes is
\emph{ind\+$\bP$} if it can be represented by a directed diagram of
morphisms $(f_\gamma\:Y_\gamma\to X_\gamma)_{\gamma\in\Gamma}$ in
$\CSch$ such that all the morphisms~$f_\gamma$ have property $\bP$ and
the transition morphisms in both the directed diagrams
$(X_\gamma)_{\gamma\in\Gamma}$ and $(Y_\gamma)_{\gamma\in\Gamma}$
are closed immersions of schemes.
\end{rem}

 Let $\ff\:\fY\rarrow\fX$ be an ind\+$\bP$ morphism of strict
ind-concentrated ind-schemes and $\fg\:\fX'\rarrow\fX$ be any morphism
of strict ind-concentrated ind-schemes.
 Then the morphism of strict ind-concentrated ind-schemes
$\ff'\:\fX'\times_\fX\fY\rarrow\fX'$ is also ind\+$\bP$.
 Indeed, once again arguing as in the proof of
Lemma~\ref{morphisms-with-property-lax-ind-P-and-strictness},
it remains to point out that the transition morphisms of schemes
$X'_\xi\times_{X_\xi} Y_\xi\rarrow X'_\upsilon\times_{X_\upsilon}
Y_\upsilon$ are closed immersions for all $\xi\le\upsilon\in\Xi$,
since the transition morphisms $X'_\xi\rarrow X'_\upsilon$, \
$X_\xi\rarrow X_\upsilon$, and $Y_\xi\rarrow Y_\upsilon$ are closed
immersions.

 Let $\fX$ and $\fY\in\Ind(\CSch)$ be two (ind-concentrated)
ind-schemes.
 One says that a morphism of ind-schemes $\fY\rarrow\fX$ is
``representable by concentrated schemes''~\cite[Lemma~1.7 and
Definition~1.8]{Rich}, \cite[Section~1.3]{Psemten} if, for any
concentrated scheme $T\in\CSch$ and any morphism of ind-schemes
$T\rarrow\fX$, the fibered product $T\times_\fX\fY$ computed in
$\Ind(\CSch)$ is an object of the full subcategory
$\CSch\subset\Ind(\CSch)$ (i.~e., a concentrated scheme).

 Let $\fX=\ilim_{\gamma\in\Gamma}X_\gamma$ be an ind-scheme
represented by a directed diagram in $\CSch$, and let
$\ff\:\fY\rarrow\fX$ be a morphism of ind-schemes that is
``representable by (concentrated) schemes''.
 Put $Y_\gamma=X_\gamma\times_\fX\fY$, and let $f_\gamma\:
Y_\gamma\rarrow X_\gamma$ be the induced morphism of schemes.
 Then we have $\fY=\ilim_{\gamma\in\Gamma}Y_\gamma$ and
$\ff=\ilim_{\gamma\in\Gamma}f_\gamma$.
 Clearly, $Y_\gamma=X_\gamma\times_{X_\delta}Y_\delta$ for all
$\gamma\le\delta\in\Gamma$.

 Conversely, let $\ff\:\fY\rarrow\fX$ be a morphism of ind-schemes
represented by a directed diagram $(f_\gamma\:Y_\gamma\to
X_\gamma)_{\gamma\in\Gamma}$ of morphisms in $\CSch$ such that
one has $Y_\gamma=X_\gamma\times_{X_\delta}Y_\delta$ for all
$\gamma\le\delta\in\Gamma$.
 Then it follows that $Y_\gamma=X_\gamma\times_\fX\fY$.
 For any concentrated scheme $T$ and any morphism of ind-schemes
$T\rarrow\fX$, there exists an index $\gamma_0\in\Gamma$ such
that the morphism $T\rarrow\fX$ factorizes as $T\rarrow X_\gamma
\rarrow\fX$.
 Then $T\times_\fX\fY=T\times_{X_{\gamma_0}}(X_{\gamma_0}\times_\fX\fY)
=T\times_{X_{\gamma_0}}Y_{\gamma_0}$ is a concentrated scheme,
so the morphism~$\ff$ is ``representable by schemes''.

 A morphism of ind-schemes~$\ff$ is said to have property~$\bP$ (as
distinguished from ind\+$\bP$) if $\ff$~is simultaneously ind\+$\bP$
and ``representable by schemes''.

 It follows immediately from the definition that the class of all
morphisms of ind-schemes that are ``representable by (concentrated)
schemes'' is stable under pullbacks.
 Hence the class of all morphisms of ind-schemes having property $\bP$
is preserved by the pullbacks, too.

\begin{lem} \label{morphisms-with-property-P}
 Let\/ $\bP$ be a property of morphisms in\/ $\CSch$ that is stable
under pullbacks.
 Then a morphism of ind-schemes\/ $\ff\:\fY\rarrow\fX$ has property\/
$\bP$ if and only if, for any concentrated scheme $T$ and any morphism
of ind-schemes $T\rarrow\fX$, the fibered product $T\times_\fX\fY$ is
a concentrated scheme and the induced morphism of schemes
$T\times_\fX\fY\rarrow T$ has property\/~$\bP$.
\end{lem}

\begin{proof}
 The ``if'' implication follows from the representation of
the morphism~$\ff$ as $\ff=\ilim_{\gamma\in\Gamma}f_\gamma$ with
$f_\gamma\:Y_\gamma=X_\gamma\times_\fX\fY\rarrow X_\gamma$, in
the notation above.

 To prove the ``only if'', pick a directed diagram of morphisms of
concentrated schemes $(f_\gamma\:Y_\gamma\to
X_\gamma)_{\gamma\in\Gamma}$ such that all the morphisms~$f_\gamma$
have property $\bP$ and $\ff=\ilim_{\gamma\in\Gamma}f_\gamma$.
 Let $T\rarrow\fX$ be a morphism into $\fX$ from a concentrated
scheme~$T$.
 Then there exists an index $\gamma_0\in\Gamma$ such that
the morphism $T\rarrow\fX$ factorizes as $T\rarrow X_{\gamma_0}
\rarrow\fX$.
 Following the discussion in
Section~\ref{fibered-products-of-ind-schemes-subsecn}, we have
$T\times_\fX\fY=\ilim_{\gamma\ge\gamma_0}
(T\times_{X_\gamma}\nobreak Y_\gamma)$.
 Since the property $\bP$ is stable under pullbacks in $\CSch$,
the morphisms of schemes $T\times_{X_\gamma}Y_\gamma\rarrow T$ have
property~$\bP$.
 By assumption, the morphism of ind-schemes~$\ff$ is ``representable
by schemes'', so we have $T\times_\fX\fY\in\CSch\subset\Ind(\CSch)$.

 According to the convention of and the definition in
Remark~\ref{ind-P-morphisms-of-strict-ind-schemes-remark}, we assume
that the transitions morphisms in the directed diagrams
$(X_\gamma)_{\gamma\in\Gamma}$ and $(Y_\gamma)_{\gamma\in\Gamma}$
are closed immersions of schemes.
 Then all the transition morphisms in the directed diagram
$(T\times_{X_\gamma}\nobreak Y_\gamma)_{\gamma\in\Gamma}$ are
closed immersions of schemes, too.
 Since $T\times_\fX\fY\in\CSch$ and the closed immersions are
monomorphisms in $\CSch$, it follows that there exists an index
$\gamma\in\Gamma$ for which $T\times_\fX\fY=
T\times_{X_\gamma}Y_\gamma$.
 Thus the morphism of schemes $T\times_\fX\fY\rarrow T$ has
property $\bP$, as desired.
\end{proof}

 We have proved that a morphism of ind-schemes~$\ff$ has property $\bP$
if and only if it can be represented by a directed diagram
$(f_\gamma\:Y_\gamma\to X_\gamma)_{\gamma\in\Gamma}$ of morphisms in
$\CSch$ such that $Y_\gamma=X_\gamma\times_\fX\fY$ for all
$\gamma\in\Gamma$ (or equivalently, $Y_\gamma=X_\gamma\times_{X_\delta}
Y_\delta$ for all $\gamma\le\delta\in\Gamma$) and all the
morphisms~$f_\gamma$ have property~$\bP$.

 In particular, taking $\bP$ to be the property of a morphism of
schemes to be affine, we obtain the notions of an \emph{affine}
and \emph{ind-affine} morphism of ind-schemes, as well as the notion
of a \emph{lax-ind-affine} morphism of ind-schemes.
 Every morphism of ind-affine ind-schemes is ind-affine, but it need
not be affine.

 Conversely, if $\ff:\fV\rarrow\fU$ is a lax-ind-affine morphism of
ind-schemes and the ind-scheme $\fU$ is ind-affine, then
the ind-scheme $\fV$ is ind-affine, too.
 Indeed, by Lemma~\ref{morphisms-with-property-lax-ind-P-and-strictness}
we have $\ff=\ilim_{\gamma\in\Gamma}f_\gamma$, where
$(f_\gamma\:V_\gamma\rarrow U_\gamma)_{\gamma\in\Gamma}$ is a directed
diagram of affine morphisms of concentrated schemes and the transition
morphisms in the direct directed diagram of schemes
$(U_\gamma)_{\gamma\in\Gamma}$ are closed immersions.
 By Corollary~\ref{Ind-of-Sch-strict-qs-qc-ss-lemma}(e), all
the schemes $U_\gamma$ are affine, and it follows that all
the schemes $V_\gamma$ are affine, too.

 Taking $\bP$ to be the property of a morphism of schemes to be
a closed immersion, we obtain the notions of a \emph{closed immersion}
and an \emph{ind-closed immersion} of ind-schemes.
 In particular, a morphism $Z\rarrow\fX$ from a concentrated scheme
to an (ind-concentrated) ind-scheme is a closed immersion if and only
if there exist a representation of $\fX$ by a directed diagram of
(closed immersions of) schemes $(X_\gamma)_{\gamma\in\Gamma}$,
an index $\gamma_0\in\Gamma$, and an isomorphism of schemes
$Z\simeq X_{\gamma_0}$ such that the triangular diagram
$Z\simeq X_{\gamma_0}\rarrow\fX$ is commutative.

 A \emph{closed subscheme} $Z\subset\fX$ of an ind-scheme $\fX$ is
an equivalence class of pairs (a concentrated scheme $Z$, a closed
immersion $Z\rarrow\fX$) with respect to the obvious equivalence
relation.
 Specifically, two closed immersions $Z'\rarrow\fX$ and $Z''\rarrow\fX$
are said to be equivalent if there exists an isomorphism $Z'\simeq Z''$
making the triangular diagram $Z'\simeq Z''\rarrow\fX$ commutative.

 Given two closed subschemes $Z'$ and $Z''$ in $\fX$, one writes
$Z'\subset Z''$ if the closed immersion morphism $Z'\rarrow\fX$
factorizes through the closed immersion morphism $Z''\rarrow\fX$.
 If this is the case, then the morphism of schemes $Z'\rarrow Z''$
is a closed immersion.

\begin{ex} \label{ind-closed-not-closed-immersion-example}
 Let $k$~be a field and $x$~be a variable.
 Then the inclusion $i\:k[x]\rarrow k[[x]]$ of the discrete ring of
polynomials $k[x]$ into the topological ring of formal power series
$k[[x]]$ (with the adic topology on $k[[x]]$) induces, by
the anti-equivalence of categories from
Lemma~\ref{limit-epimorphic-ind-affine-ind-schemes}, a morphism
of ind-affine ind-schemes $\fii:\Spf k[[x]]\rarrow\Spec k[x]$.
 The morphism of ind-affine ind-schemes~$\fii$ is an ind-closed
immersion, but \emph{not} a closed immersion.
 In fact, if an ind-scheme $\fX$ admits a closed immersion (or more
generally, a ``morphism representable by schemes'') into a concentrated
scheme, then $\fX$ is a concentrated scheme.
\end{ex}

\subsection{Open ind-subschemes}
 In the context of
Section~\ref{classes-of-morphisms-of-ind-schemes-subsecn},
take $\bP$ to be the property of a morphism of schemes to be
an open immersion.
 This produces the notions of an \emph{open immersion} and
an \emph{ind-open immersion} of ind-schemes.

\begin{ex} \label{ind-open-but-not-open-immersion-counterex}
 The following counterexample shows that an ind-open immersion
of (ind-concentrated) ind-schemes need not be an open immersion.
 Let $\Xi$ be an infinite set and $(U_\xi=\Spec R_\xi)_{\xi\in\Xi}$
be a family of nonempty affine schemes indexed by $\Xi$, as in
Remark~\ref{ind-objects-in-schemes-problematic-remark}.
 Following the notation of that remark, put
$X=\coprod_{\xi\in\Xi}U_\xi$ and $\fX=\ilim_{\gamma\in\Gamma}
\coprod_{\xi\in\gamma}U_\xi$, where $\Gamma$ is the directed poset
of all finite subsets of~$\Xi$.
 Then the non-quasi-compact scheme $X$ represents the ind-concentrated
ind-scheme~$\fX$.

 Consider the infinite product ring $S=\prod_{\xi\in\Xi}R_\xi$ and
the affine scheme $V=\Spec S$.
 Then, for every index $\xi\in\Xi$, the affine scheme $U_\xi$ is
an open-closed subscheme in the affine scheme~$V$; so $V$ is
the disjoint/disconnected union of the affine schemes $U_\xi$ and
$V\setminus U_\xi$.
 It follows that $X$ is an open subscheme in~$V$.
 The natural morphism of ind-schemes $\fX\rarrow V$ is an ind-open
immersion, but it is not an open immersion in our terminology, as it is
not ``representable by (concentrated) schemes''.
\end{ex}

 An \emph{open ind-subscheme} $\fU\subset\fX$ of an ind-scheme $\fX$ is
an equivalence class of pairs (an ind-scheme $\fU$, an open
immersion $\fU\rarrow\fX$) with respect to the obvious equivalence
relation.
 Specifically, two open immersions $\fU'\rarrow\fX$ and
$\fU''\rarrow\fX$ are said to be equivalent if there exists
an isomorphism $\fU'\simeq\fU''$ making the triangular diagram
$\fU'\simeq\fU''\rarrow\fX$ commutative.
 (Cf.~\cite[Definition~1.13]{Rich}.)

 The inclusion order on the set of all open ind-schemes $\fU\subset\fX$
is also defined in the obvious way.
 One has $\fU'\subset\fU''$ if the open immersion $\fU'\rarrow\fX$
factorizes through the open immersion $\fU''\rarrow\fX$.
 In this case, the morphism $\fU'\rarrow\fU''$ is an open immersion
of ind-schemes as well.

 For any scheme $X$, one can consider its underlying topological space
(which we denote also by~$X$).
 Given a (strict ind-concentrated) ind-scheme $\fX$, its underlying
topological space $\fX$ is constructed as follows.

 Let $(X_\gamma)_{\gamma\in\Gamma}$ be a directed diagram of
(concentrated) schemes and closed immersions of schemes
representing~$\fX$.
 Then the topological space $\fX$ is the inductive limit (the union)
of the topological spaces $X_\gamma$, endowed with the inductive
limit topology.
 So a subset $\fV\subset\fX$ is said to be open if $\fV\cap X_\gamma$
is an open subset of $X_\gamma$ for every $\gamma\in\Gamma$.
 Clearly, this definition of the underlying topological space of $\fX$
does not depend on the choice of a directed diagram of closed
immersions of schemes representing~$\fX$.

 Let us say that an open subset $\fV\subset\fX$ is
\emph{ind-quasi-compact} if the open subset $\fV\cap X_\gamma
\subset X_\gamma$ is quasi-compact for every $\gamma\in\Gamma$.

\begin{lem} \label{open-ind-subschemes=qc-open-subsets}
 For any ind-scheme\/ $\fX$, there is a natural order-preserving
bijection between the poset of all open ind-subschemes in\/ $\fX$
and the poset of all ind-quasi-compact open subsets in the underlying
topological space of\/ $\fX$, both ordered by inclusion.
 To any open ind-subscheme\/ $\fj\:\fU\rarrow\fX$, the open subset\/
$\fj(\fU)\subset\fX$ of the underlying topological space of\/ $\fX$
is assigned.
\end{lem}

\begin{proof}
 The proof is left to the reader.
\end{proof}

\begin{rem}
 The underlying topological spaces of affine schemes are quasi-compact,
but the underlying topological spaces of ind-affine ind-schemes are
\emph{almost never} quasi-compact.
 A general discussion can be found in~\cite[Section~1.4]{Rich}.
 Let us present one explicit counterexample.

 To keep things uncomplicated, let $k$~be a field of cardinality
at least~$3$.
 Let $R_n=k[x_1,\dotsc,x_n]$, \,$n\ge0$, be the rings of polynomials
in $n$~variables over~$k$.
 Then we have the natural isomorphisms of $k$\+algebras
$R_{n+1}/(x_{n+1})\simeq R_n$ taking $x_i$ to~$x_i$ for all
$1\le i\le n$.
 This produces a projective system of surjective homomorphisms
of commutative rings $R_0\larrow R_1\larrow R_2\larrow\dotsb$.
 Put $X_n=\Spec R_n$, and let $X_0\rarrow X_1\rarrow X_2\rarrow\dotsb$
be the corresponding inductive system of closed immersions of
affine schemes.
 Then the underlying topological space of the strict ind-affine
$\aleph_0$\+ind-scheme $\fX=\ilim_{n\ge0}X_n$ is \emph{not}
quasi-compact.
 
 Indeed, let $a$, $b\in k$ be any two different nonzero elements,
$a\ne b$, \,$a\ne0\ne b$.
 Consider the Cartesian power set $F=\{a,b\}^{\boZ_{\ge1}}$.
 The elements of $F$ are arbitrary functions $\phi\:\boZ_{\ge1}
\rarrow\{a,b\}$.
 Given such a function~$\phi$, put
$R_{n,\phi}=R_n[(x_1-\phi(1))^{-1},\dotsc,(x_n-\phi(n))^{-1}]$.
 Here $x_i-\phi(i)\in k[x_i]\subset R_n$ (for $1\le i\le n$) is
a polynomial of degree~$1$.
 So $X_{n,\phi}=\Spec R_{n,\phi}$ is a principal affine open subscheme
of the affine scheme~$X_n$.

 Similarly to the above, we have natural isomorphism of $k$\+algebras
$R_{n+1,\phi}/(x_{n+1})\simeq R_{n,\phi}$ taking $x_i$ to~$x_i$ for all
$1\le i\le n$, and the related projective system of surjective
homomorphisms of commutative rings $R_{0,\phi}\larrow R_{1,\phi}
\larrow R_{2,\phi}\larrow\dotsb$.
 Let $X_{0,\phi}\rarrow X_{1,\phi}\rarrow X_{2,\phi}\rarrow\dotsb$
be the corresponding inductive system of closed immersions of
affine schemes.
 Then $\fX_\psi=\ilim_{n\ge0}X_{n,\phi}$ is an open ind-subscheme
in~$\fX$.

 Let $C\subset F$ be any subset such that, for any function
$\phi\in F$ and any integer $n\ge1$, there exists a function
$\psi\in C$ for which $\phi(i)=\psi(i)$ for all $1\le i\le n$.
 One can easily make the set $C$ countable, but it cannot be made
finite.
 Now $\fX=\bigcup_{\psi\in C}\fX_\psi$ is an infinite open covering
of $\fX$ that has no finite subcovering.
 Indeed, one has $X_n=\bigcup_{\psi\in C} X_{n,\psi}$ for every
$n\ge0$, and this open covering of the affine scheme $X_n$ has
\emph{no} subcovering consisting of fewer than $2^n$ open subsets.

 Another counterexample of non-quasi-compact strict ind-affine
$\aleph_0$\+ind-scheme is provided by the ind-scheme $\fX$ from
Remark~\ref{ind-objects-in-schemes-problematic-remark} and
Example~\ref{ind-open-but-not-open-immersion-counterex}.
\end{rem}

\subsection{Ind-affine-finite ind-schemes}
 Let $\fX$ be a strict ind-concentrated ind-scheme.
 We will say that the ind-scheme $\fX$ is \emph{ind-affine-finite}
if there exists a finite collection of open ind-subschemes
$\fU_1$,~\dots, $\fU_N\subset\fX$ such that the underlying
topological space of $\fX$ is the union of the underlying topological
spaces of $\fU_\alpha$, that is $\fX=\bigcup_{\alpha=1}^N\fU_\alpha$.

 The latter condition can be also expressed in the language of
schemes and ind-schemes as follows.
 In the context of
Section~\ref{classes-of-morphisms-of-ind-schemes-subsecn}, consider
the property $\bP$ of a morphism of schemes to be surjective.
 The class of all surjective morphisms of schemes is closed under
pullbacks by~\cite[Lemma Tag~01S1]{SP}.
 This produces the notions of a \emph{surjective} and
an \emph{ind-surjective} morphism of ind-schemes.

 Given a finite collection of open ind-subschemes
$\fU_1$,~\dots, $\fU_N\subset\fX$, consider the disjoint/disconnected
union of ind-schemes $\fU=\coprod_{\alpha=1}^N\fU_\alpha$.
 Then we have a natural morphism of ind-schemes $\fU\rarrow\fX$,
which is obviously ``representable by (concentrated) schemes''
(since so are the open immersions $\fU_\alpha\rarrow\fX$).
 One has $\fX=\bigcup_{\alpha=1}^N\fU_\alpha$ if and only if
the morphism of ind-schemes $\fU\rarrow\fX$ is surjective, or
equivalently, ind-surjective.

\begin{ex} \label{ind-affine-over-ind-affine-finite-example}
 Let $\fX$ be an ind-affine-finite ind-scheme and $\ff\:\fY\rarrow\fX$
be a lax-ind-affine morphism of ind-schemes.
 (In particular, any ind-affine morphism of ind-schemes is
lax-ind-affine.)
 Then the ind-scheme $\fY$ is also ind-affine-finite.

 Specifically, let $\fX=\bigcup_{\alpha=1}^N\fU_\alpha$ be a finite
ind-affine open covering of~$\fX$.
 Put $\fV_\alpha=\fU_\alpha\times_\fX\fY$.
 Following the discussion of the preservation by pullbacks in
Section~\ref{classes-of-morphisms-of-ind-schemes-subsecn},
the morphisms of ind-schemes $\fV_\alpha\rarrow\fY$ are open immersions,
while the morphisms $\fV_\alpha\rarrow\fU_\alpha$ are lax-ind-affine.
 According to the discussion of lax-ind-affine morphisms
in Section~\ref{classes-of-morphisms-of-ind-schemes-subsecn}, it
follows that the ind-schemes $\fV_\alpha$ are ind-affine.
 Finally, one has $\fY=\bigcup_{\alpha=1}^N\fV_\alpha$ (e.~g.,
since the classes of surjective and ind-surjective morphisms of
ind-schemes are closed under pullbacks according to
Section~\ref{classes-of-morphisms-of-ind-schemes-subsecn}).
\end{ex}

 For any two categories $\sK$ and $\sL$, any functor
$F\:\sK\rarrow\sL$ induces a functor on the categories of ind-objects
$\Ind F\:\Ind\sK\rarrow\Ind\sL$ defined by the obvious rule
$(\Ind F)(\ilim_{\gamma\in\Gamma}X_\gamma)=
\ilim_{\gamma\in\Gamma}F(X_\gamma)$.
 In particular, for any concentrated scheme $X$, the maximal reduced
closed subscheme $X_\red\subset X$ is a concentrated scheme
depending functorially on~$X$.
 Using the construction above, the functor $X\longmapsto X_\red\:
\CSch\rarrow\CSch$ induces a functor $\Ind(\CSch)\rarrow\Ind(\CSch)$,
which we also denote by $\fX\longmapsto\fX_\red$.

 Clearly, the functor $X\longmapsto X_\red$ takes closed immersions
of schemes to closed immersions of schemes; hence, for any strict
ind-scheme $\fX$, the ind-scheme $\fX_\red$ is also strict.
 Similarly, the functor $\fX\longmapsto\fX_\red$ preserves
the ind-affineness, ind-semi-separatedness, and ind-Noetherianity
of ind-schemes, and takes $\aleph_0$\+ind-schemes to
$\aleph_0$\+ind-schemes.

 In the context of the theory of ind-schemes, a (strict
ind-concentrated) ind-scheme $\fX$ is said to be a \emph{formal
scheme} if the ind-scheme $\fX_\red$ is a concentrated scheme,
i.~e., $\fX_\red\in\CSch\subset\Ind(\CSch)$.
 Given a representation $\fX=\ilim_{\gamma\in\Gamma}X_\gamma$ of
$\fX$ by a directed diagram of closed immersions of concentrated
schemes $(X_\gamma)_{\gamma\in\Gamma}$, the ind-scheme $\fX$ is
a formal scheme if and only if there exists an index
$\gamma_0\in\Gamma$ for which the morphism $(X_{\gamma_0})_\red
\rarrow\fX_\red$ is an isomorphism.
 Equivalently, the morphisms $X_{\gamma_0})_\red\rarrow(X_\gamma)_\red$
are isomorphisms for all $\gamma\ge\gamma_0$.
 So an ind-scheme is a formal scheme if and only if it can be
represented by a directed diagram of \emph{bijective} closed
immersions of concentrated schemes.

 A \emph{Noetherian formal scheme} is an ind-scheme that is
simultaneously an ind-Noetherian ind-scheme and a formal scheme.
 Similarly, an \emph{affine formal scheme} is an ind-scheme that is
simultaneously an ind-affine ind-scheme and a formal scheme.
 The definition of a \emph{semi-separated} formal scheme is also
similar.
 So a formal scheme is Noetherian if and only if it can be represented
by a directed diagram of bijective closed immersions of Noetherian
schemes, and similarly for affine/semi-separated formal schemes
(by Corollary~\ref{Ind-of-Sch-strict-qs-qc-ss-lemma}(b,c,e) and
the subsequent paragraph).

\begin{ex} \label{Noetherian-formal-schemes-are-ind-affine-finite-ex}
 Any Noetherian formal scheme $\fX$ is ind-affine-finite.
 Indeed, let $(X_\gamma)_{\gamma\in\Gamma}$ be a directed diagram of
bijective closed immersions of Noetherian schemes representing~$\fX$.
 So we have $(X_\gamma)_\red=X=\fX_\red$ for all $\gamma\in\Gamma$.

 Let $X=\bigcup_{\alpha=1}^N U_\alpha$ be a finite affine open
covering of the Noetherian scheme~$X$.
 For every index $\gamma\in\Gamma$, the underlying topological space
of $X_\gamma$ is naturally identified with the underlying
topological space of~$X$; and $U_\alpha$ is an open subset of
this topological space.
 So we have a uniquely defined open subscheme $U_{\alpha,\gamma}\subset
X_\gamma$ corresponding to this open subset.
 In other words, $(U_{\alpha,\gamma})_\red=U_\alpha$.
 The scheme $U_{\alpha,\gamma}$ is Noetherian, and
by~\cite[Lemma Tag~01YQ]{SP} it follows that $U_{\alpha,\gamma}$
is affine.

 It remains to put
$\fU_\alpha=\ilim_{\gamma\in\Gamma}U_{\alpha,\gamma}$.
 Then $\fU_\alpha$ is an ind-affine open ind-subscheme in $\fX$,
and $\fX=\bigcup_{\alpha=1}^N\fU_\alpha$.
\end{ex}

\begin{rem} \label{ind-qc-need-not-be-ind-affine-finite}
 A strict ind-quasi-compact ind-(quasi)-separated $\aleph_0$\+ind-scheme
\emph{need not} be ind-affine-finite.
 Indeed, let $\fX=\bigcup_{\alpha=1}^N\fU_\alpha$ be
an ind-affine-finite strict ind-concentrated ind-scheme, covered by $N$
ind-affine open ind-subschemes~$\fU_\alpha$.
 Let $Z\subset\fX$ be a closed subscheme.
 Then the schemes $Z\times_\fX\fU_\alpha$ are affine as closed
subschemes in the ind-affine ind-schemes~$\fU_\alpha$.
 So $Z=\bigcup_\alpha(Z\times_\fX\fU_\alpha)$ is an affine open covering
of $Z$ by $N$ affine open subschemes.
 Thus the minimal number of covering open affines for closed subschemes
$Z\subset\fX$ is uniformly bounded by one constant $N$ common for all
the closed subschemes $Z\subset\fX$.

 For a specific counterexample, let $k$~be a field and $V=k^{(\omega)}$
be a $k$\+vector space of countably infinite dimension.
 Then the projectization $\fP(V)$ is a strict ind-quasi-compact
ind-separated $\aleph_0$\+ind-scheme such that, for every field
extension $K/k$, the $K$\+points of $\fP(V)$ correspond bijectively to
one-dimensional vector subspaces in the $K$\+vector space $K\ot_kV$.
 Simply put, $\fP(V)$ is the ind-scheme represented by the sequence of
standard coordinate closed immersions of projective spaces of increasing
dimension $\Spec k=\boP^0_k\rarrow\boP^1_k\rarrow\boP^2_k\rarrow\dotsb$.
 We refer to~\cite[Example~1.3(2)]{Rich} or~\cite[Section~11.2]{Psemten}
for a more detailed and formal definition of the ind-scheme~$\fP(V)$.

 Now the $n$\+dimensional projective space $\boP^n_k$ over~$k$ cannot
be covered by less than $n+1$ open affines.
 Indeed, for any quasi-compact semi-separated scheme $Z$ covered by
$n$~open affines and every quasi-coherent sheaf $\cM$ on $Z$ one has
$H^i(Z,\cM)=0$ for all $i\ge n$, since the sheaf cohomology
$H^*(Z,\cM)$ can be computed by the finite \v Cech complex~\cite[Lemmas
Tags~01XD and~01FM]{SP}, which vanishes in the cohomological
degrees~$\ge n$.
 On the other hand, one has $H^n(\boP^n_k,\cO(-n-1))=k\ne0$
\,\cite[Lemma Tag~01XT]{SP}.
 For every integer $n\ge0$, the scheme $\boP^n_k$ is a closed
subscheme of $\fP(V)$.
 Hence the ind-scheme $\fP(V)$ is \emph{not} ind-affine-finite.
\end{rem}

\subsection{Reasonable ind-schemes}
\label{reasonable-ind-schemes-subsecn}
 Let $X$ be a concentrated scheme.
 A quasi-coherent sheaf $\cM$ on $X$ is said to be (\emph{locally})
\emph{finitely generated} if, for every affine open subscheme
$U\subset X$, the $\cO_X(U)$\+module $\cM(U)$ is finitely generated.
 It suffices to check this condition for affine open subschemes
$U\subset X$ belonging to any chosen affine open covering of~$X$
(see, e.~g., \cite[Lemma~3.1(a)]{PS6}).

 Simply put, a quasi-coherent sheaf $\cM$ on $X$ is finitely
generated if and only if, as a sheaf of $\cO_X$\+modules, it is
generated by a finite set of local sections.

 Let $i\:Z\rarrow X$ be a closed immersion of concentrated schemes,
and let $\cI\subset\fO_X$ be the corresponding quasi-coherent sheaf
of ideals on~$X$.
 Following the terminology of~\cite[Sections~7.11.1 and~7.11.4]{BD2}
and~\cite[Section~2.1]{Psemten}, we say that the closed immersion
morphism~$i$ is \emph{reasonable} if the quasi-coherent sheaf $\cI$
on $X$ is finitely generated.
 If this is the case, $Z$ is also said to be a \emph{reasonable closed
subscheme} in~$X$.

 Let $\fX$ be a (strict ind-concentrated) ind-scheme.
 A closed subscheme $Z\subset\fX$ is said to be \emph{reasonable} if,
for every closed subscheme $Y\subset\fX$ such that $Z\subset Y$,
the closed subscheme $Z$ is reasonable in~$Y$.
 Given a representation of $\fX$ by a directed diagram of closed
immersions of concentrated schemes $(X_\gamma)_{\gamma\in\Gamma}$,
a closed subscheme $Z\subset\fX$ is reasonable if and only if,
for every index $\gamma\in\Gamma$ such that $Z\subset X_\gamma$,
the closed subscheme $Z$ is reasonable in~$X_\gamma$
\,\cite[Section~2.1]{Psemten}.

 An ind-scheme $\fX$ is said to be \emph{reasonable} if it is
a filtered inductive limit of its reasonable closed subschemes.
 Equivalently, $\fX$ is reasonable if and only if it can be
represented by a directed diagram $(X_\gamma)_{\gamma\in\Gamma}$
of reasonable closed immersions of concentrated schemes~$X_\gamma$.

\begin{ex} \label{ind-noetherian-are-reasonable-example}
 On a Noetherian scheme $X$, any quasi-coherent sheaf of ideals is
locally finitely generated (and in fact, coherent).
 Therefore, any strict ind-Noetherian ind-scheme is reasonable.
 Moreover, any closed subscheme of an ind-Noetherian ind-scheme
is reasonable~\cite[Section~4.2]{Psemten}.
\end{ex}

\begin{ex} \label{repr-by-schemes-over-reasonable-example}
 Let $\fY\rarrow\fX$ be a morphism of ind-schemes that is
``representable by (concentrated) schemes''.
 Assume that the ind-scheme $\fX$ is reasonable.
 Then the ind-scheme $\fY$ is reasonable as
well~\cite[Section~2.1]{Psemten}.
\end{ex}

\Section{Flat and Very Flat Pro-Quasi-Coherent Pro-Sheaves}
\label{flat-and-vfl-proqcoh-prosheaves-secn}

 Following the conventions of the end of
Section~\ref{basics-of-ind-schemes-subsecn} and
Remark~\ref{ind-P-morphisms-of-strict-ind-schemes-remark}, for
the rest of this paper an \emph{ind-scheme} means an object of
$\Ind(\CSch)$ represented by a directed diagram of closed immersions
of concentrated schemes.

 According to the discussion in~\cite[Sections~7.11.3\+-4]{BD2},
there are two kinds of quasi-coherent sheaf-type objects over
an ind-scheme.
 In the terminology of~\cite[Chapters~2\+-3]{Psemten}, they are
called the \emph{quasi-coherent torsion sheaves} and
the \emph{pro-quasi-coherent pro-sheaves}.
 In this paper, we only deal with pro-quasi-coherent pro-sheaves.

\subsection{Flat pro-quasi-coherent pro-sheaves}
\label{flat-proqcoh-prosh-subsecn}
 Following the discussion in~\cite[Section~7.11.3]{BD2}
and~\cite[Examples~3.1 and~3.8(1)]{Psemten}, the category of arbitrary
pro-quasi-coherent pro-sheaves on an ind-scheme is usually \emph{not}
abelian and generally \emph{not} well-behaved homologically.
 The category of \emph{flat} pro-quasi-coherent pro-sheaves is
a well-behaved exact category (in the sense of Quillen~\cite{Bueh}).
 For this reason, we restrict ourselves to flat pro-quasi-coherent
pro-sheaves in this paper.
 We suggest the survey paper~\cite{Bueh} as the reference source
on exact categories.

 Given a scheme $X$, we denote by $X\Qcoh$ the Grothendieck abelian
category of quasi-coherent sheaves on~$X$.
 A quasi-coherent sheaf $\cF$ on $X$ is said to be \emph{flat} if,
for every affine open subscheme $U\subset X$, the $\cO_X(U)$\+module
$\cF(U)$ is flat.
 It suffices to check this condition for affine open subschemes
$U\subset X$ belonging to any chosen affine open covering of~$X$.
 The notation $X\Qcoh_\fl\subset X\Qcoh$ stands for the full
subcategory of flat quasi-coherent sheaves in $X\Qcoh$.
 The full subcategory $X\Qcoh_\fl$ is closed under extensions, kernels
of epimorphisms, infinite coproducts, and directed colimits in
$X\Qcoh$; so $X\Qcoh_\fl$ inherits an exact category structure from
the abelian exact structure of $X\Qcoh$.

 For every morphism of concentrated schemes $f\:Y\rarrow X$, we have
a pair of adjoint functors of \emph{direct image} $f_*\:Y\Qcoh\rarrow
X\Qcoh$ and \emph{inverse image} $f^*\:X\Qcoh\rarrow Y\Qcoh$.
 The inverse image functor~$f^*$ is left adjoint to the direct image
functor~$f_*$.
 The inverse image functor~$f^*$ always takes flat quasi-coherent
sheaves on $X$ to flat quasi-coherent sheaves on~$Y$.
 The resulting functor $f^*\:X\Qcoh_\fl\rarrow Y\Qcoh_\fl$ is exact
(as a functor of exact categories).

 Let $\fX$ be a (strict ind-concentrated) ind-scheme.
 A \emph{flat pro-quasi-coherent pro-sheaf} $\rF$ on $\fX$ is
the following set of data:
\begin{enumerate}
\renewcommand{\theenumi}{\roman{enumi}}
\item to every closed subscheme $Z\subset\fX$, a flat quasi-coherent
sheaf $\rF^{(Z)}$ on $Z$ is assigned;
\item to every pair of closed subschemes $Z\subset Y$ in $\fX$ with
the closed immersion morphism $i_{ZY}\:Z\rarrow Y$, an isomorphism
of flat quasi-coherent sheaves $i_{ZY}^*\rF^{(Y)}\simeq\rF^{(Z)}$
on $Z$ is assigned.
\end{enumerate}
 The following compatibility axiom for triples of closed subschemes
is imposed:
\begin{enumerate}
\renewcommand{\theenumi}{\roman{enumi}}
\setcounter{enumi}{2}
\item for every triple of closed subschemes $W\subset Z\subset Y$
in $\fX$, the triangular diagram $i_{WY}^*\rF^{(Y)}\rarrow
i_{WZ}^*\rF^{(Z)}\rarrow\rF^{(W)}$ of isomorphisms in $W\Qcoh_\fl$
is commutative.
\end{enumerate}

 Let $\fX=\ilim_{\gamma\in\Gamma}X_\gamma$ be a representation of $\fX$
by a directed diagram of closed immersions of concentrated schemes.
 Then, in order to define a flat pro-quasi-coherent pro-sheaf $\rF$
on $\fX$, it suffices to specify flat pro-quasi-coherent pro-sheaves
$\rF^{(\gamma)}=\rF^{(X_\gamma)}$ on the closed subschemes
$X_\gamma\subset\fX$ for all indices $\gamma\in\Gamma$, together with
isomorphisms $i^*_{\gamma\delta}\rF^{(\delta)}\simeq\rF^{(\gamma)}$
for all pair of indices $\gamma\le\delta\in\Gamma$.
 Here the notation is $i_{\gamma\delta}=i_{X_\gamma X_\delta}\:
X_\gamma\rarrow X_\delta$.
 Assuming that the compatibility axiom~(iii) is satisfied for all
triples of closed subschemes $W=X_\beta$, \ $Z=X_\gamma$, and
$Y=X_\delta$ with $\beta\le\gamma\le\delta\in\Gamma$, such a datum
of quasi-coherent sheaves $\rF^{(\gamma)}$ and related isomorphisms
$i^*_{\gamma\delta}\rF^{(\delta)}\simeq\rF^{(\gamma)}$ can be extended
to a flat pro-quasi-coherent pro-sheaf $\rF$ on $\fX$ in a unique way.
 The flat quasi-coherent sheaves $\rF^{(Z)}$ for all other closed
subschemes $Z\subset\fX$ together with the related isomorphisms~(ii)
can be uniquely recovered in such a way that the compatibility
axoim~(iii) is satisfied for all $W\subset Z\subset Y\subset\fX$.

 Morphisms of flat pro-quasi-coherent pro-sheaves on $\fX$ are defined
in the obvious way.
 We denote the additive category of flat pro-quasi-coherent pro-sheaves
on an ind-scheme $\fX$ by $\fX\Proqc_\fl$.

 A short sequence of flat pro-quasi-coherent pro-sheaves $0\rarrow\rF
\rarrow\rG\rarrow\rH\rarrow0$ on $\fX$ is said to be (\emph{admissible})
\emph{exact} if the short sequence of flat quasi-coherent sheaves
$0\rarrow\rF^{(Z)}\rarrow\rG^{(Z)}\rarrow\rH^{(Z)}\rarrow0$ is exact
in $Z\Qcoh_\fl$ for every closed subscheme $Z\subset\fX$.
 Given a directed diagram of closed immersions of concentrated schemes
$(X_\gamma)_{\gamma\in\Gamma}$ representing the ind-scheme $\fX$, it
suffices to check the latter exactness condition for the closed
subschemes $Z=X_\gamma\subset\fX$, \,$\gamma\in\Gamma$.

 The class of short exact sequences defined in the previous paragraph
defines an exact category structure on the additive category
$\fX\Proqc_\fl$ of flat pro-quasi-coherent pro-sheaves on~$\fX$
\,\cite[Proposition~3.5]{Psemten}.
 The fact that the inverse image functors $i_{ZY}^*\:Y\Qcoh_\fl\rarrow
X\Qcoh_\fl$ are exact plays a crucial role here.

 All set-indexed coproducts and directed colimits exist in the additive
category $\fX\Proqc_\fl$.
 Given a infinite family $(\rF_\xi)_{\xi\in\Xi}$ of flat
pro-quasi-coherent pro-sheaves on $\fX$, indexed by some set $\Xi$,
the coproduct $\rG=\coprod_{\xi\in\Xi}\rF_\xi$ is constructed by
the obvious rule $\rG^{(Z)}=\coprod_{\xi\in\Xi}\rF_\xi^{(Z)}$ for all
closed subschemes $Z\subset\fX$.
 Similarly, given a directed diagram
$(\rF_\upsilon)_{\upsilon\in\Upsilon}$ of flat pro-quasi-coherent
pro-sheaves on $\fX$, indexed by some directed poset $\Upsilon$,
the directed colimit $\rG=\varinjlim_{\upsilon\in\Upsilon}\rF_\upsilon$
is constructed by the obvious rule
$\rG^{(Z)}=\varinjlim_{\upsilon\in\Upsilon}\rF_\upsilon^{(Z)}$ for
all closed subschemes $Z\subset\fX$ \,\cite[Section~3.5]{Psemten}.

 In other words, the functors $\rF\longmapsto\rF^{(Z)}\:
\fX\Proqc_\fl\rarrow Z\Qcoh_\fl$ preserve coproducts and directed
colimits for all closed subschemes $Z\subset\fX$.
 The functors of directed colimit are exact in the exact
category $\fX\Proqc_\fl$; in particular, this applies to
the infinite coproduct functors.

\begin{ex} \label{ind-affine-ind-schemes-and-flat-contramodules}
 Let $\fU=\ilim_{\gamma\in\Gamma}U_\gamma$ be a strict ind-affine
$\aleph_0$\+ind-scheme represented by a countable directed diagram
of closed immersions of affine schemes.
 Let $\fR=\varprojlim_{\gamma\in\Gamma}\cO(U_\gamma)$ be
the corresponding complete, separated topological commutative ring
with a countable base of neighborhoods of zero consisting of open
ideals, as per Lemma~\ref{limit-epimorphic-ind-affine-ind-schemes}.
 In this case, the exact category of flat pro-quasi-coherent
pro-sheaves $\fU\Proqc_\fl$ is naturally embedded into an abelian
category, which is called the category of \emph{$\fR$\+contramodules}
and denoted by $\fR\Contra$.

 We skip the definition and introductory discussion of contramodules
over a topological ring, which can be found in
the papers~\cite[Sections~1.1\+-1.2 and~5]{PR},
\cite[Sections~6\+-7]{PS1}, \cite[Section~2]{Pproperf},
\cite[Section~7]{Pflcc}, or~\cite[Section~3]{Pbc}.
 For any complete, separated topological associative rings with
a countable base of neighborhoods of zero consisting of open right
ideals, the category of left $\fR$\+contramodules $\fR\Contra$ is
a locally $\aleph_1$\+presentable (in the sense
of~\cite[Definition~1.17 and Theorem~1.20]{AR}) abelian category
with enough projective objects.

 The relevant full exact subcategory in $\fR\Contra$ is the full
subcategory of \emph{flat} $\fR$\+contramodules $\fR\Contra_\fl
\subset\fR\Contra$ \,\cite[Section~E.1]{Pcosh},
\cite[Sections~6\+-7]{PR}, \cite[Section~3]{Pproperf},
\cite[Section~8]{Pflcc}.
 A relevant description of the category of flat contramodules over
a topological associative ring with a countable base of two-sided
ideals can be found in~\cite[Lemma~8.6(b)]{Pflcc} (see
also~\cite[Lemma~8.4]{Pflcc} for background results providing
the description of the exact category structure on $\fR\Contra_\fl$).

 In the context of an ind-affine $\aleph_0$\+ind-scheme $\fU$ and
the corresponding topological ring $\fR=\fO(\fU)$, \ $\fU=\Spf\fR$,
there is a natural equivalence of exact categories
$$
 \fU\Proqc_\fl\simeq\fR\Contra_\fl.
$$
 To a flat pro-quasi-coherent pro-sheaf $\rF$ on $\fU$,
the flat $\fR$\+contramodule $\fF=\varprojlim_{\gamma\in\Gamma}
\rF^{(U_\gamma)}(U_\gamma)$ is assigned.
 To a flat $\fR$\+contramodule $\fF$, the flat pro-quasi-coherent
pro-sheaf $\rF$ with the quasi-coherent sheaf $\rF^{(U_\gamma)}$
on $U_\gamma$ corresponding to the $\cO(U_\gamma)$\+module
$\cO(U_\gamma)\ocn_\fR\fF=\fF/(\fI_\gamma\tim\fF)$ is
assigned~\cite[Examples~3.8(1\+-2)]{Psemten}.
 Here $\fI_\gamma\subset\fR$ is the open kernel ideal of the natural
surjective ring homomorphism $\fR\rarrow\cO(U_\gamma)$; so
$\cO(U_\gamma)\simeq\fR/\fI_\gamma$.
 See the references above for the explanations of the notation
$\ocn_\fR$ (for the \emph{contratensor product} over~$\fR$)
and $\fI\tim\fF\subset\fF$.

 It is also worth noticing that the full subcategory of flat
contramodules $\fR\Contra_\fl$ is closed under directed colimits
(as well as extensions and kernels of epimorphisms) in the locally
presentable abelian category $\fR\Contra$.
 So the directed colimits in $\fU\Proqc_\fl$ agree with those in
$\fR\Contra_\fl$ and in $\fR\Contra$.
\end{ex}

\subsection{Inverse and direct images}
 Let $\ff\:\fY\rarrow\fX$ be a morphism of ind-schemes.
 For any closed subscheme $W\subset\fY$, there exists a closed
subscheme $Z\subset\fX$ such that the composition $W\rarrow\fY
\rarrow\fX$ factorizes (always uniquely) as $W\rarrow Z\rarrow\fX$.
 So we have a commutative square diagram of morphisms of schemes
and ind-schemes
\begin{equation} \label{morphism-of-indsch-closed-subsches-commutative}
\begin{gathered}
 \xymatrix{
  W \ar@{>->}[r] \ar[d]_{f_{WZ}} & \fY \ar[d]^\ff \\
  Z \ar@{>->}[r] & \fX 
 }
\end{gathered}
\end{equation}
with the horizontal arrows with tails showing closed immersion
morphisms.

 Let $\rF$ be a flat pro-quasi-coherent pro-sheaf on~$\fX$.
 Then, with the notation of
the diagram~\eqref{morphism-of-indsch-closed-subsches-commutative},
the rule $(\ff^*\rF)^{(W)}=f_{WZ}^*(\rF^{(Z)})$ defines a flat
pro-quasi-coherent pro-sheaf $\ff^*\rF$ on~$\fY$.
 The resulting functor
\begin{equation} \label{inverse-image-of-flat-proqc-prosh}
 \ff^*\:\fX\Proqc_\fl\lrarrow\fY\Proqc_\fl
\end{equation}
is called the functor of \emph{inverse image} of flat
pro-quasi-coherent pro-sheaves.
 For any morphism of ind-schemes~$\ff$, the inverse image
functor $\ff^*$ \,\eqref{inverse-image-of-flat-proqc-prosh}
is exact (as a functor of exact categories).

 In particular, if $\fii\:Z\rarrow\fX$ is a closed immersion of
a concentrated scheme $Z$ into an ind-scheme $\fX$, then for
any flat pro-quasi-coherent pro-sheaf $\rF$ on $\fX$ one has
$\fii^*\rF=\rF^{(Z)}$.

 A morphism of schemes $f\:Y\rarrow X$ is said to be \emph{flat} if,
for every pair of affine open subschemes $V\subset Y$ and $U\subset X$
such that $f(V)\subset U$, the related homomorphism of commutative
rings $\cO_X(U)\rarrow\cO_Y(V)$ makes $\cO_Y(V)$ a flat
$\cO_X(U)$\+module.
 Given two of affine open coverings $Y=\bigcup_\alpha V_\alpha$
and $X=\bigcup_\alpha U_\alpha$ of the schemes $Y$ and $X$, indexed
by the same set of indices~$\{\alpha\}$ and such that
$f(V_\alpha)\subset U_\alpha$ for all~$\alpha$, it suffices to check
the flatness condition above for the ring homomorphisms
$\cO_X(U_\alpha)\rarrow\cO_Y(V_\alpha)$ \,\cite[Lemma Tag~01U5]{SP}.
 The class of all flat morphisms of schemes is closed under
pullbacks~\cite[Lemma Tag~01U9]{SP}.

 For any flat affine morphism of schemes $f\:Y\rarrow X$, the direct
image functor~$f_*$ takes flat quasi-coherent sheaves on $Y$ to flat
quasi-coherent sheaves on~$X$.
 The resulting functor $f_*\:Y\Qcoh_\fl\rarrow X\Qcoh_\fl$ is exact
(as a functor of exact categories).

 In the context of
Section~\ref{classes-of-morphisms-of-ind-schemes-subsecn},
take $\bP$ to be the property of a morphism of schemes to be flat.
 This produces the notions of a \emph{flat} and an \emph{ind-flat}
morphism of ind-schemes.

 Let $\ff\:\fY\rarrow\fX$ be a flat affine morphism of ind-schemes.
 Then, for every closed subscheme $Z\subset\fX$, the ind-scheme
$W=Z\times_\fX\fY$ is a concentrated scheme, and the induced morphism
of schemes $f_Z\:W\rarrow Z$ is flat and affine (while the morphism
$W\rarrow\fY$ is a closed immersion).
 We have a pullback diagram of morphisms of schemes and ind-schemes
\begin{equation} \label{morphism-of-indsch-closed-subsches-pullback}
\begin{gathered}
 \xymatrix{
  W \ar@{>->}[r] \ar[d]_{f_Z} & \fY \ar[d]^\ff \\
  Z \ar@{>->}[r] & \fX 
 }
\end{gathered}
\end{equation}
which is a special case of
the diagram~\eqref{morphism-of-indsch-closed-subsches-commutative}.

 Let $\rG$ be a flat pro-quasi-coherent pro-sheaf on~$\fY$.
 Then, with the notation of the pullback
diagram~\eqref{morphism-of-indsch-closed-subsches-pullback},
the rule $(\ff_*\rG)^{(Z)}=(f_Z)_*(\rG^{(W)})$ defines a flat
pro-quasi-coherent pro-sheaf $\ff_*\rG$ on~$\fX$.
 One needs to use the affine base change lemma for quasi-coherent
sheaves~\cite[Lemma Tag~02KG]{SP}, \cite[Lemma~3.3(a)]{Psemten}
in order to construct the isomorphisms~(ii) from
Section~\ref{flat-proqcoh-prosh-subsecn} above for the flat
quasi-coherent sheaves $(\ff_*\rG)^{(Z)}$ defined in this way.

 The resulting functor
\begin{equation} \label{direct-image-of-flat-proqc-prosh}
 \ff_*\:\fY\Proqc_\fl\lrarrow\fX\Proqc_\fl
\end{equation}
is called the functor of \emph{direct image} of flat
pro-quasi-coherent pro-sheaves.
 For any flat affine morphism of ind-schemes~$\ff$, the direct image
functor $\ff_*$ \,\eqref{direct-image-of-flat-proqc-prosh}
is exact (as a functor of exact categories).

 For any flat affine morphism of ind-schemes~$\ff$, the direct image
functor~$\ff_*$ \,\eqref{direct-image-of-flat-proqc-prosh}
is right adjoint to the inverse image functor
$\ff^*$~\,\eqref{inverse-image-of-flat-proqc-prosh}
\,\cite[Sections~3.3 and~3.4]{Psemten}.

\subsection{Very flat modules and ring homomorphisms}
\label{very-flat-modules-and-ring-maps-subsecn}
 The concept of very flatness was introduced in the present author's
preprint on contraherent cosheaves~\cite[Sections~1.1, 1.10,
and~4.1]{Pcosh}.
 Further significant developments in this context are
the papers~\cite{ST}, \cite{PSl1}, and~\cite{ES}.
 For discussions of the complementary concept of contraadjustedness,
see~\cite[Sections~1.1, 2.5, and~4.1]{Pcosh} and
the papers~\cite{ST,Pcta,Pal}.
 Another reference source is the preprint~\cite{Pphil}.

 All rings and modules in this paper are associative and unital.
 Given an associative ring $R$, we denote by $R\Modl$ the abelian
category of left $R$\+modules.
 The full subcategory of flat left $R$\+modules is denoted by
$R\Modl_\fl\subset R\Modl$.

 A left $R$\+module $C$ is called \emph{cotorsion} (in
the sense of Enochs~\cite{En2}) if $\Ext^1_R(F,C)=0$ for all flat
left $R$\+modules $C$, or equivalently, $\Ext^n_R(F,C)=0$ for all
flat left $R$\+modules $C$ and all $n\ge1$.
 The class of cotorsion $R$\+modules is closed under extensions,
cokernels of injective morphisms, and infinite direct products in
$R\Modl$.
 We denote the full subcategory of cotorsion left $R$\+modules by
$R\Modl^\cot\subset R\Modl$.
 The full subcategories $R\Modl_\fl$ and $R\Modl^\cot$ inherit
exact category structures from the abelian exact structure of $R\Modl$.

 Let $R$ be a commutative ring.
 Given an element $s\in R$, we denote by $R[s^{-1}]$ the ring obtained
from $R$ by formally inverting the element~$s$.
 In other words, $R[s^{-1}]$ is the ring of regular functions on
the principal affine open subscheme in $\Spec R$ corresponding to~$s$,
or the localization of $R$ at the multiplicative subset
$\{1,s,s^2,s^3,\dotsc\}$.
 The $R$\+module $R[s^{-1}]$ is flat of projective dimension
at most~$1$ \,\cite[proof of Lemma~2.1]{Pcta}, \cite[Remark~5.1]{Pphil}.

 An $R$\+module $C$ is said to be \emph{contraadjusted} if
$\Ext^1_R(R[s^{-1}],C)=0$ for all $s\in R$.
 The class of contraadjusted $R$\+modules is closed under extensions
and infinite direct products in $R\Modl$.
 Any quotient $R$\+module of a contraadjusted $R$\+module is
contraadjusted.
 All cotorsion $R$\+modules are contraadjusted.
 We denote the full subcategory of contraadjusted $R$\+modules
by $R\Modl^\cta\subset R\Modl$.

 An $R$\+module $F$ is said to be \emph{very flat} if
$\Ext^1_R(F,C)=0$ for all contraadjusted $R$\+modules~$C$.
 The class of very flat $R$\+modules is closed under extensions,
kernels of surjective morphisms, and infinite direct sums in
$R\Modl$.
 All very flat $R$\+modules have projective dimensions at most~$1$
\,\cite[Section~1.1]{Pcosh}.
 All very flat $R$\+modules are flat.
 We denote the full subcategory of very flat $R$\+modules by
$R\Modl_\vfl\subset R\Modl_\fl\subset R\Modl$.
 The full subcategories $R\Modl^\cta$ and $R\Modl_\vfl$ inherit
exact category structures from the abelian exact structure of $R\Modl$.

 The following characterization of very flat
$R$\+modules~\cite[Corollary~1.1.4]{Pcosh}, based on the results of
the paper~\cite{ET}, may be more intuitive than the abstract
homological algebra definition above.
 An $R$\+module $G$ is very flat if and only if it is a direct
summand of an $R$\+module $F$ admitting an increasing filtration
$(F_i\subset F)_{i\in\alpha}$, indexed by some ordinal~$\alpha$,
with the successive quotient $R$\+modules $F_{i+1}/F_i$ of the form
$F_{i+1}/F_i\simeq R[s_i^{-1}]$ for some elements $s_i\in R$,
\,$0\le i<\alpha$.
 Here the notation presumes that the filtration is given by
a smooth chain of $R$\+submodules with $F_0=0$, \
$F_i\subset F_j$ for all $0\le i\le j\le\alpha$,
\ $F_j=\bigcup_{i<j}F_i$ for all limit ordinals $j\le\alpha$,
and $F_\alpha=F$.

 For any affine scheme $U$ with an affine open subscheme $V\subset U$,
the ring $\cO(V)$ is very flat as
an $\cO(U)$\+module~\cite[Lemma~1.2.4]{Pcosh}.
 Much more generally, for any homomorphism of commutative rings
$R\rarrow S$ such that $S$ is a finitely presented $R$\+algebra
and a flat $R$\+module, the $R$\+module $S$ is actually very
flat~\cite[Main Theorem~1.1]{PSl1}.

 A homomorphism of commutative rings $R\rarrow S$ is said to be
\emph{very flat} if, for every element $s\in S$, the ring $S[s^{-1}]$
is very flat as an $R$\+module~\cite[Section~9]{PSl1},
\cite[Section~6.1]{Pphil}.
 If this is the case, the ring $S$ is also said to be
a \emph{very flat commutative $R$\+algebra}.
 This condition first appears in~\cite[Lemma~1.2.3]{Pcosh}.
 Requiring very flatness of the $R$\+module $S$ itself is \emph{not}
sufficient here~\cite[Example~9.7]{PSl1}, \cite[Section~6.1]{Pphil}.

\subsection{Very flat sheaves and morphisms of schemes}
 A quasi-coherent sheaf $\cF$ on a scheme $X$ is said to be
\emph{very flat} if, for every affine open subscheme $U\subset X$,
the $\cO_X(U)$\+module $\cF(U)$ is very flat.
 It suffices to check this condition for affine open subschemes $U$
belonging to any chosen affine open covering of the scheme~$X$
\,\cite[Lemma~1.2.6(a) and Section~1.10]{Pcosh},
\cite[Example~2.5]{Pal}.

 We denote the full subcategory of very flat quasi-coherent sheaves
by $X\Qcoh_\vfl\subset X\Qcoh_\fl\subset X\Qcoh$.
 The full subcategory $X\Qcoh_\vfl$ is closed under extensions,
kernels of epimorphisms, and infinite coproducts in the abelian
category $X\Qcoh$ (and in its full subcategory $X\Qcoh_\fl$);
so $X\Qcoh_\fl$ inherits an exact category structure from the abelian
exact structure of $X\Qcoh$.
 
 For any morphism of schemes $f\:Y\rarrow X$, the inverse image functor
$f^*\:X\Qcoh\allowbreak\rarrow Y\Qcoh$ takes very flat quasi-coherent
sheaves on $X$ to very flat quasi-coherent sheaves on~$Y$
\,\cite[Lemma~1.2.2(b) and Section~1.10]{Pcosh}.

 A morphism of schemes $f\:Y\rarrow X$ is said to be \emph{very flat}
if, for every pair of affine open subschemes $V\subset Y$ and
$U\subset X$ such that $f(V)\subset U$, the $\cO_X(U)$\+algebra
$\cO_Y(V)$ is very flat.
 Requiring this for all pairs of affine open subschemes $V\subset Y$
and $U\subset X$ as above is equivalent to requiring
the $\cO_X(U)$\+module $\cO_Y(V)$ to be very flat for all such
pairs of affine open subschemes.

 Let $X=\bigcup_\alpha U_\alpha$ be an open covering of~$X$.
 Then a morphism of schemes $f\:Y\rarrow X$ is very flat if and only
if, for every index~$\alpha$ and every affine open subscheme
$V\subset Y$ such that $f(V)\subset U_\alpha$,
the $\cO_X(U_\alpha)$\+module $\cO_Y(V)$ is very
flat~\cite[Corollary~1.2.5(b), Lemma~1.2.6(a), and Section~1.10]{Pcosh}.

 Let $X=\bigcup_\alpha U_\alpha$ and $Y=\bigcup_\alpha V_\alpha$
be affine open coverings of the schemes $X$ and $Y$, indexed by
the same set of indices~$\alpha$, such that one has $f(V_\alpha)
\subset U_\alpha$ for every~$\alpha$.
 Then the morphism of schemes~$f$ is very flat if and only if, for
every index~$\alpha$, the $\cO_X(U_\alpha)$\+algebra $\cO_Y(V_\alpha)$
is very flat~\cite[Lemma~1.8.6 and Section~1.10]{Pcosh}.

 For any very flat affine morphism of schemes $f\:Y\rarrow X$,
the direct image functor $f_*\:Y\Qcoh\rarrow X\Qcoh$ takes very flat
quasi-coherent sheaves on $Y$ to very flat quasi-coherent sheaves
on~$X$ \,\cite[Lemma~1.2.3(b) and Section~1.10]{Pcosh}.

\subsection{Very flat pro-quasi-coherent pro-sheaves}
 Let $\fX$ be an ind-scheme.
 A flat pro-quasi-coherent pro-sheaf $\rF$ on $\fX$ is said to be
\emph{very flat} if, for every closed subscheme $Z\subset\fX$,
the flat quasi-coherent sheaf $\rF^{(Z)}$ on $Z$ is very flat.

 Given a representation of $\fX$ by a directed diagram of closed
immersions of concentrated schemes $(X_\gamma)_{\gamma\in\Gamma}$,
it suffices to check the very flatness condition above for the closed
subschemes $Z=X_\gamma\subset\fX$, \,$\gamma\in\Gamma$.
 This follows from the stability of the class of very flat
quasi-coherent sheaves under inverse images with respect to
morphisms of schemes.

 We denote the full subcategory of very flat pro-quasi-coherent
pro-sheaves by $\fX\Proqc_\vfl\subset\fX\Proqc_\fl$.
 The full subcategory of very flat pro-quasi-coherent pro-sheaves
$\fX\Proqc_\vfl$ is closed under extensions, kernels of admissible
epimorphisms, and infinite coproducts in the exact category of flat
pro-quasi-coherent pro-sheaves $\fX\Proqc_\fl$; so the full subcategory
$\fX\Proqc_\vfl$ inherits an exact category structure from
$\fX\Proqc_\fl$.

 For any morphism of ind-schemes $\ff\:\fY\rarrow\fX$, the inverse
image functor $\ff^*$ \,\eqref{inverse-image-of-flat-proqc-prosh}
takes very flat pro-quasi-coherent pro-sheaves on $\fX$ to very flat
pro-quasi-coherent pro-sheaves on~$\fY$.
 So we have an exact inverse image functor
\begin{equation} \label{inverse-image-of-very-flat-proqc-prosh}
 \ff^*\:\fX\Proqc_\vfl\lrarrow\fY\Proqc_\vfl.
\end{equation}

 Our next aim is to define the notion of a very flat morphism of
ind-schemes.
 The property of a morphism of schemes to be very flat is \emph{not}
known to be stable under pullbacks; so the theory of
Section~\ref{classes-of-morphisms-of-ind-schemes-subsecn} is
not applicable.
 However, the very flatness of morphisms of schemes is preserved by
the pullbacks with respect to locally closed immersions of
schemes~\cite[Lemma~2.10.1]{Pcosh}.

 Let $\ff:\fY\rarrow\fX$ be a morphism of ind-schemes that is
``representable by (concentrated) schemes''.
 So, for every closed subscheme $Z\subset\fX$, the ind-scheme
$W=Z\times_\fX\fY$ is a concentrated scheme (and a closed subscheme
in~$\fY$).
 We will say that the morphism of ind-schemes~$\ff$ is \emph{very flat}
if, for every closed subscheme $Z\subset\fX$, the morphism of
schemes $W=Z\times_\fX\fY\rarrow Z$ is very flat.

 Once again, given a representation of $\fX$ by a directed diagram of
closed immersions of concentrated schemes
$(X_\gamma)_{\gamma\in\Gamma}$, it suffices to check the very flatness
condition above for the closed subschemes $Z=X_\gamma\subset\fX$,
\,$\gamma\in\Gamma$.
 This follows from the stability of the class of very flat morphisms
of schemes under the pullbacks by closed immersions.

 In particular, all open immersions of ind-schemes are very flat
morphisms (since all open immersions of schemes are very flat
morphisms).

 For any very flat affine morphism of ind-schemes $\ff\:\fY\rarrow\fX$,
the direct image functor $\ff_*$
\,\eqref{direct-image-of-flat-proqc-prosh} takes very flat
pro-quasi-coherent pro-sheaves on $\fY$ to very flat pro-quasi-coherent
pro-sheaves on~$\fX$.
 So we have an exact direct image functor
\begin{equation} \label{direct-image-of-very-flat-proqc-prosh}
 \ff_*\:\fY\Proqc_\vfl\lrarrow\fX\Proqc_\vfl.
\end{equation}

\begin{ex} \label{ind-affine-ind-schemes-and-very-flat-contramodules}
 We continue the discussion of ind-affine $\aleph_0$\+ind-schemes
started in Lemma~\ref{limit-epimorphic-ind-affine-ind-schemes}
and Example~\ref{ind-affine-ind-schemes-and-flat-contramodules},
and keep the notation introduced there.

 Let $\fR$ be a complete, separated topological commutative ring with
a countable base of neighborhoods of zero consisting of open ideals.
 Then an $\fR$\+contramodule $\fF$ is said to be \emph{very flat} if,
for every open ideal $\fI\subset\fR$, the $\fR/\fI$\+module
$\fR/\fI\ocn_\fR\fF\simeq\fF(\fI\tim\fF)$ is very
flat~\cite[Section~E.4]{Pcosh}, \cite[Example~7.12(3)]{PR}
(cf.~\cite[Section~2.12]{Pform} for a special case).
 It suffices to check this condition for open ideals $\fI$ belonging
to any chosen base of neighborhoods of zero in $\fR$ consisting of
open ideals.

 All very flat contramodules are flat (since all very flat modules
over commutative rings are flat).
 The full subcategory of very flat $\fR$\+contramodules
$\fR\Contra_\vfl$ is closed under extensions and kernels of
epimorphisms in $\fR\Contra$ (or, which is the same, closed under
extensions and kernels of admissible epimorphisms in $\fR\Contra_\fl$),
so it inherits an exact category structure.

 Let $\fU=\Spf\fR$ be the strict ind-affine $\aleph_0$\+ind-scheme
corresponding to~$\fR$.
 Then the equivalence of exact categories from
Example~\ref{ind-affine-ind-schemes-and-flat-contramodules}
restricts to an equivalence of exact subcategories
$$
 \fU\Proqc_\vfl\simeq\fR\Contra_\vfl.
$$
\end{ex}

\Section{Accessibility and Pseudolimits}  \label{accessibility-secn}

\subsection{Background on accessible categories}
 Let $\kappa$~be a regular cardinal.
 We refer to the book~\cite[Theorem~1.5, Definition~1.13, and
Remark~1.21]{AR} and the paper~\cite{AN} for a discussion of
the definitions of, the difference between, and the essential
equivalence of the \emph{$\kappa$\+directed} and
\emph{$\kappa$\+filtered} colimits in categories.
 Let us only mention, to avoid any confusion, that a poset $\Gamma$
is said to be $\kappa$\+directed if any subposet of cardinality
smaller than~$\kappa$ has an upper bound in $\Gamma$.

 Let $\sK$ be a category with $\kappa$\+directed colimits.
 An object $K\in\sK$ is said to be \emph{$\kappa$\+presentable}
if the functor $\Mor_\sK(K,{-})\:\sK\rarrow\Sets$ preserves
$\kappa$\+directed (or equivalently, $\kappa$\+filtered) colimits.

 The category $\sK$ is said to be \emph{$\kappa$\+accessible} if
there exists a \emph{set} (rather than a proper class) of objects
$\sS\subset\sK$ such that all the objects of $\sS$ are
$\kappa$\+presentable in $\sK$ and every object of $\sK$ is
a $\kappa$\+directed colimit of objects from~$\sS$
\,\cite[Definition~2.1]{AR}.
 If this is the case, then the $\kappa$\+presentable objects of
$\sK$ are precisely all the retracts of the objects from~$\sS$.

 A category is said to be \emph{accessible} if it is
$\kappa$\+accessible for some regular cardinal~$\kappa$
\,\cite[Remark~2.2(1)]{AR}.
 The $\aleph_0$\+accessible categories are called \emph{finitely
accessible}, and the $\aleph_1$\+accessible categories can be
similarly called \emph{countably accessible}.
 The same terminology applies to objects: the $\aleph_0$\+presentable
objects are called \emph{finitely presentable}, and
the $\aleph_1$\+presentable objects are likewise called
\emph{countably presentable}.

 In a $\kappa$\+accessible category $\sK$, every object $K\in\sK$ has
a presentability rank~\cite[[Remark~2.2(3)]{AR}.
 This means that for every $K\in\sK$ there exists a regular cardinal
$\rho\ge\kappa$ such that the object $K$ is $\rho$\+presentable.

 A $\kappa$\+accessible category where all colimits exist is said to
be \emph{locally\/ $\kappa$\+presentable}~\cite[Definition~1.17
and Theorem~1.20]{AR}.
 A category is said to be \emph{locally presentable} if it is
locally $\kappa$\+presentable for some regular cardinal~$\kappa$.
 The locally $\aleph_0$\+presentable categories are called
\emph{locally finitely presentable}~\cite[Definition~1.9 and
Theorem~1.11]{AR}.

\subsection{Background of pseudolimits}
\label{pseudolimits-subsecn}
 Accessibility of pseudolimits of accessible categories is
a result from the book~\cite[Theorem~5.1.6]{MP}.
 More precise bounds on the accessibility ranks of the pseudolimits
are provided by results of the unpublished old
preprint~\cite[Theorem~3.8, Corollary~3.9, and Remark~3.11(II)]{Ulm},
rediscovered in the recent paper~\cite{Pacc}.
 An introductory discussion can be found in~\cite[Section~2.H]{AR}.

 In this section we do not discuss any weighted pseudolimits/bilimits
generally, but rather concentrate on the concrete constructions of
the products, \emph{inserters}, and \emph{equifiers}, as well as
the \emph{pseudopullbacks} and what we call the \emph{isomorphers}
(otherwise known as \emph{iso-inserters}) of categories.

\begin{prop} \label{product-accessible-prop}
 Let\/ $\kappa$~be a regular cardinal, $\Xi$ be a set of cardinality
smaller than\/~$\kappa$, and $(\sK_\xi)_{\xi\in\Xi}$ be
a\/ $\Xi$\+indexed family of\/ $\kappa$\+accessible categories.
 Then the Cartesian product category\/ $\sK=\prod_{\xi\in\Xi}\sK_\xi$
is\/ $\kappa$\+accessible.
 The\/ $\kappa$\+presentable objects of\/ $\sK$ are precisely all
the families of objects $(S_\xi\in\sK_\xi)_{\xi\in\Xi}$ such that
the object $S_\xi$ is\/ $\kappa$\+presentable in\/ $\sK_\xi$ for
every index\/ $\xi\in\Xi$.
\end{prop}

\begin{proof}
 This is a corrected version of~\cite[Proposition~2.67]{AR},
as per~\cite[Proposition~2.1]{Pacc}.
\end{proof}

 Let $\lambda$~be an ordinal and $\sK$ be a category.
 A \emph{$\lambda$\+indexed chain} (of objects and morphisms) in $\sK$
is a functor $\lambda\rarrow\sK$.
 Here we are using the observations that any ordinal can be viewed as
a poset, and any poset $\Gamma$ can be viewed as a category with
a morphism $\gamma'\rarrow\gamma''$ for every pair of elements
$\gamma'\le\gamma''\in\Gamma$; so an ordinal~$\lambda$ can be viewed
as a category.
 In particular, this applies to any cardinal~$\lambda$.

 Let $\sK$ and $\sL$ be two categories, and let $F$, $G\:\sK
\rightrightarrows\sL$ be two functors from $\sK$ to~$\sL$.
 In this context, the \emph{inserter category} $\sD$ is defined
as the category of pairs $(K,\theta)$, where $K\in\sK$ is
an object and $\theta\:F(K)\rarrow G(K)$ is a morphism in~$\sL$.
 Morphisms in $\sD$ are defined in the obvious way.

\begin{prop} \label{inserter-accessible-prop}
 Let\/ $\kappa$ be a regular cardinal and\/ $\lambda$ be an infinite
cardinal smaller than\/~$\kappa$.
 Let\/ $\sK$ and\/ $\sL$ be two\/ $\kappa$\+accessible categories
where colimits of all\/ $\lambda$\+indexed chains exist.
  Let $F$, $G\:\sK\rightrightarrows\sL$ be a parallel pair of functors
preserving\/ $\kappa$\+directed colimits and colimits of
$\lambda$\+indexed chains, and taking\/ $\kappa$\+presentable objects
to\/ $\kappa$\+presentable objects.
 Then the inserter category\/ $\sD$ is\/ $\kappa$\+accessible.
 An object $(K,\theta)\in\sD$ is\/ $\kappa$\+presentable if and only
if the object $K$ is\/ $\kappa$\+presentable in\/~$\sK$.
\end{prop}

\begin{proof}
 This is~\cite[Theorem~3.8, Corollary~3.9, and Remark~3.11(II)]{Ulm}
or~\cite[Theorem~4.1]{Pacc}.
\end{proof}

 Let $\sD$ and $\sM$ be two categories, let
$F$, $G\:\sD\rightrightarrows\sM$ be two functors, and let
$\phi$, $\psi\:F\rightrightarrows G$ be two natural transformations.
 In this context, the \emph{equifier category} $\sE$ is the full
subcategory $\sE\subset\sD$ consisting of all the objects $D\in\sD$
for which the two morphisms $\phi_D$ and $\psi_D\:F(D)\rightrightarrows
G(D)$ are equal to each other.

\begin{prop} \label{equifier-accessible-prop}
 Let\/ $\kappa$ be a regular cardinal and\/ $\lambda$ be an infinite
cardinal smaller than\/~$\kappa$.
 Let\/ $\sD$ and\/ $\sM$ be two\/ $\kappa$\+accessible categories
where colimits of all\/ $\lambda$\+indexed chains exist.
  Let $F$, $G\:\sD\rightrightarrows\sM$ be a parallel pair of functors
preserving\/ $\kappa$\+directed colimits and colimits of
$\lambda$\+indexed chains, and taking\/ $\kappa$\+presentable objects
to\/ $\kappa$\+presentable objects.
 Let\/ $\phi$, $\psi\:F\rightrightarrows G$ be two natural
transformations.
 Then the equifier category\/ $\sE$ is\/ $\kappa$\+accessible.
 An object of\/ $\sE$ is\/ $\kappa$\+presentable if and only if it is\/
$\kappa$\+presentable in\/~$\sD$.
\end{prop}

\begin{proof}
 This is~\cite[Theorem~3.8, Corollary~3.9, and Remark~3.11(II)]{Ulm}
or~\cite[Theorem~3.1]{Pacc}.
\end{proof}

 Let $\sK$, $\sL$, and $\sM$ be three categories, and let
$F\:\sK\rarrow\sM$ and $G\:\sL\rarrow\sM$ be two functors.
 In this context, the \emph{pseudopullback category} $\sP$ is
defined as the category of triples $(K,L,\theta)$, where $K\in\sK$
and $L\in\sL$ are two objects and $\theta\:F(K)\simeq G(L)$ is
an isomorpsm in\/~$\sM$.
 Morphisms in $\sP$ are defined in the obvious way.

\begin{cor} \label{pseudopullback-accessible-cor}
 Let\/ $\kappa$ be a regular cardinal and\/ $\lambda$ be an infinite
cardinal smaller than\/~$\kappa$.
 Let\/ $\sK$, $\sL$, and\/ $\sM$ be three\/ $\kappa$\+accessible
categories where colimits of all\/ $\lambda$\+indexed chains exist.
 Let $F\:\sK\rarrow\sM$ and $G\:\sL\rarrow\sM$ be two functors
preserving\/ $\kappa$\+directed colimits and colimits of
$\lambda$\+indexed chains, and taking\/ $\kappa$\+presentable objects
to\/ $\kappa$\+presentable objects.
 Then the pseupopullback category\/ $\sP$ is\/ $\kappa$\+accessible.
 An object $(K,L,\theta)\in\sP$ is\/ $\kappa$\+presentable if and only
if the object $K$ is\/ $\kappa$\+presentable in\/ $\sK$ and
the object $L$ is\/ $\kappa$\+presentable in\/~$\sL$.
\end{cor}

\begin{proof}
 This is~\cite[Remark~3.2(I), Theorem~3.8, Corollary~3.9, and
Remark~3.11(II)]{Ulm}, or~\cite[Pseudopullback Theorem~2.2]{RR}
(based on~\cite[Proposition~3.1]{CR}), or~\cite[Corollary~5.1]{Pacc}.
 This corollary can be easily deduced from
Propositions~\ref{product-accessible-prop},
\ref{inserter-accessible-prop}, and~\ref{equifier-accessible-prop}.
\end{proof}

 Let $\sK$ and $\sL$ be two categories, and let $F$, $G\:\sK
\rightrightarrows\sL$ be two functors from $\sK$ to~$\sL$.
 In this context, the \emph{isomorpher} (or \emph{iso-inserter})
\emph{category} $\sD$ is defined as the category of pairs $(K,\theta)$,
where $K\in\sK$ is an object and $\theta\:F(K)\simeq G(K)$ is
an isomorphism in~$\sL$.
 Morphisms in $\sD$ are defined in the obvious way.

\begin{cor} \label{isomorpher-accessible-cor}
 Let\/ $\kappa$ be a regular cardinal and\/ $\lambda$ be an infinite
cardinal smaller than\/~$\kappa$.
 Let\/ $\sK$ and\/ $\sL$ be two\/ $\kappa$\+accessible categories
where colimits of all\/ $\lambda$\+indexed chains exist.
  Let $F$, $G\:\sK\rightrightarrows\sL$ be a parallel pair of functors
preserving\/ $\kappa$\+directed colimits and colimits of
$\lambda$\+indexed chains, and taking\/ $\kappa$\+presentable objects
to\/ $\kappa$\+presentable objects.
 Then the isomorpher category\/ $\sD$ is\/ $\kappa$\+accessible.
 An object $(K,\theta)\in\sD$ is\/ $\kappa$\+presentable if and only
if the object $K$ is\/ $\kappa$\+presentable in\/~$\sK$.
\end{cor}

\begin{proof}
 This is~\cite[Remark~3.2(I), Theorem~3.8, Corollary~3.9, and
Remark~3.11(II)]{Ulm} or~\cite[Remark~5.2]{Pacc}.
 This corollary is essentially an equivalent restatement of
the previous one.
\end{proof}

 Given an additive category $\sK$, we denote by $\Com(\sK)$
the additive category of (unbounded) complexes in~$\sK$.
 An additive category is said to be \emph{$\kappa$\+accessible}
if it is $\kappa$\+accessible as an abstract category.

\begin{prop} \label{category-of-complexes-accessible-prop}
 Let\/ $\kappa$~be a regular cardinal and\/ $\lambda$ be an infinite
cardinal smaller than\/~$\kappa$.
 Let\/ $\sK$ be a\/ $\kappa$\+accessible additive category where
colimits of all\/ $\lambda$\+indexed chains exist.
 Then the additive category\/ $\Com(\sK)$ of complexes in\/ $\sK$ is\/
$\kappa$\+accessible.
 A complex $K^\bu$ in\/ $\sK$ is\/ $\kappa$\+presentable as an object
of\/ $\Com(\sK)$ if and only if all its terms $K^n$, \,$n\in\boZ$,
are\/ $\kappa$\+presentable in\/~$\sK$.
\end{prop}

\begin{proof}
 This is a special case of~\cite[Theorem~6.2]{Pacc};
cf.~\cite[proof of Corollary~10.4]{Pacc}.
\end{proof}

\begin{lem} \label{category-of-complexes-locally-presentable-lem}
 Let\/ $\kappa$~be a regular cardinal and\/ $\sK$ be a locally\/
$\kappa$\+presentable additive category.
 Then the additive category\/ $\Com(\sK)$ of complexes in\/ $\sK$ is
locally\/ $\kappa$\+presentable.
 In the case of an uncountable regular cardinal\/ $\kappa\ge\aleph_1$,
a complex $K^\bu$ in\/ $\sK$ is\/ $\kappa$\+presentable as an object
of\/ $\Com(\sK)$ if and only if all its terms $K^n$, \,$n\in\boZ$,
are\/ $\kappa$\+presentable in\/~$\sK$.
 In the case of the countable cardinal\/ $\kappa=\aleph_0$,
a complex $K^\bu$ in\/ $\sK$ is finitely presentable as an object
of\/ $\Com(\sK)$ if and only if $K^\bu$ \emph{is bounded and} all its
terms $K^n$, \,$n\in\boZ$, are finitely presentable in\/~$\sK$.
\end{lem}

\begin{proof}
 The case of $\kappa\ge\aleph_1$ is covered by
Proposition~\ref{category-of-complexes-accessible-prop}
(take $\lambda=\aleph_0$).
 In the case of $\kappa=\aleph_0$, the first assertion of
the lemma is covered by~\cite[Lemma~6.3]{PS4}, and the second
assertion can be discerned from the same argument
using~\cite[Lemma~A.3(2)]{Sto0}.
 Another relevant reference is~\cite[Theorem~1.2]{Hen}.
\end{proof}

\subsection{Background on accessibility over schemes}
 Let $R$ be an associative ring.
 Then the abelian category $R\Modl$ is locally finitely presentable,
hence also locally $\kappa$\+presentable for all regular
cardinals~$\kappa$.
 An $R$\+module is $\kappa$\+presentable if and only if it is
(isomorphic to) the cokernel of a morphism of two free $R$\+modules
with fewer than~$\kappa$ generators.

 Let $X$ be a scheme.
 A quasi-coherent sheaf $\cM$ on $X$ is said to be \emph{locally\/
$\kappa$\+presentable} if, for every affine open subscheme
$U\subset X$, the $\cO_X(U)$\+module $\cM(U)$ is $\kappa$\+presentable.
 It suffices to check this condition for affine open subschemes $U$
belonging to any chosen affine open covering of the scheme~$X$
\,\cite[Lemma~3.1]{PS6}.

 For any morphism of schemes $f\:Y\rarrow X$, the inverse image
functor $f^*\:X\Qcoh\allowbreak\rarrow Y\Qcoh$ takes locally
$\kappa$\+presentable quasi-coherent sheaves on $X$ to locally
$\kappa$\+presentable quasi-coherent sheaves on~$Y$.

 Let $X$ be a concentrated scheme.
 Then the abelian category $X\Qcoh$ is locally finitely
presentable~\cite[0.5.2.5 and Corollaire~I.6.9.12]{EGA1},
\cite[Definition Tag~01BN and Lemma Tag~01PJ or~01PK]{SP},
hence also locally $\kappa$\+presentable for all regular
cardinals~$\kappa$.
 An quasi-coherent sheaf on $X$ is $\kappa$\+presentable if and only
if it is locally $\kappa$\+presentable in the sense of the previous
paragraph (see also~\cite[Remark~4.6]{PS6}).

\begin{prop} \label{flat-qcoh-sheaves-accessible}
 Let $X$ be a concentrated scheme and\/ $\kappa\ge\aleph_1$ be
an uncountable regular cardinal.
 Then the category $X\Qcoh_\fl$ of flat quasi-coherent sheaves on
$X$ is\/ $\kappa$\+accessible.
 A flat quasi-coherent sheaf on $X$ is\/ $\kappa$\+presentable as
an object of $X\Qcoh_\fl$ if and only if it is locally\/
$\kappa$\+presentable (or equivalently, $\kappa$\+presentable as
an object of $X\Qcoh$).
\end{prop}

\begin{proof}
 This is essentially~\cite[Theorem~3.5]{PS6} (see
also~\cite[Theorem~4.5]{PS6}).
\end{proof}

 For an exact category $\sE$, we denote by $\Ac(\sE)\subset\Com(\sE)$
the full subcategory of acyclic complexes in\/~$\sE$.
 In particular, the notation $\Ac(X\Qcoh_\fl)$ stands for the category
of all acyclic complexes in the exact category $X\Qcoh$ of flat
quasi-coherent sheaves on~$X$.
 Notice that a complex in $X\Qcoh_\fl$ that is acyclic in $X\Qcoh$
\emph{need not} be acyclic in $X\Qcoh_\fl$, that is
$\Ac(X\Qcoh_\fl)\ne\Com(X\Qcoh_\fl)\cap\Ac(\Qcoh)$.
 The complexes acyclic in $X\Qcoh_\fl$ are otherwise known as
\emph{pure acyclic complexes} of flat quasi-coherent sheaves.

\begin{prop} \label{complexes-of-flat-qcoh-sheaves-accessible}
 Let $X$ be a concentrated scheme and\/ $\kappa\ge\aleph_1$ be
an uncountable regular cardinal.
 Then the category\/ $\Com(X\Qcoh_\fl)$ of complexes of flat
quasi-coherent sheaves on $X$ is\/ $\kappa$\+accessible.
 A complex\/ $\cF^\bu$ of flat quasi-coherent sheaves on $X$ is\/
$\kappa$\+presentable as an object of\/ $\Com(X\Qcoh_\fl)$ if and
only if all its terms\/ $\cF^n$, \,$n\in\boZ$, are locally\/
$\kappa$\+presentable quasi-coherent sheaves.
\end{prop}

\begin{proof}
 This is essentially~\cite[Theorem~5.1]{PS6}.
 Alternatively, the assertion is provable by comparing
Propositions~\ref{category-of-complexes-accessible-prop}
and~\ref{flat-qcoh-sheaves-accessible}.
\end{proof}

\begin{prop} \label{pure-acycl-cmplxs-of-flat-qcoh-sheaves-accessible}
 Let $X$ be a concentrated scheme and\/ $\kappa\ge\aleph_1$ be
an uncountable regular cardinal.
 Then the category\/ $\Ac(X\Qcoh_\fl)$ of pure acyclic complexes of
flat quasi-coherent sheaves on $X$ is\/ $\kappa$\+accessible.
 A pure acyclic complex\/ $\cF^\bu$ of flat quasi-coherent sheaves on
$X$ is\/ $\kappa$\+presentable as an object of\/ $\Ac(X\Qcoh_\fl)$
if and only if all its terms\/ $\cF^n$, \,$n\in\boZ$, are locally\/
$\kappa$\+presentable quasi-coherent sheaves.
\end{prop}

\begin{proof}
  This is essentially~\cite[Theorem~5.2]{PS6}.
\end{proof}

\subsection{Accessibility of flat pro-quasi-coherent pro-sheaves}
\label{accessibility-of-flat-proqcoh-prosheaves-subsecn}
 Let $\mu$ be an infinite cardinal.
 An (ind-concentrated) ind-scheme $\fX$ is said to be
an \emph{$\mu$\+ind-scheme} if $\fX$ can be represented by a directed
diagram of concentrated schemes $(X_\gamma)_{\gamma\in\Gamma}$ indexed
by a poset $\Gamma$ of cardinality smaller than or equal to~$\mu$.

 The obvious version of Lemma~\ref{strict-aleph0-ind-schemes-lemma}
is valid for directed (sub)posets of cardinality~$\le\mu$ in place
of countable directed (sub)posets.
 The proof is the same, based on the fact that for any directed
poset $\Lambda$ and any subset $\Lambda''\subset\Lambda$ of
cardinality~$\le\mu$ there is a directed subposet
$\Lambda'\subset\Lambda$ of cardinality~$\le\mu$ such that for
every $\lambda''\in\Lambda''$ there exists $\lambda'\in\Lambda'$ with
$\lambda''\le\lambda'$ in~$\Lambda$.
 In this sense, the notion of a strict $\mu$\+ind-scheme is
well-behaved.

 The following theorem is the main result of
Section~\ref{accessibility-secn}.

\begin{thm} \label{flat-proqcoh-prosheaves-accessible}
 Let\/ $\mu$ be an infinite cardinal, $\kappa\ge\aleph_1$, $\kappa>\mu$
be a regular cardinal, and\/ $\fX$ be an (ind-concentrated)
$\mu$\+ind-scheme.
 Then the category\/ $\fX\Proqc_\fl$ of flat pro-quasi-coherent
pro-sheaves on\/ $\fX$ is\/ $\kappa$\+accessible.
 A flat pro-quasi-coherent pro-sheaf\/ $\rF$ on\/ $\fX$ is
$\kappa$\+presentable as an object of\/ $\fX\Proqc_\fl$ if and only
if, for every closed subscheme $Z\subset\fX$, the quasi-coherent
sheaf\/ $\rF^{(Z)}$ on $Z$ is locally\/ $\kappa$\+presentable
(i.~e., $\kappa$\+presentable as an object of $X\Qcoh_\fl$, or
equivalently, as an object of $X\Qcoh$).
\end{thm}

\begin{proof}
 This is a far-reaching partial generalization of the discussion of
accessibility of flat contramodules in~\cite[Corollary~7.5]{PR}
and~\cite[Theorem~10.1]{Pflcc}; cf.\
Example~\ref{ind-affine-ind-schemes-and-flat-contramodules} above.
 Let $\fX=\ilim_{\gamma\in\Gamma}X_\gamma$ be a representation of $\fX$
by a directed diagram of concentrated schemes, indexed by a directed
poset $\Gamma$ of cardinality at most~$\mu$.
 Consider the Cartesian products of categories
$\sK=\prod_{\gamma\in\Gamma}X_\gamma\Qcoh_\fl$ and
$\sL=\prod_{\gamma<\delta\in\Gamma}X_\gamma\Qcoh_\fl$.
 Here, in the construction of the category $\sL$, the Cartesian
product is taken over the set of all pairs of indices
$(\gamma,\delta)\in\Gamma^2$ such that $\gamma<\delta$.
 Propositions~\ref{product-accessible-prop}
and~\ref{flat-qcoh-sheaves-accessible} tell us that the categories
$\sK$ and $\sL$ are $\kappa$\+accessible, and provide a description
of the full subcategories of $\kappa$\+presentable objects in $\sK$
and~$\sL$.

 Consider the following pair of functors
$F$, $G\:\sK\rightrightarrows\sL$.
 To a collection of flat quasi-coherent sheaves
$(\cF_\gamma)_{\gamma\in\Gamma}\in\sK$, the functor $F$ assigns
the collection of flat quasi-coherent sheaves
$(\cH_{\gamma\delta})_{\gamma<\delta\in\Gamma}\in\sL$ with
$\cH_{\gamma\delta}=\cF_\gamma$ for all $\gamma<\delta\in\Gamma$.
 To the same collection of flat quasi-coherent sheaves
$(\cF_\gamma)_{\gamma\in\Gamma}\in\sK$, the functor $G$ assigns
the collection of flat quasi-coherent sheaves
$(\cG_{\gamma\delta})_{\gamma<\delta\in\Gamma}\in\sL$ with
$\cG_{\gamma\delta}=i_{\gamma\delta}^*\cF_\delta$
for all $\gamma<\delta\in\Gamma$.
 Here $i_{\gamma\delta}\:X_\gamma\rarrow X_\delta$ is the notation
for the transition morphism in the direct system of schemes
$(X_\gamma)_{\gamma\in\Gamma}$.

 Notice that all directed colimits exist in the categories $\sK$
and $\sL$, and all of them are preserved by the functors $F$ and~$G$.
 Furthermore, the functors $F$ and $G$ take $\kappa$\+presentable
objects of $\sK$ to $\kappa$\+presentable objects of~$\sL$.
 Take $\lambda$~to be the countable cardinal, $\lambda=\aleph_0$.
 Then the assumptions of Corollary~\ref{isomorpher-accessible-cor} are
satisfied for the pair of functors $F$, $G\:\sK\rightrightarrows\sL$.
 Thus the isomorpher category $\sD$ is $\kappa$\+accessible, and
we also obtain a description of $\kappa$\+presentable objects in~$\sD$.

 The isomorpher category $\sD$ is a rough version of the category of
flat pro-quasi-coherent pro-sheaves $\fX\Proqc_\fl$ as defined in
Section~\ref{flat-proqcoh-prosh-subsecn}, with an unimportant
difference that only the closed subschemes $Z=X_\gamma\subset\fX$,
\ $\gamma\in\Gamma$, are involved, and an important difference that
the triple compatibility condition~(iii) is not imposed.
 In order to impose the triple compatibility condition~(iii) from
Section~\ref{flat-proqcoh-prosh-subsecn}, we use the equifier
construction.

 Consider the Cartesian product of categories
$\sM=\prod_{\beta<\gamma<\delta\in\Gamma}X_\beta\Qcoh_\fl$, taken over
the set of all triples of indices $(\beta,\gamma,\delta)\in\Gamma^3$
such that $\beta<\gamma<\delta$.
 Propositions~\ref{product-accessible-prop}
and~\ref{flat-qcoh-sheaves-accessible} tell us that the category $\sM$
is $\kappa$\+accessible, and provide a description of the full
subcategory of $\kappa$\+presentable objects in~$\sM$.

 Consider the following pair of functors
$F$, $G\:\sD\rightrightarrows\sM$.
 An object of $\sD$ is a collection of flat quasi-coherent sheaves
$(\cF_\gamma)_{\gamma\in\Gamma}\in\sK$ together with a collection
of isomorphisms $\theta_{\gamma\delta}\:i_{\gamma\delta}^*\cF_\delta
\simeq\cF_\gamma$ defined for all pairs of indices
$\gamma<\delta\in\Gamma$.
 To such a datum, the functor $F$ assigns the collection of flat
quasi-coherent sheaves
$(\cH_{\beta\gamma\delta})_{\beta<\gamma<\delta\in\Gamma}\in\sM$
with $\cH_{\beta\gamma\delta}=\cF_\beta$ for all
$\beta<\gamma<\delta\in\Gamma$.
 To the same datum, the functor $G$ assigns the collection of
flat quasi-coherent sheaves
$(\cG_{\beta\gamma\delta})_{\beta<\gamma<\delta\in\Gamma}\in\sM$
with $\cG_{\beta\gamma\delta}=i_{\beta\delta}^*\cF_\delta$
for all $\beta<\gamma<\delta\in\Gamma$.

 Consider the following pair of natural transformations
$\phi$, $\psi\:F\rightrightarrows G$.
 Let $D=((\cF_\gamma)_{\gamma\in\Gamma},
(\theta_{\gamma\delta})_{\gamma<\delta\in\Gamma})$ be
an object of~$\sD$.
 For any triple of indices $\beta<\gamma<\delta\in\Gamma$,
the corresponding component of the morphism $\phi_D\:F(D)\rarrow G(D)$
is the isomorphism $\theta_{\beta\delta}^{-1}\:\cF_\beta\rarrow
i^*_{\beta\delta}\cF_\delta$.
 For the same triple of indices $\beta<\gamma<\delta\in\Gamma$,
the corresponding component of the morphism $\psi_D\:F(D)\rarrow G(D)$
is the isomorphism $\theta_{\beta\gamma}^{-1}\circ
i^*_{\beta\gamma}(\theta_{\gamma\delta}^{-1})\:\cF_\beta\rarrow
i^*_{\beta\delta}\cF_\delta$.

 Notice that all directed colimits exist in the categories $\sD$
and $\sM$, and all of them are preserved by the functors $F$ and~$G$.
 Furthermore, the functors $F$ and $G$ take $\kappa$\+presentable
objects of $\sD$ to $\kappa$\+presentable objects of~$\sM$.
 Once again, we take $\lambda$~to be the countable cardinal,
$\lambda=\aleph_0$.
 Then the assumptions of Proposition~\ref{equifier-accessible-prop} are
satisfied for the pair of functors $F$, $G\:\sD\rightrightarrows\sM$.
 Thus the equifier category $\sE$ is $\kappa$\+accessible, and
we also obtain a description of $\kappa$\+presentable objects in~$\sE$.

 Comparing the discussion in Section~\ref{flat-proqcoh-prosh-subsecn}
with the constructions above, one can immediately see that the equifier
category $\sE$ is naturally equivalent to $\fX\Proqc_\fl$.
\end{proof}

\subsection{Accessibility for complexes}
 In this section we extend the result of
Theorem~\ref{flat-proqcoh-prosheaves-accessible} to complexes of
flat pro-quasi-coherent pro-sheaves.
 The results that we obtain are a far-reaching partial generalization
of the results about accessibility of various classes of complexes of
flat contramodules from~\cite[Proposition~10.2 and
Corollary~11.4]{Pflcc}; cf.\
Example~\ref{ind-affine-ind-schemes-and-flat-contramodules} above.

\begin{thm} \label{complexes-of-flat-proqcoh-prosheaves-accessible}
 Let\/ $\mu$ be an infinite cardinal, $\kappa\ge\aleph_1$, $\kappa>\mu$
be a regular cardinal, and\/ $\fX$ be an (ind-concentrated)
$\mu$\+ind-scheme.
 Then the category\/ $\Com(\fX\Proqc_\fl)$ of complexes in
the additive category\/ $\fX\Proqc_\fl$ is\/ $\kappa$\+accessible.
 A complex of flat pro-quasi-coherent pro-sheaves\/ $\rF^\bu$ on\/ $\fX$
is\/ $\kappa$\+presentable as an object of\/ $\Com(\fX\Proqc_\fl)$
if and only if all its terms\/ $\rF^n$, \,$n\in\boZ$, are\/
$\kappa$\+presentable as objects of\/ $\fX\Proqc_\fl$.
 Equivalently, a complex\/ $\rF^\bu$ is\/ $\kappa$\+presentable in\/
$\Com(\fX\Proqc_\fl)$ if and only if, for every closed subscheme
$Z\subset\fX$, the complex of flat quasi-coherent sheaves
$\rF^{\bu\,(Z)}$ on $Z$ is\/ $\kappa$\+presentable as an object of\/
$\Com(Z\Qcoh_\fl)$.
\end{thm}

\begin{proof}
 The proof is similar to that of
Theorem~\ref{flat-proqcoh-prosheaves-accessible} and using
Proposition~\ref{complexes-of-flat-qcoh-sheaves-accessible}.
 One just needs to perform all the constructions from the proof of
Theorem~\ref{flat-proqcoh-prosheaves-accessible} working with
complexes of (flat quasi-coherent) sheaves on the schemes $X_\gamma$
instead of just sheaves.
 Basically, the only difference is that, in the context of the proof
of the present theorem, one needs to consider the Cartesian product
categories $\sK=\prod_{\gamma\in\Gamma}\Com(X_\gamma\Qcoh_\fl)$, \
$\sL=\prod_{\gamma<\delta\in\Gamma}\Com(X_\gamma\Qcoh_\fl)$, and
$\sM=\prod_{\beta<\gamma<\delta\in\Gamma}\Com(X_\beta\Qcoh_\fl)$.

 Alternatively, the assertions of the theorem can be deduced from
Theorem~\ref{flat-proqcoh-prosheaves-accessible} by applying
Proposition~\ref{category-of-complexes-accessible-prop}.
\end{proof}

\begin{thm} \label{purac-cmplxs-of-flat-proqcoh-prosheaves-accessible}
 Let\/ $\mu$ be an infinite cardinal, $\kappa\ge\aleph_1$, $\kappa>\mu$
be a regular cardinal, and\/ $\fX$ be an (ind-concentrated)
$\mu$\+ind-scheme.
 Then the category\/ $\Ac(\fX\Proqc_\fl)$ of acyclic complexes in 
the exact category\/ $\fX\Proqc_\fl$ is\/ $\kappa$\+accessible.
 A (pure) acyclic complex of flat pro-quasi-coherent pro-sheaves\/
$\rF^\bu$ on\/ $\fX$ is\/ $\kappa$\+presentable as an object of\/
$\Ac(\fX\Proqc_\fl)$ if and only if all its terms\/ $\rF^n$,
\,$n\in\boZ$, are\/ $\kappa$\+presentable as objects of\/
$\fX\Proqc_\fl$.
 Equivalently, a complex\/ $\rF^\bu$ is\/ $\kappa$\+presentable in\/
$\Ac(\fX\Proqc_\fl)$ if and only if, for every closed subscheme
$Z\subset\fX$, the pure acyclic complex of flat quasi-coherent sheaves
$\rF^{\bu\,(Z)}$ on $Z$ is\/ $\kappa$\+presentable as an object of\/
$\Ac(Z\Qcoh_\fl)$.
\end{thm}

\begin{proof}
 Similar to the proofs of
Theorems~\ref{flat-proqcoh-prosheaves-accessible}
and~\ref{complexes-of-flat-proqcoh-prosheaves-accessible}, and using
Proposition~\ref{pure-acycl-cmplxs-of-flat-qcoh-sheaves-accessible}.
 The only difference is that, in the context of the proof
of the present theorem, one needs to consider the Cartesian product
categories $\sK=\prod_{\gamma\in\Gamma}\Ac(X_\gamma\Qcoh_\fl)$, \
$\sL=\prod_{\gamma<\delta\in\Gamma}\Ac(X_\gamma\Qcoh_\fl)$, and
$\sM=\prod_{\beta<\gamma<\delta\in\Gamma}\Ac(X_\beta\Qcoh_\fl)$.
\end{proof}

\begin{rem}
 Similarly to (and rather simpler than)
Theorem~\ref{purac-cmplxs-of-flat-proqcoh-prosheaves-accessible},
one can prove that, under the same assumptions, the category of
(admissible) short exact sequences in $\fX\Proqc_\fl$ is
$\kappa$\+accessible, and its $\kappa$\+presentable objects are
the short exact sequences of $\kappa$\+presentable objects of
$\fX\Proqc_\fl$.
 (See~\cite[Corollary~10.13]{Pacc} or~\cite[Proposition~4.6]{Plce}
for a similar result for flat modules.)
 In the terminology of the paper~\cite{Plce}, this assertion can be
restated by saying that $\fX\Proqc_\fl$ is a \emph{locally\/
$\kappa$\+coherent exact category} (see~\cite[proof of
Lemma~1.2]{Plce}).
\end{rem}

\Section{Deconstructibility and the Hill Lemma}
\label{deconstructibility-secn}

\subsection{Background on deconstructibility}
 We suggest the papers~\cite[Sections~1\+-2]{SaoSt},
\cite[Section~4]{PR}, \cite[Section~7]{PS6}, \cite[Section~1]{Pfltp},
the lecture notes~\cite[Section~3]{Sto-ICRA}, and
the book~\cite[Section~B.2]{Pcosh} as reference sources on
deconstructibility in exact categories.

 Let $\sK$ be a category and $\lambda$~be an ordinal.
 A \emph{$\lambda$\+indexed chain} in $\sK$ is a functor
$\lambda\rarrow\sK$, in the sense explained in
Section~\ref{pseudolimits-subsecn}.
 In other words, a $\lambda$\+indexed chain is a directed diagram
$(K_i\to K_j)_{0\le i<j<\lambda}$ of objects and morphisms in $\sK$
indexed by the ordered set~$\lambda$.
 A $\lambda$\+indexed chain in $\sK$ is said to be \emph{smooth} if
$K_j=\varinjlim_{i<j}K_i$ in $\sK$ for all limit ordinals $0<j<\lambda$.

 Let $\sE$ be an exact category.
 A \emph{$\lambda$\+indexed filtration} $(F_i)_{0\le i\le\lambda}$
(on the object $F=F_\lambda$) in $\sE$ is a smooth
$(\lambda+\nobreak1)$\+indexed chain in $\sE$ such that $F_0=F$
and the transition morphisms $F_i\rarrow F_{i+1}$ are admissible
monomorphisms for all $0\le i<\lambda$.

 If this is the case, denote by $S_i\in\sE$, \,$0\le i<\lambda$,
the cokernels of the admissible monomorphisms $F_i\rarrow F_{i+1}$.
 Then the object $F=F_\lambda$ is said to be \emph{filtered by}
the objects $S_i$ in~$\sE$.
 In a different terminology, the object $F$ is said to be
a \emph{transfinitely iterated extension} (\emph{in the sense of
the directed colimit}) of the objects $\{S_i\}_{0\le i<\lambda}$.

 Let $\sS\subset\sE$ be a class of objects.
 The class of all objects (that can be) filtered by (objects
isomorphic to) objects from $\sS$ in $\sE$ is denoted by
$\Fil(\sS)\subset\sE$.
 A class of objects $\sF\subset\sE$ is said to be
\emph{deconstructible} if there exists a \emph{set} (rather than
a proper class) of objects $\sS\subset\sE$ such that $\sF=\Fil(\sS)$.

\subsection{Background on the Hill lemma in Grothendieck categories}
\label{background-on-hill-lemma-subsecn}
 The material of this section is an extraction from
the paper~\cite{Sto0}.

 Let $\kappa$~be a regular cardinal.
 Given a category $\sK$ with $\kappa$\+directed colimits, an object
$K\in\sK$ is said to be \emph{$\kappa$\+generated} if the functor
$\Mor_\sK(K,{-})\:\sK\rarrow\Sets$ preserves the colimits of all
the $\kappa$\+directed \emph{diagrams of monomorphisms} in~$\sK$.
 In a locally $\kappa$\+presentable abelian category $\sK$,
an object is $\kappa$\+generated if and only if it is an epimorphic
image of a $\kappa$\+presentable object~\cite[Proposition~1.69]{AR},
\cite[Lemma~A.3]{Sto0}.

 For any Grothendieck category $\sK$, there exists a regular
cardinal~$\kappa$ such that the category $\sK$ is
$\kappa$\+presentable~\cite[Lemma~A.1]{Sto0}.
 If this is the case, then the classes of $\rho$\+generated
and $\rho$\+presentable objects are closed under extensions
in $\sK$ for every regular cardinal $\rho\ge\kappa$
\,\cite[Lemma~A.4]{Sto0}.

 As the class of $\rho$\+presentable objects is preserved by
colimits of $\rho$\+small diagrams (in any category with
$\rho$\+directed colimits)~\cite[Proposition~1.16]{AR}, it follows
that the class of $\rho$\+presentable objects in $\sK$ is closed
under transfinitely iterated extensions indexed by ordinals
$\lambda<\rho$ \,\cite[Corollary~A.5]{Sto0}.

 A \emph{complete lattice} $L$ is a poset where every subset has
a least upper bound (the \emph{join}) and a greatest lower bound
(the \emph{meet}).
 A \emph{complete sublattice} in a lattice $L$ is a subset of $L$
closed under the joins and the meets of all subsets.
 A \emph{homomorphism} of complete lattices is a map preserving
all the joins and the meets.

 Given a set $\Lambda$, we denote by $P(\Lambda)=
(P(\Lambda),\cup,\cap)$ the complete lattice of all subsets
of~$\Lambda$ (ordered by inclusion).
 Given a complete, cocomplete abelian category $\sK$ and an object
$M\in\sK$, we denote by $\Sub(M)=(\Sub(M),\Sigma,\cap)$ the complete
lattice of all subobjects of~$M$ (ordered by inclusion).

\begin{thm} \label{hill-lemma-theorem}
 Let\/ $\kappa$~be a regular cardinal and\/ $\sK$ be a locally\/
$\kappa$\+presentable Grothendieck abelian category.
 Let\/ $\lambda$~be an ordinal and $M\in\sK$ be an object endowed
with a\/ $\lambda$\+indexed filtration
$(M_i\subset M)_{0\le i\le\lambda}$ such that the successive
quotient objects $S_i=M_{i+1}/M_i$ are\/ $\kappa$\+presentable
in\/ $\sK$ for all\/ $0\le i<\lambda$.
 Then there exists a complete sublattice $L\subset P(\lambda)$
and a map $\ell\:L\rarrow\Sub(M)$ with the following properties:
\begin{enumerate}
\renewcommand{\theenumi}{H\arabic{enumi}}
\item For every ordinal\/ $0\le\alpha\le\lambda$ viewed as a subset\/
$\alpha\subset\lambda$, one has\/ $\alpha\in L$ and
$\ell(\alpha)=M_\alpha$.
\item The map~$\ell$ is a homomorphism of complete lattices.
 So $\ell$~takes the unions of arbitrary families of subsets
of\/~$\lambda$ belonging to $L$ to the sums of subobjects in $M$,
and the intersections of arbitrary families of subsets of\/~$\lambda$
belonging to $L$ to the intersections of subobjects in~$M$.
\item For any two subsets $A$, $B\subset\lambda$ such that
$A$, $B\in L$ and $A\subset B$, the object $N=\ell(B)/\ell(A)\in\sK$
has a filtration $(N_j)_{0\le j\le\nu}$ indexed by some ordinal\/~$\nu$
endowed with a bijection $b\:B\setminus A\rarrow\nu$ such that
there is an isomorphism $M_{i+1}/M_i\simeq N_{b(i)+1}/N_{b(i)}$
in\/ $\sK$ for every $i\in B\setminus A$.
\item For any\/ $\kappa$\+generated subobject $K\subset M$,
there exists a subset $C\subset\lambda$, \,$C\in L$, with
the cardinality of $C$ smaller than~$\kappa$, such that
$K\subset\ell(C)\subset M$.
\end{enumerate}
\end{thm}

\begin{proof}
 This is the statement of the Hill lemma in Grothendieck
categories~\cite[Theorem~2.1]{Sto0}.
\end{proof}

 Given a class of objects $\sF$ in an additive category $\sK$,
we denote by $\sF^\oplus\subset\sK$ the class of all direct
summands of objects from $\sF$ in~$\sK$.

 Let $\kappa$~be a regular cardinal and $\sK$ be a locally
$\kappa$\+presentable Grothendieck category.
 A class of objects $\sF\subset\sK$ is said to be
\emph{$\kappa$\+deconstructible} if $\sF=\Fil(\sS)\subset\sK$
for some set of $\kappa$\+presentable objects $\sS$ in~$\sK$.

\begin{prop} \label{kaplansky-direct-summand-prop}
 Let\/ $\kappa\ge\aleph_1$ be an uncountable regular cardinal and\/
$\sK$ be a locally\/ $\kappa$\+presentable Grothendieck category.
 Then, for any\/ $\kappa$\+deconstructible class of objects\/
$\sF\subset\sK$, the class of objects\/ $\sF^\oplus\subset\sK$
is also\/ $\kappa$\+deconstructible.
\end{prop}

\begin{proof}
 This is~\cite[Proposition~2.9(1)]{Sto0}.
\end{proof}

 Given an exact category $\sK$ and a full subcategory $\sE\subset\sK$
closed under extensions in $\sK$, we always endow $\sE$ with
the inherited exact category structure.
 We also endow the category of complexes $\Com(\sK)$ with
the \emph{termwise exact structure}; so a short sequence of
complexes in $\sK$ is (admissible) exact if and only if it is
exact in every degree.
 Then the termwise exact structure on $\Com(\sE)$ is inherited from
the termwise exact structure from $\Com(\sK)$.

 For an abelian category $\sK$ endowed with the abelian exact
structure, the termwise exact structure on $\Com(\sK)$ coincides
with the abelian exact structure.

\begin{prop} \label{class-of-all-complexes-deconstructible-prop}
 Let\/ $\kappa\ge\aleph_1$ be an uncountable regular cardinal and\/
$\sK$ be a locally\/ $\kappa$\+presentable Grothendieck category.
 Let\/ $\sF\subset\sK$ be a\/ $\kappa$\+deconstructible class of
objects in\/~$\sK$ (viewed as a full subcategory in\/~$\sK$).
 Then the class of objects\/ $\Com(\sF)$ is\/ $\kappa$\+deconstructible
in the Grothendieck category\/ $\Com(\sK)$.
\end{prop}

\begin{proof}
 This is~\cite[Proposition~4.3]{Sto0}.
 It only needs to be mentioned that the Grothendieck category
$\Com(\sK)$ is locally $\kappa$\+presentable, and
the $\kappa$\+presentable objects of $\Com(\sK)$ are precisely
all the complexes of $\kappa$\+presentable objects in $\sK$,
by Lemma~\ref{category-of-complexes-locally-presentable-lem} above.
\end{proof}

\begin{prop} \label{class-of-acyclic-complexes-deconstructible-prop}
 Let\/ $\kappa\ge\aleph_1$ be an uncountable regular cardinal and\/
$\sK$ be a locally\/ $\kappa$\+presentable Grothendieck category.
 Let\/ $\sF\subset\sK$ be a\/ $\kappa$\+deconstructible class of
objects in\/~$\sK$ (viewed as a full subcategory in\/ $\sK$ with
the inherited exact structure).
 Then the class of objects\/ $\Ac(\sF)$ is\/ $\kappa$\+deconstructible
in the Grothendieck category\/ $\Com(\sK)$.
\end{prop}

\begin{proof}
 This is~\cite[Proposition~4.4]{Sto0}. 
\end{proof}

\subsection{Background on deconstructibility over schemes}
 Let $\kappa$~be a cardinal and $X$ be a scheme with an affine open
covering $X=\bigcup_\alpha U_\alpha$.
 We will say that the affine open covering $X=\bigcup_\alpha U_\alpha$
\emph{can be extended to a $\kappa$\+small affine open covering of
the intersections} if there exists an affine open covering
$X=\bigcup_{\xi\in\Xi}U_\xi$ of the scheme $X$ such that
\begin{itemize}
\item the set of indices $\{\alpha\}$ is a subset of~$\Xi$;
\item for any two indices $\upsilon'$, $\upsilon''\in\Xi$,
the intersection $U_{\upsilon'}\cap U_{\upsilon''}\subset X$
is the union of those open subschemes $U_\xi$, \,$\xi\in\Xi$,
that are contained in $U_{\upsilon'}\cap U_{\upsilon''}$;
\item the cardinality of the set of indices $\Xi$ is smaller
than~$\kappa$.
\end{itemize}

 One can easily see that any finite affine open covering of
a concentrated scheme can be extended to a countable
(i.~e., $\aleph_1$\+small) affine open covering of the intersections.
 In fact, the same applies to any countable affine open covering of
a countably quasi-compact countably quasi-separated scheme
in the sense of~\cite[Section~4]{PS6}.
 Moreover, given a regular cardinal~$\kappa$, a scheme $X$ admits
an affine open covering that can be extended to a $\kappa$\+small
affine open covering of the intersections if and only if the scheme $X$
is $\kappa$\+quasi-compact and $\kappa$\+quasi-separated
in the sense of~\cite[Remark~4.6]{PS6}.

\begin{prop} \label{deconstructibility-for-qcoh-sheaves-prop}
 Let\/ $\kappa\ge\aleph_1$ be an uncountable regular cardinal and $X$
be a scheme with an affine open covering $X=\bigcup_\alpha U_\alpha$
that can be extended to a\/ $\kappa$\+small affine open covering of
the intersections.
 Suppose that we are given, for every index\/~$\alpha$,
a\/ $\kappa$\+deconstructible class of modules\/ $\sF_\alpha\subset
\cO_X(U_\alpha)\Modl$.
 Denote by\/ $\sF_X\subset X\Qcoh$ the class of all quasi-coherent
sheaves\/ $\cF$ on $X$ such that\/ $\cF(U_\alpha)\in\sF_\alpha$
for every index\/~$\alpha$.
 Then the class\/ $\sF_X$ is\/ $\kappa$\+deconstructible in\/ $X\Qcoh$.
\end{prop}

\begin{proof}
 This is essentially~\cite[Notation~3.12 and Theorem~3.13]{EGPT}.
 To arrive at our formulation of this result from~\cite{EGPT}, one
needs to use the facts that the Grothendieck category $X\Qcoh$ is
locally $\kappa$\+presentable, and a quasi-coherent sheaf on $X$ is
$\kappa$\+presentable as an object of $X\Qcoh$ if and only if it
is locally $\kappa$\+presentable; this is
essentially~\cite[Remark~4.6]{PS6}.
 For another formulation of a comparable assertion,
see~\cite[Lemma~B.3.5]{Pcosh}.
\end{proof}

\begin{cor} \label{deconstructibility-for-cmplxes-qcoh-sheaves-cor}
 In the context and notation of
Proposition~\ref{deconstructibility-for-qcoh-sheaves-prop},
consider the Grothendieck category\/ $\Com(X\Qcoh)$ of complexes
of quasi-coherent sheaves on $X$ and the class\/
$\Com(\sF_X)\subset\Com(X\Qcoh)$ of all complexes with the terms
belonging to\/~$\sF_X$.
 Then the class\/ $\Com(\sF_X)$ is\/ $\kappa$\+deconstructible
in\/ $\Com(X\Qcoh)$.
\end{cor}

\begin{proof}
 Compare Propositions~\ref{class-of-all-complexes-deconstructible-prop}
and~\ref{deconstructibility-for-qcoh-sheaves-prop}.
\end{proof}

\begin{cor} \label{deconstructibty-for-acycl-cmplxs-qcoh-sheaves-cor}
 In the context and notation of
Proposition~\ref{deconstructibility-for-qcoh-sheaves-prop},
consider the Grothendieck category\/ $\Com(X\Qcoh)$ and the class\/
$\Ac(\sF_X)\subset\Com(X\Qcoh)$ of all acyclic complexes of
quasi-coherent sheaves on $X$ with the terms belonging to\/ $\sF_X$
and the sheaves of cocycles also belonging to\/~$\sF_X$.
 Then the class\/ $\Ac(\sF_X)$ is\/ $\kappa$\+deconstructible
in\/ $\Com(X\Qcoh)$.
\end{cor}

\begin{proof}
 Compare
Propositions~\ref{class-of-acyclic-complexes-deconstructible-prop}
and~\ref{deconstructibility-for-qcoh-sheaves-prop}.
\end{proof}

\begin{lem} \label{flat-very-flat-modules-deconstructible}
\textup{(a)} Let\/ $\kappa\ge\aleph_1$ be an uncountable regular
cardinal and $R$ be an associative ring of cardinality smaller
than\/~$\kappa$.
 Then the class of all flat left $R$\+modules $R\Modl_\fl$ is\/
$\kappa$\+deconstructible in $R\Modl$. \par
\textup{(b)} Let $R$ be a commutative ring.
 Then the class of all very flat $R$\+modules $R\Modl_\vfl$ is\/
$\aleph_1$\+deconstructible in $R\Modl$.
\end{lem}

\begin{proof}
 Part~(a) is~\cite[Lemma~1 and proof of Proposition~2]{BBE}.
 Part~(b) is~\cite[Lemma~2.2]{ST} (the assertion follows
from~\cite[Corollary~1.1.4]{Pcosh}, see the discussion in
Section~\ref{very-flat-modules-and-ring-maps-subsecn} above,
by virtue of Proposition~\ref{kaplansky-direct-summand-prop}).
\end{proof}

\begin{cor} \label{flat-and-very-flat-sheaves-deconstructible}
 Let\/ $\kappa\ge\aleph_1$ be an uncountable regular cardinal and $X$
be a scheme with an affine open covering $X=\bigcup_\alpha U_\alpha$
that can be extended to a\/ $\kappa$\+small affine open covering of
the intersections.
 In this setting: \par
\textup{(a)} Assume that the cardinalities of all the rings\/
$\cO_X(U_\alpha)$ are smaller than\/~$\kappa$.
 Then the class $X\Qcoh_\fl$ of all flat quasi-coherent sheaves on $X$
is\/ $\kappa$\+deconstructible in $X\Qcoh$. \par
\textup{(b)} The class $X\Qcoh_\vfl$ of all very flat quasi-coherent
sheaves on $X$ is\/ $\kappa$\+de\-con\-struct\-ible in $X\Qcoh$.
\end{cor}

\begin{proof}
 This was mentioned, in some form,
in~\cite[proof of Lemma~5.1.1]{Pcosh}.
 To prove the corollary, compare
Proposition~\ref{deconstructibility-for-qcoh-sheaves-prop}
with Lemma~\ref{flat-very-flat-modules-deconstructible}.
\end{proof}

\begin{cor} \label{cmplxs-of-flat-very-flat-sheaves-deconstructible}
 Let\/ $\kappa\ge\aleph_1$ be an uncountable regular cardinal and $X$
be a scheme with an affine open covering $X=\bigcup_\alpha U_\alpha$
that can be extended to a\/ $\kappa$\+small affine open covering of
the intersections.
 In this setting: \par
\textup{(a)} Assume that the cardinalities of all the rings\/
$\cO_X(U_\alpha)$ are smaller than\/~$\kappa$.
 Then the class\/ $\Com(X\Qcoh_\fl)$ of all complexes of flat
quasi-coherent sheaves on $X$ is\/ $\kappa$\+deconstructible in\/
$\Com(X\Qcoh)$. \par
\textup{(b)} The class\/ $\Com(X\Qcoh_\vfl)$ of all complexes of
very flat quasi-coherent sheaves on $X$ is\/ $\kappa$\+deconstructible
in\/ $\Com(X\Qcoh)$.
\end{cor}

\begin{proof}
 This was mentioned, in some form,
in~\cite[proof of Lemma~5.1.9]{Pcosh}.
 To prove the corollary, compare
Corollary~\ref{deconstructibility-for-cmplxes-qcoh-sheaves-cor}
with Lemma~\ref{flat-very-flat-modules-deconstructible}.
\end{proof}

\begin{cor} \label{pac-cmplxs-of-flat-vfl-sheaves-deconstructible}
 Let\/ $\kappa\ge\aleph_1$ be an uncountable regular cardinal and $X$
be a scheme with an affine open covering $X=\bigcup_\alpha U_\alpha$
that can be extended to a\/ $\kappa$\+small affine open covering of
the intersections.
 In this setting: \par
\textup{(a)} Assume that the cardinalities of all the rings\/
$\cO_X(U_\alpha)$ are smaller than\/~$\kappa$.
 Then the class\/ $\Ac(X\Qcoh_\fl)$ of all acyclic complexes of flat
quasi-coherent sheaves on $X$ with flat quasi-coherent sheaves of
cocycles is\/ $\kappa$\+deconstructible in\/ $\Com(X\Qcoh)$. \par
\textup{(b)} The class\/ $\Ac(X\Qcoh_\vfl)$ of all acyclic complexes of
very flat quasi-coherent sheaves on $X$ with very flat quasi-coherent
sheaves of cocycles is\/ $\kappa$\+deconstructible in\/ $\Com(X\Qcoh)$.
\end{cor}

\begin{proof}
 This was also mentioned, in some form,
in~\cite[proof of Lemma~5.1.1]{Pcosh}.
 To prove the corollary, compare
Corollary~\ref{deconstructibty-for-acycl-cmplxs-qcoh-sheaves-cor}
with Lemma~\ref{flat-very-flat-modules-deconstructible}.
\end{proof}

\subsection{Deconstructing flat and very flat pro-quasi-coherent
pro-sheaves}
 We start with a simple category-theoretic lemma.
 Given a regular cardinal~$\kappa$, a cocomplete abelian category
$\sK$ is said to be \emph{locally\/ $\kappa$\+generated} if $\sK$
has a set of generators consisting of $\kappa$\+generated objects.
 For example, any locally $\kappa$\+presentable abelian category is
locally $\kappa$\+generated.

 In any locally $\kappa$\+generated cocomplete abelian category,
the functors of $\kappa$\+directed colimits are exact.
 This is a version of~\cite[Corollary~9.6]{PS6} for
an arbitrary regular cardinal~$\kappa$ in place of the countable
cardinal, provable in the same way.
 Alternatively, instead of using this fact, one could simply
assume the category $\sL$ in part~(b) of the next lemma to be
Grothendieck, which is sufficient for our purposes.

\begin{lem} \label{kappa-generated-subobject-lemma}
\textup{(a)} Let\/ $\sK$ be locally\/ $\kappa$\+generated cocomplete
abelian category.
 Let $N\subset M\in\sK$ be an object with a subobject such that
the quotient object $M/N$ is\/ $\kappa$\+generated.
 Then there exists a\/ $\kappa$\+generated subobject $K\subset M$
such that $M=N+K$. \par
\textup{(b)} Let\/ $\sK$ and\/ $\sL$ be locally\/ $\kappa$\+generated
cocomplete abelian categories and $F\:\sK\rarrow\sL$ be a functor
preserving epimorphisms and\/ $\kappa$\+directed colimits.
 Let $N\subset M\in\sK$ be an object with a subobject, and let
$H\subset F(M/N)\in\sL$ be a $\kappa$\+generated subobject in $F(M/N)$.
 Then there exists a\/ $\kappa$\+generated subobject $K\subset M$
such that the image of the induced morphism $F(K)\rarrow F(M/N)$
contains~$H$.
\end{lem}
 
\begin{proof}
 Let us prove part~(b).
 By assumption, there exists a set $\Xi$ and a $\Xi$\+indexed family
of $\kappa$\+generated objects $S_\xi\in\sK$ together with
a morphism $f_\xi\:S_\xi\rarrow M$ for each $\xi\in\Xi$ such that
the induced morphism $f\:\coprod_{\xi\in\Xi}S_\xi\rarrow M$ is
an epimorphism in~$\sK$.
 Denote by $\Upsilon$ the poset of all subsets of $\Xi$ of
cardinality smaller than~$\kappa$, ordered by inclusion.
 Put $K_\upsilon=\coprod_{\xi\in\upsilon}S_\xi$ for every
$\upsilon\in\Upsilon$.
 Then $\Upsilon$ is $\kappa$\+directed poset, $K_\upsilon\in\sK$
are $\kappa$\+generated objects, and $\coprod_{\xi\in\Xi}S_\xi=
\varinjlim_{\upsilon\in\Upsilon}K_\upsilon$.

 So we have an epimorphism $f\:\varinjlim_{\upsilon\in\Upsilon}
K_\upsilon\rarrow M$ in $\sK$, and it follows that
$F(f)\:\varinjlim_{\upsilon\in\Upsilon}F(K_\upsilon)\rarrow F(M)$
is an epimorphism in~$\sL$.
 The morphism $F(M)\rarrow F(M/N)$ is an epimorphism in~$\sL$, too.
 Denote by $L_\upsilon\subset F(M/N)$ the image of the morphism
$F(K_\upsilon)\rarrow F(M/N)$.
 Then the epimorphism $\varinjlim_{\upsilon\in\Upsilon}F(K_\upsilon)
\rarrow F(M/N)$ factorizes as $\varinjlim_{\upsilon\in\Upsilon}
F(K_\upsilon)\rarrow\varinjlim_{\upsilon\in\Upsilon}L_\upsilon
\rarrow F(M/N)$.  {\hbadness=1325\par}

 By a suitable version of~\cite[Proposition~9.5]{PS6}, the morphism
$\varinjlim_{\upsilon\in\Upsilon}L_\upsilon \rarrow F(M/N)$ is
a monomorphism.
 Thus the morphism $\varinjlim_{\upsilon\in\Upsilon}L_\upsilon
\rarrow F(M/N)$ is actually an isomorphism.
 Consider the intersections of subobjects $L_\upsilon\cap H$
in $F(M/N)$.
 By suitable versions of~\cite[Proposition~9.5]{PS6}
and~\cite[Proposition~III.1.2]{Mit}, it follows that the morphism
$\varinjlim_{\upsilon\in\Upsilon}(L_\upsilon\cap H)\rarrow H$
is an isomorphism.
 Since the object $H$ is $\kappa$\+generated and the poset
$\Upsilon$ is $\kappa$\+directed, we can conclude that there exists
$\upsilon\in\Upsilon$ such that $H=L_\upsilon\cap H$.
 In other words, we have $H\subset L_\upsilon$.
 It remains to denote by $K$ the image of the morphism $K_\upsilon
\rarrow M$.

 Part~(a) is a special case of part~(b) for the categories $\sK=\sL$,
the identity functor $F=\mathrm{Id}$, and the subobject $H=M/N$.
\end{proof}

 The definition of a \emph{$\mu$\+ind-scheme} for an infinite
cardinal~$\mu$ was given in
Section~\ref{accessibility-of-flat-proqcoh-prosheaves-subsecn}.
 Given an ind-scheme $\fX=\ilim_{\gamma\in\Gamma}X_\gamma$
and a flat pro-quasi-coherent pro-sheaf $\rF$ on $\fX$, we denote
by $\rF^{(\gamma)}=\rF^{(X_\gamma)}\in X_\gamma\Qcoh_\fl$ the inverse
images of $\rF$ to~$X_\gamma$.

 Let us extend our terminology of ``$\kappa$\+deconstructible
classes'' to $\kappa$\+accessible exact categories.
 Let $\kappa$~be a regular cardinal and $\sE$ be an exact category
whose underlying (additive) category $\sE$ is $\kappa$\+accessible.
 Then we will say that a class of objects $\sF\subset\sE$ is
\emph{$\kappa$\+deconstructible} if $\sF=\Fil(\sS)$, where $\sS$
is a set of objects in $\sE$ and all the objects from $\sS$ are
$\kappa$\+presentable in\/~$\sE$.

 The following theorem is the first main result of
Section~\ref{deconstructibility-secn}.

\begin{thm} \label{deconstrbty-for-proqcoh-prosheaves}
 Let\/ $\kappa\ge\aleph_1$ be an uncountable regular cardinal and\/
$\fX=\ilim_{\gamma\in\Gamma}X_\gamma$ be an ind-scheme represented
by a directed diagram of concentrated schemes indexed by a directed
poset\/ $\Gamma$ of cardinality smaller than\/~$\kappa$.
 Suppose that we are given, for every index\/ $\gamma\in\Gamma$,
a class of flat quasi-coherent sheaves $\sF_\gamma\subset
X_\gamma\Qcoh_\fl$ on $X_\gamma$ such that the class\/ $\sF_\gamma$
is\/ $\kappa$\+deconstructible in $X_\gamma\Qcoh$ (or equivalently,
in $X_\gamma\Qcoh_\fl$).
 Denote by\/ $\sF_\fX\subset\fX\Proqc_\fl$ the class of all
flat pro-quasi-coherent pro-sheaves\/ $\rF$ on\/ $\fX$ such that\/
$\rF^{(\gamma)}\in\sF_\gamma$ for every index\/ $\gamma\in\Gamma$.
 Then the class\/ $\sF_\fX$ is $\kappa$\+deconstructible in\/
$\fX\Proqc_\fl$.
\end{thm}

\begin{proof}
 This is a far-reaching partial generalization of the discussion of
deconstructibility of classes of contramodules
in~\cite[Corollary~7.6 and Examples~7.12(2\+-3)]{PR}; cf.\
Examples~\ref{ind-affine-ind-schemes-and-flat-contramodules}
and~\ref{ind-affine-ind-schemes-and-very-flat-contramodules} above.
 We refer to Theorem~\ref{flat-proqcoh-prosheaves-accessible} for
the underlying assertion that the category $\fX\Proqc_\fl$ is
$\kappa$\+accessible, and for the description of
the $\kappa$\+presentable objects of $\fX\Proqc_\fl$.
 The proof of the theorem is based on the Hill lemma in
the Grothendieck categories $X_\gamma\Qcoh$
(Theorem~\ref{hill-lemma-theorem}).

 To begin with, recall our usual notation $i_{\gamma\delta}\:X_\gamma
\rarrow X_\delta$, \,$\gamma<\delta\in\Gamma$, for the transition
morphisms in the directed diagram of schemes
$(X_\gamma)_{\gamma\in\Gamma}$.

 Denote by $\sS_\gamma\subset X_\gamma\Qcoh_\fl$ the class (up
to isomorphism, a set) of all $\kappa$\+presentable objects in
the class $\sF_\gamma\subset X_\gamma\Qcoh_\fl$.
 Denote by $\sS_\fX$ the class (up to isomorphism, a set) of all
flat pro-quasi-coherent pro-sheaves $\rS$ on $\fX$ such that
$\rS^{(\gamma)}\in\sS_\gamma$ for all $\gamma\in\Gamma$.
 We claim that $\sF_\fX=\Fil(\sS_\fX)$.
 The inclusion $\Fil(\sS_\fX)\subset\sF_\fX$ is obvious; our task
is to prove that $\sF_\fX\subset\Fil(\sS_\fX)$.

 More specifically, let $\rF\in\sF_\fX$ be a flat pro-quasi-coherent
pro-sheaf and $\lambda$~be an ordinal such that, for every
$\gamma\in\Gamma$, we have a $\lambda$\+indexed filtration
$(G^{(\gamma)}_i\rF^{(\gamma)})_{0\le i\le\lambda}$ on $\rF^{(\gamma)}$
with the successive quotients $G^{(\gamma)}_{i+1}\rF^{(\gamma)}/
G^{(\gamma)}_i\rF^{(\gamma)}$, \ $0\le i<\lambda$, belonging
to~$\sS_\gamma$.
 We will construct a filtration
$(G^{(\fX)}_\alpha\rF)_{0\le \alpha\le\lambda}$
on the object $\rF$ of the exact category $\fX\Proqc_\fl$ with
the successive quotients $G^{(\fX)}_{\alpha+1}\rF/G^{(\fX)}_\alpha\rF$,
\ $0\le\alpha<\lambda$, belonging to~$\sS_\fX$.

 For this purpose, for each index $\gamma\in\Gamma$, we consider
the complete lattice $L_\gamma\subset P(\lambda)$ and the homomorphism
of complete lattices $\ell_\gamma\:L_\gamma\rarrow\Sub(\rF^{(\gamma)})$
provided by Theorem~\ref{hill-lemma-theorem} for the filtration
$(G^{(\gamma)}_i\rF^{(\gamma)})_{0\le i\le\lambda}$ on the object
$\rF^{(\gamma)}$ of the Grothendieck category $X_\gamma\Qcoh$.
 Recall that the category $X_\gamma\Qcoh$ is locally finitely
presentable (hence locally $\kappa$\+presentable)
by~\cite[0.5.2.5 and Corollaire~I.6.9.12]{EGA1}
or~\cite[Definition Tag~01BN and Lemma Tag~01PJ or~01PK]{SP},
since $X_\gamma$ is a concentrated scheme by assumption.
 Even more specifically, for every ordinal $0\le\alpha\le\lambda$
and every index $\gamma\in\Gamma$, we will construct a subset
$D_\alpha^\gamma\subset\lambda$, \,$D_\alpha^\gamma\in L_\gamma$,
for which we will have $(G^{(\fX)}_\alpha\rF)^{(\gamma)}=
\ell_\gamma(D_\alpha^\gamma)\subset\rF^{(\gamma)}$ and
$D_\alpha^\gamma\subset D_\beta^\gamma$ for all
$0\le\alpha\le\beta\le\lambda$.

 We proceed by transfinite induction on the ordinals
$0\le\alpha\le\lambda$.
 For $\alpha=0$, put $D_\alpha^\gamma=\varnothing$ for all
$\gamma\in\Gamma$; so $G^{(\fX)}_0\rF=0$.
 For a limit ordinal $0<\beta\le\lambda$, put
$D_\beta^\gamma=\bigcup_{\alpha<\beta}D_\alpha^\gamma$ for all
$\gamma\in\Gamma$; so $G^{(\fX)}_\beta\rF=
\varinjlim_{\alpha<\beta}G^{(\fX)}_\alpha\rF$, as required.
 A nontrivial construction is needed in the case of a successor
ordinal $\beta=\alpha+1\le\lambda$.

 For all ordinals $0\le\alpha\le\lambda$ and all indices
$\gamma\in\Gamma$, we will have $\alpha\subset D_\alpha^\gamma$
(where $\alpha$~is viewed as a subset of~$\lambda$).
 Now, it order to pass from $\alpha$ to $\beta=\alpha+1$, we
start with setting $E_0^\gamma=D_\alpha^\gamma\cup\beta$ for
all $\gamma\in\Gamma$.
 So the complement $E_0^\gamma\setminus D_\alpha^\gamma$ is
either empty or a singleton set.

 Let $\mu$ be the cardinality of~$\Gamma$ (or the countable
cardinality if $\Gamma$ is finite).
 So we have $\mu<\kappa$.
 Pick a map of sets $\psi\:\mu\rarrow\Gamma$ such that, for every
$\gamma\in\Gamma$, the subset $\psi^{-1}(\gamma)\subset\mu$
is cofinal in~$\mu$.
 For example, it suffices to make~$\psi$ the first coordinate of
a surjective map $(\psi,\phi)\:\mu\rarrow\Gamma\times\mu$.

 Proceeding by secondary transfinite induction on ordinals
$1\le j\le\mu$, for every index $\gamma\in\Gamma$ we will construct
a subset $E_j^\gamma\subset\lambda$, \,$E_j^\gamma\in L_\gamma$.
 The following properties will be preserved throughout the construction
of the sets~$E_j^\gamma$.
 For all $0\le j\le k\le\mu$ and all $\gamma\in\Gamma$, we will have
$D_\alpha^\gamma\subset E_j^\gamma\subset E_k^\gamma$, and
the cardinality of the complement $E_j^\gamma\setminus D_\alpha^\gamma$
will be smaller than~$\kappa$.

 For a limit ordinal $0<k\le\mu$, put
$E_k^\gamma=\bigcup_{j<k}E_j^\gamma$ for every $\gamma\in\Gamma$.
 Once again, a substantial construction is needed in the case of
a successor ordinal $k=j+1<\mu$.

 Put $\cH_j^\gamma=\ell_\gamma(E_j^\gamma)\subset\rF^{(\gamma)}$ for
every $\gamma\in\Gamma$, and $\gamma_j=\psi(j)$.
 By assertion~(H3) of Theorem~\ref{hill-lemma-theorem},
the quasi-coherent sheaves $\cH_j^\gamma$ and
$\rF^{(\gamma)}/\cH_j^\gamma$ on $X_\gamma$ are flat for all
$\gamma\in\Gamma$ (since so are the quasi-coherent sheaves
$G^{(\gamma)}_{i+1}\rF^{(\gamma)}/G^{(\gamma)}_i\rF^{(\gamma)}$
for all $0\le i<\lambda$).

 Furthermore, we have $(G^{(\fX)}_\alpha\rF)^{(\gamma)}\subset
\cH_j^\gamma$ (since $D_\alpha^\gamma\subset E_j^\gamma$).
 By virtue of the same assertion~(H3), we also know that
the quasi-coherent sheaf
$\cH_j^\gamma/((G^{(\fX)}_\alpha\rF)^{(\gamma)})$
on $X_\gamma$ is (locally) $\kappa$\+presentable, since 
the cardinality of $E_j^\gamma\setminus D_\alpha^\gamma$ is
smaller than~$\kappa$.

 Given an index $\gamma\in\Gamma$, we consider four cases.

 Case~I: $\gamma=\gamma_j$.
 In this case, we put $E_k^{\gamma_j}=E_j^{\gamma_j}$.

 Case~II: the elements $\gamma$ and $\gamma_j$ are incomparable
in~$\Gamma$ ($\gamma\not>\gamma_j$ and $\gamma\not<\gamma_j$).
 In this case, we also put $E_k^\gamma=E_j^\gamma$.

 Case~III: $\gamma<\gamma_j$.
 So we have the transition morphism of schemes $i_{\gamma,\gamma_j}\:
X_\gamma\rarrow X_{\gamma_j}$.
 Applying the inverse image functor $i^*_{\gamma,\gamma_j}\:
X_\gamma\Qcoh\rarrow X_{\gamma_j}\Qcoh$, we obtain a quasi-coherent
subsheaf $i_{\gamma,\gamma_j}^*\cH_j^{\gamma_j}\subset\rF^{(\gamma)}$.

 Furthermore, $G^{(\fX)}_\alpha\rF$ is a flat pro-quasi-coherent
pro-sheaf on~$\fX$ (and an admissible subobject of~$\rF$) by
the assumption of induction on~$\alpha$.
 So we have $i^*_{\gamma,\gamma_j}((G^{(\fX)}_\alpha\rF)^{(\gamma_j)})
=(G^{(\fX)}_\alpha\rF)^{(\gamma)}$.
 It follows that
$$
 \bigl(i_{\gamma,\gamma_j}^*\cH_j^{\gamma_j}\bigr)
 \big/(G^{(\fX)}_\alpha\rF)^{(\gamma)}\simeq
 i_{\gamma,\gamma_j}^*\bigl(\cH_j^{\gamma_j}/
 ((G^{(\fX)}_\alpha\rF)^{(\gamma_j)})\bigr)
$$
is a (locally) $\kappa$\+presentable quasi-coherent sheaf on~$X_\gamma$.
 In particular, the latter quasi-coherent quotient sheaf is
$\kappa$\+generated as an object of $X_\gamma\Qcoh$.

 By Lemma~\ref{kappa-generated-subobject-lemma}(a), it follows that
there exists a $\kappa$\+generated quasi-coherent subsheaf
$\cK\subset i_{\gamma,\gamma_j}^*\cH_j^{\gamma_j}$ such that
$i_{\gamma,\gamma_j}^*\cH_j^{\gamma_j}=
(G^{(\fX)}_\alpha\rF)^{(\gamma)}+\cK$.
 Applying statement~(H4) of Theorem~\ref{hill-lemma-theorem} to
the $\kappa$\+generated subobject $\cK\subset\rF^{(\gamma)}$,
we obtain a subset $C\subset\lambda$, \,$C\in L_\gamma$, with
the cardinality of $C$ smaller than~$\kappa$, such that
$\cK\subset\ell_\gamma(C)\subset\rF^{(\gamma)}$.

 Now we put $E_k^\gamma=E_j^\gamma\cup C$.
 Setting $\cH_k^\gamma=\ell_\gamma(E_k^\gamma)\subset\rF^{(\gamma)}$,
we have $i_{\gamma,\gamma_j}^*\cH_j^{\gamma_j}\subset\cH_k^\gamma$.

 Case~IV: $\gamma>\gamma_j$.
 So we have the transition morphism of schemes $i_{\gamma_j,\gamma}\:
X_{\gamma_j}\rarrow X_\gamma$.
 Once again, we have
$i^*_{\gamma_j,\gamma}((G^{(\fX)}_\alpha\rF)^{(\gamma)})
=(G^{(\fX)}_\alpha\rF)^{(\gamma_j)}$, hence
$i^*_{\gamma_j,\gamma}(\rF^{(\gamma)}/(G^{(\fX)}_\alpha\rF)^{(\gamma)})
=\rF^{(\gamma_j)}/(G^{(\fX)}_\alpha\rF)^{(\gamma_j)}$.
 The quasi-coherent sheaf
$\cH_j^{\gamma_j}/((G^{(\fX)}_\alpha\rF)^{(\gamma_j)})$
on $X_{\gamma_j}$ is (locally) $\kappa$\+presentable, hence
$\kappa$\+generated as an object of $X_{\gamma_j}\Qcoh$.

 By Lemma~\ref{kappa-generated-subobject-lemma}(b), it follows that
there exists a $\kappa$\+generated quasi-coherent subsheaf
$\cK\subset\rF^{(\gamma)}$ for which the image of the induced
morphism $i^*_{\gamma_j,\gamma}\cK\rarrow
\rF^{(\gamma_j)}/((G^{(\fX)}_\alpha\rF)^{(\gamma_j)})$
in $X_{\gamma_j}\Qcoh$ contains
$\cH_j^{\gamma_j}/((G^{(\fX)}_\alpha\rF)^{(\gamma_j)})$.
 Once again, we apply statement~(H4) of
Theorem~\ref{hill-lemma-theorem} to the $\kappa$\+generated subobject
$\cK\subset\rF^{(\gamma)}$ and obtain a subset $C\subset\lambda$,
\,$C\in L_\gamma$, with the cardinality of $C$ smaller than~$\kappa$,
such that $\cK\subset\ell_\gamma(C)\subset\rF^{(\gamma)}$.

 Once again, we put $E_k^\gamma=E_j^\gamma\cup C$.
 Setting $\cH_k^\gamma=\ell_\gamma(E_k^\gamma)\subset\rF^{(\gamma)}$,
we have $\cH_j^{\gamma_j}\subset i_{\gamma_j,\gamma}^*\cH_k^\gamma
\subset\rF^{(\gamma_j)}$.

 The construction of the subsets $E_j^\gamma\subset\lambda$,
\,$E_j^\gamma\in L_\gamma$, for all ordinals $0\le j\le\mu$ and
all indices $\gamma\in\Gamma$, is finished.
 The cardinality of the complement $E_j^\gamma\setminus
D_\alpha^\gamma$ always remains smaller than~$\kappa$, since
the union of at most~$\mu$ sets of cardinality smaller than~$\kappa$
has cardinality smaller than~$\kappa$.

 Now we put $\cF^\gamma=\ell_\gamma(E_\mu^\gamma)\subset\rF^{(\gamma)}$
for all $\gamma\in\Gamma$.
 So we have $\cF^\gamma=\varinjlim_{j<\mu}\cH_j^\gamma$.
 Let us show that, for every pair of indices $\delta<\epsilon\in\Gamma$,
we have $i^*_{\delta\epsilon}\cF^\epsilon=\cF^\delta
\subset\rF^{(\delta)}$.

 Indeed, in order to prove the inclusion
$i^*_{\delta\epsilon}\cF^\epsilon\subset\cF^\delta$, consider
the cofinal set of ordinals $0\le j<\mu$ for which $\psi(j)=\epsilon$.
 So we have $\gamma_j=\psi(j)=\epsilon$; put $k=j+1$ and
$\gamma=\delta$.
 According to Case~III above, we have
$i_{\delta\epsilon}^*\cH_j^\epsilon\subset\cH_k^\delta$.
 Passing to the directed colimit over the ordinals $0\le j<\mu$
with $\psi(j)=\epsilon$, we conclude that
$i^*_{\delta\epsilon}\cF^\epsilon\subset\cF^\delta$.

 In order to prove the converse inclusion $\cF^\delta\subset
i^*_{\delta\epsilon}\cF^\epsilon$, consider the cofinal set of ordinals
$0\le j<\mu$ for which $\psi(j)=\delta$.
 So we have $\gamma_j=\psi(j)=\delta$; put $k=j+1$ and
$\gamma=\epsilon$.
 According to Case~IV above, we have $\cH_j^\delta\subset
i_{\delta\epsilon}^*\cH_k^\epsilon$.
 Passing to the directed colimit over the ordinals $0\le j<\mu$
with $\psi(j)=\delta$, we conclude that
$\cF^\delta\subset i^*_{\delta\epsilon}\cF^\epsilon$.

 It remains to put $D_\beta^\gamma=E_\mu^\gamma$ for all
$\gamma\in\Gamma$.
 According to the argument above, the quasi-coherent subsheaves
$\cF^\gamma=\ell_\gamma(D_\beta^\gamma)\subset\rF^{(\gamma)}$,
taken together for all indices $\gamma\in\Gamma$, form an admissible
subobject $G^{(\fX)}_\beta\rF\subset\rF$ in the exact category
$\fX\Proqc_\fl$.
 So we have $(G^{(\fX)}_\beta\rF)^{(\gamma)}=\cF^\gamma$ for all
$\gamma\in\Gamma$.

 The admissible quotient object
$G^{(\fX)}_\beta\rF/G^{(\fX)}_\alpha\rF$ belongs to $\sS_\fX$,
since the cardinalities of the sets $D_\beta^\gamma\setminus
D_\alpha^\gamma$ are smaller than~$\kappa$ for all $\gamma\in\Gamma$.
 Finally, we have $\alpha\subset D_\alpha^\gamma\subset\lambda$ for
all ordinals $0\le\alpha\le\lambda$ and indices $\gamma\in\Gamma$
by construction; so $D_\lambda^\gamma=\lambda$, hence
$\ell_\gamma(D_\lambda^\gamma)=\rF^{(\gamma)}$ and
$G^{(\fX)}_\lambda\rF=\rF$.
\end{proof}

\begin{cor} \label{flat-vfl-proqcoh-prosheaves-deconstructible}
 Let\/ $\kappa\ge\aleph_1$ be an uncountable regular cardinal and\/
$\fX=\ilim_{\gamma\in\Gamma}X_\gamma$ be an ind-scheme represented
by a directed diagram of concentrated schemes indexed by a directed
poset\/ $\Gamma$ of cardinality smaller than\/~$\kappa$.
 In this setting: \par
\textup{(a)} Assume further that, for every index\/ $\gamma\in\Gamma$
and every affine open subscheme $U\subset X_\gamma$, the cardinality
of the ring\/ $\cO_{X_\gamma}(U)$ is smaller than\/~$\kappa$.
 Then the exact category of flat pro-quasi-coherent pro-sheaves\/
$\fX\Proqc_\fl$ is\/ $\kappa$\+deconstructible in itself (i.~e.,
as a class of objects in\/ $\fX\Proqc_\fl$). \par
\textup{(b)} The class of all very flat pro-quasi-coherent
pro-sheaves\/ $\fX\Proqc_\vfl$ is\/ $\kappa$\+de\-con\-struct\-ible
in the exact category\/ $\fX\Proqc_\fl$.
\end{cor}

\begin{proof}
 Compare Corollary~\ref{flat-and-very-flat-sheaves-deconstructible}
with Theorem~\ref{deconstrbty-for-proqcoh-prosheaves}.
\end{proof}

\subsection{Deconstructibility for complexes}
  The following theorem is the second main result of
Section~\ref{deconstructibility-secn}.

\begin{thm} \label{deconstrbty-for-complexes-of-proqcoh-prosheaves}
 Let\/ $\kappa\ge\aleph_1$ be an uncountable regular cardinal and\/
$\fX=\ilim_{\gamma\in\Gamma}X_\gamma$ be an ind-scheme represented
by a directed diagram of concentrated schemes indexed by a directed
poset\/ $\Gamma$ of cardinality smaller than\/~$\kappa$.
 Suppose that we are given, for every index\/ $\gamma\in\Gamma$,
a class of complexes of flat quasi-coherent sheaves $\sF_\gamma\subset
\Com(X_\gamma\Qcoh_\fl)$ on $X_\gamma$ such that the class\/
$\sF_\gamma$ is\/ $\kappa$\+deconstructible in\/ $\Com(X_\gamma\Qcoh)$,
or equivalently, in\/ $\Com(X_\gamma\Qcoh_\fl)$.
 Denote by\/ $\sF_\fX\subset\Com(\fX\Proqc_\fl)$ the class of all
complexes of flat pro-quasi-coherent pro-sheaves\/ $\rF^\bu$ on\/ $\fX$
such that\/ $\rF^{\bu\,(\gamma)}\in\sF_\gamma$ for every index\/
$\gamma\in\Gamma$.
 Then the class\/ $\sF_\fX$ is $\kappa$\+deconstructible in\/
$\Com(\fX\Proqc_\fl)$.
\end{thm}

\begin{proof}
 We refer to
Theorem~\ref{complexes-of-flat-proqcoh-prosheaves-accessible} for
the underlying assertion that the category $\Com(\fX\Proqc_\fl)$
is $\kappa$\+accessible, and for the description of
the $\kappa$\+presentable objects of $\Com(\fX\Proqc_\fl)$.
 The proof of the present theorem is completely similar to
that of Theorem~\ref{deconstrbty-for-proqcoh-prosheaves}.
 One just needs to replace all the mentions of the Grothendieck
categories $X_\gamma\Qcoh$ with the Grothendieck categories
$\Com(X_\gamma\Qcoh)$, and all the mentions of the exact subcategories
$X_\gamma\Qcoh_\fl\subset X_\gamma\Qcoh$ with the exact subcategories
$\Com(X_\gamma\Qcoh_\fl)\subset\Com(X_\gamma\Qcoh)$.
 The Grothendieck category $\Com(X_\gamma\Qcoh)$ is actually
locally finitely presentable by
Lemma~\ref{category-of-complexes-locally-presentable-lem} (since
$X_\gamma$ is a concentrated scheme by assumption), hence also
locally $\kappa$\+presentable.
\end{proof}

\begin{cor} \label{complexes-of-flat-vfl-proqcoh-prosh-deconstrble}
 Let\/ $\kappa\ge\aleph_1$ be an uncountable regular cardinal and\/
$\fX=\ilim_{\gamma\in\Gamma}X_\gamma$ be an ind-scheme represented
by a directed diagram of concentrated schemes indexed by a directed
poset\/ $\Gamma$ of cardinality smaller than\/~$\kappa$.
 In this setting: \par
\textup{(a)} Assume further that, for every index\/ $\gamma\in\Gamma$
and every affine open subscheme $U\subset X_\gamma$, the cardinality
of the ring\/ $\cO_{X_\gamma}(U)$ is smaller than\/~$\kappa$.
 Then the exact category of complexes of flat pro-quasi-coherent
pro-sheaves\/ $\Com(\fX\Proqc_\fl)$ is\/ $\kappa$\+deconstructible in
itself, i.~e., as a class of objects in\/ $\Com(\fX\Proqc_\fl)$. \par
\textup{(b)} The class of all complexes of very flat pro-quasi-coherent
pro-sheaves\/ $\Com(\fX\Proqc_\vfl)$ is\/
$\kappa$\+deconstructible in the exact category\/
$\Com(\fX\Proqc_\fl)$.  \hbadness=1900
\end{cor}

\begin{proof}
 Compare
Corollary~\ref{cmplxs-of-flat-very-flat-sheaves-deconstructible}
with Theorem~\ref{deconstrbty-for-complexes-of-proqcoh-prosheaves}.
\end{proof}

\begin{cor} \label{acycl-cmplxs-of-flat-vfl-proqcoh-prosh-deconstrble}
 Let\/ $\kappa\ge\aleph_1$ be an uncountable regular cardinal and\/
$\fX=\ilim_{\gamma\in\Gamma}X_\gamma$ be an ind-scheme represented
by a directed diagram of concentrated schemes indexed by a directed
poset\/ $\Gamma$ of cardinality smaller than\/~$\kappa$.
 In this setting: \par
\textup{(a)} Assume further that, for every index\/ $\gamma\in\Gamma$
and every affine open subscheme $U\subset X_\gamma$, the cardinality
of the ring\/ $\cO_{X_\gamma}(U)$ is smaller than\/~$\kappa$.
 Then the class\/ $\Ac(\fX\Proqc_\fl)\subset\Com(\fX\Proqc_\fl)$ of all
acyclic complexes in the exact category\/ $\fX\Proqc_\fl$ of flat
pro-quasi-coherent pro-sheaves on\/ $\fX$ is\/ $\kappa$\+deconstructible
in the exact category\/ $\Com(\fX\Proqc_\fl)$. \par
\textup{(b)} The class\/ $\Ac(\fX\Proqc_\vfl)\subset\Com(\fX\Proqc_\fl)$
of all acyclic complexes in the exact category\/ $\fX\Proqc_\vfl$ of
very flat pro-quasi-coherent pro-sheaves on\/ $\fX$ is\/
$\kappa$\+deconstructible in the exact category\/ $\Com(\fX\Proqc_\fl)$.
\end{cor}

\begin{proof}
 Compare
Corollary~\ref{pac-cmplxs-of-flat-vfl-sheaves-deconstructible}
with Theorem~\ref{deconstrbty-for-complexes-of-proqcoh-prosheaves}.
\end{proof}

\Section{The $\Ext^1$-Orthogonality Argument on Ind-Schemes}
\label{orthogonality-secn}

 In this section we recall and generalize an argument
from~\cite[Section~8.5]{Pform} proving $\Ext^1$\+orthogonality of
some flat pro-quasi-coherent pro-sheaves $\rF$ and $\rC$ on
a reasonable $\aleph_0$\+ind-scheme $\fX$ under certain conditions on
the inverse images of $\rF$ and $\rC$ to closed subschemes of~$\fX$.

\subsection{Background on the $\Ext$ notation}
 As a general rule, we denote the abelian groups of morphisms in
an additive category $\sK$ by $\Hom_\sK({-},{-})$ and the Yoneda
$\Ext$ groups in an exact category $\sE$ by $\Ext^*_\sE({-},{-})$.

 Given a scheme $X$, we denote by $\Hom_X({-},{-})$ the abelian
groups of morphisms in the abelian category of quasi-coherent
sheaves $X\Qcoh$.
 Given an ind-scheme $\fX$, we denote by $\Hom_{(\fX)}({-},{-})$
the abelian groups of morphisms in the exact category of flat
pro-quasi-coherent pro-sheaves $X\Proqc_\fl$.
 
 Furthermore, we denote by $\Ext^*_X({-},{-})$ the $\Ext$ groups
computed in the abelian category $X\Qcoh$.
 We denote by $\Ext^*_{(\fX)}({-},{-})$ the $\Ext$ groups computed
in the exact category $\fX\Proqc_\fl$.

 Notice that, generally speaking, there is \emph{no reason} for
the $\Ext$ groups computed in the exact category of flat
quasi-coherent sheaves $X\Qcoh_\fl$ to agree with the $\Ext$ groups
in the abelian category $X\Qcoh$.
 However, the groups ($\Hom=\Ext^0$ and) $\Ext^1$ computed in
$X\Qcoh_\fl$ and in $X\Qcoh$ always agree, as $X\Qcoh_\fl$ is
a full subcategory closed under extensions in $X\Qcoh$.

 The full subcategory $X\Qcoh_\fl$ is also always closed under
kernels of epimorphisms in $X\Qcoh$.
 On a quasi-compact semi-separated scheme $X$, any quasi-coherent
sheaf is a quotient sheaf of a flat quasi-coherent
sheaf~\cite[Section~2.4]{M-n}, \cite[Section~3.2]{M-th},
\cite[Lemma~A.1]{EP}, \cite[Lemma~4.1.1 or~4.1.8]{Pcosh}.
 It follows that all the groups $\Ext^*$ computed in $X\Qcoh_\fl$
agree with the ones computed in $X\Qcoh$ in this
case (see~\cite[Theorem~12.1(b)]{Kel}, \cite[Proposition~13.2.2(i)]{KS},
\cite[Proposition~A.2.1 or~A.3.1(a)]{Pcosh},
or~\cite[Lemma~6.3]{Pfltp}).

\subsection{The general formulation of orthogonality}
 We start with restating a lemma from~\cite[Section~8.5]{Pform}.

\begin{lem} \label{tensor-with-sheaf-of-ideals-adjunction-sequence}
 Let $i\:Z\rarrow X$ be a closed immersion of schemes, and let\/
$\cI\subset\cO_X$ be the corresponding quasi-coherent sheaf of
ideals on~$X$.
 Then, for any flat quasi-coherent sheaf\/ $\cF$ on $X$, there is
a natural short exact sequence of quasi-coherent sheaves on~$X$
$$
 0\lrarrow\cI\ot_{\cO_X}\cF\lrarrow\cF\lrarrow i_*i^*\cF\lrarrow0,
$$
where $\cF\rarrow i_*i^*\cF$ is the adjunction unit for
the pair of adjoint functors of inverse and direct image
$i^*\:X\Qcoh\rarrow Z\Qcoh$ and $i_*\:Z\Qcoh\rarrow X\Qcoh$.
\end{lem}

\begin{proof}
 This is~\cite[Lemma~8.5.5]{Pform}.
\end{proof}

 Recall the shorthand notation $\rG^{(\gamma)}=\rG^{(X_\gamma)}\in
X_\gamma\Qcoh_\fl$ for a flat pro-quasi-coherent pro-sheaf $\rG$
on an ind-scheme $\fX=\ilim_{\gamma\in\Gamma}X_\gamma$.

\begin{thm} \label{flat-proqcoh-prosheaves-Ext1-orthogonality-thm}
 Let\/ $\fX=\ilim_{\gamma\in\Gamma}X_\gamma$ be a strict\/
$\aleph_0$\+ind-scheme represented by a directed diagram of closed
immersions of concentrated schemes indexed by the ordered set of
positive integers\/ $\Gamma=\boZ_{\ge1}$.
 Let\/ $\rF$ and\/ $\rC$ be flat pro-quasi-coherent pro-sheaves
on\/~$\fX$.
 Assume that
\begin{enumerate}
\renewcommand{\theenumi}{\roman{enumi}}
\item $\Ext^1_{X_\gamma}(\rF^{(\gamma)},\rC^{(\gamma)})=0$
for all\/ $\gamma\in\boZ_{\ge1}$;
\item $\Ext^1_{X_\gamma}(\rF^{(\gamma)},\>
\cI_\gamma\ot_{\cO_{X_\gamma}}\rC^{(\gamma)})=0$
for all\/ $\gamma\in\boZ_{\ge2}$, where\/
$\cI_{\gamma}\subset\cO_{X_\gamma}$ denotes the quasi-coherent sheaf
of ideals of the closed subscheme $X_{\gamma-1}\subset X_\gamma$
for\/ $\gamma\ge2$.
\end{enumerate}
 Then\/ $\Ext^1_{(\fX)}(\rF,\rC)=0$.
\end{thm}

\begin{proof}
 This is a generalization of~\cite[Lemma~8.5.7]{Pform}, and the proof
is similar.
 Denote by $i_\gamma\:X_\gamma\rarrow X_{\gamma+1}$,
\,$\gamma\in\boZ_{\ge1}$, the transition morphisms of schemes in
the directed diagram $(X_\gamma)_{\gamma\in\Gamma}$.
 We need to show that any (admissible) short exact sequence
$0\rarrow\rC\rarrow\rG\rarrow\rF\rarrow0$ in $\fX\Proqc_\fl$ splits.
 Denote by $j\:\rC\rarrow\rG$ the admissible monomorphism of flat
pro-quasi-coherent pro-sheaves in our short exact sequence and by
$j^{(\gamma)}\:\rC^{(\gamma)}\rarrow\rG^{(\gamma)}$ its inverse
images to $X_\gamma\subset\fX$.

 For any index $\gamma\in\boZ_{\ge1}$, the short exact sequence
$0\rarrow\rC^{(\gamma)}\rarrow\rG^{(\gamma)}\rarrow\rF^{(\gamma)}
\rarrow0$ in $X_\gamma\Qcoh_\fl$ splits by condition~(i).
 So we have a splitting morphism $s_\gamma\:\rF^{(\gamma)}\rarrow
\rG^{(\gamma)}$ in $X_\gamma\Qcoh_\fl$.
 To prove the theorem, we need to show that the splitting
morphisms~$s_\gamma$ can be chosen in such a way that they are
compatible as $\gamma$~varies, in the sense that
$s_\gamma=i_\gamma^*(s_{\gamma+1})$ for all $\gamma\in\boZ_{\ge1}$.

 Proceeding by induction on~$\gamma$, we can assume that we have
already chosen a splitting morphism~$s_\gamma$, and we need to find
a splitting morphism~$s_{\gamma+1}$ satisfying the equation
$s_\gamma=i_\gamma^*(s_{\gamma+1})$.
 Pick an arbitrary splitting morphism $s'_{\gamma+1}\:
\rF^{(\gamma+1)}\rarrow\rG^{(\gamma+1)}$ for the split short exact
sequence $0\rarrow\rC^{(\gamma+1)}\rarrow\rG^{(\gamma+1)}\rarrow
\rF^{(\gamma+1)}\rarrow0$ in $X_{\gamma+1}\Qcoh_\fl$.
 Then the difference $s_\gamma-i_\gamma^*(s'_{\gamma+1})\:\rF^{(\gamma)}
\rarrow\rG^{(\gamma)}$ factorizes through the admissible monomorphism
$\rC^{(\gamma)}\rarrow\rG^{(\gamma)}$, so we obtain a morphism
$t_\gamma\:\rF^{(\gamma)}\rarrow\rC^{(\gamma)}$ in $X_\gamma\Qcoh_\fl$.
 Specifically, we have $s_\gamma-i_\gamma^*(s'_{\gamma+1})=
j^{(\gamma)}\circ t_\gamma$.

 All we need to do is to lift the morphism~$t_\gamma$ to
$X_{\gamma+1}\Qcoh_\fl$, i.~e., find a morphism
$t_{\gamma+1}\:\rF^{(\gamma+1)}\rarrow\rC^{(\gamma+1)}$ such
that $t_\gamma=i_\gamma^*(t_{\gamma+1})$.
 After this is done, the desired splitting morphism
$s_{\gamma+1}\:\rF^{(\gamma+1)}\rarrow\rG^{(\gamma+1)}$ compatible
with~$s_\gamma$ will be constructed as the sum ot two morphisms
$s_{\gamma+1}=s'_{\gamma+1}+j^{(\gamma+1)}\circ t_{\gamma+1}$.

 The morphism $t_\gamma\:i_\gamma^*\rF^{(\gamma+1)}\rarrow
i_\gamma^*\rC^{(\gamma+1)}$ in $X_\gamma\Qcoh_\fl$ corresponds by
adjunction to a morphism $\tilde t_\gamma\:\rF^{(\gamma+1)}\rarrow
i_\gamma{}_*i_\gamma^*\rC^{(\gamma+1)}$ in $X_{\gamma+1}\Qcoh_\fl$.
 Lifing~$t_\gamma$ to a morphism $t_{\gamma+1}\:\rF^{(\gamma+1)}
\rarrow\rC^{(\gamma+1)}$ such that $t_\gamma=i_\gamma^*(t_{\gamma+1})$
is equivalent to lifting~$\tilde t_\gamma$ to a morphism
$t_{\gamma+1}\:\rF^{(\gamma+1)}\rarrow\rC^{(\gamma+1)}$ such that
the composition of~$t_{\gamma+1}$ with the adjunction morphism
$\rC^{(\gamma+1)}\rarrow i_\gamma{}_*i_\gamma^*\rC^{(\gamma+1)}$
is equal to~$\tilde t_\gamma$.

 Finally, by
Lemma~\ref{tensor-with-sheaf-of-ideals-adjunction-sequence} we have
a short exact sequence $0\rarrow\cI_{\gamma+1}
\ot_{\cO_{X_{\gamma+1}}}\rC^{(\gamma+1)}\rarrow\rC^{(\gamma+1)}
\rarrow i_\gamma{}_*i_\gamma^*\rC^{(\gamma+1)}\rarrow0$ in
$X_{\gamma+1}\Qcoh$.
 By condition~(ii), the group
$\Ext^1_{X_{\gamma+1}}(\rF^{(\gamma+1)},\>
\cI_{\gamma+1}\ot_{\cO_{X_{\gamma+1}}}\rC^{(\gamma+1)})$ vanishes.
 Hence any morphism $\rF^{(\gamma+1)}\rarrow
i_\gamma{}_*i_\gamma^*\rC^{(\gamma+1)}$ in the abelian category
$X_{\gamma+1}\Qcoh$ can be lifted to a morphism
$\rF^{(\gamma+1)}\rarrow\rC^{(\gamma+1)}$, and we are done.
\end{proof}

\subsection{Background on contraadjusted and cotorsion
quasi-coherent sheaves}
 This section is an extraction from~\cite[Section~4.1]{Pcosh}.

 Let $X$ be a scheme.
 A quasi-coherent sheaf $\cC$ on $X$ is said to be
\emph{cotorsion}~\cite{EE} if $\Ext^1_X(\cF,\cC)=0$ for all flat
quasi-coherent sheaves $\cF$ on~$X$.
 Similarly, a quasi-coherent sheaf $\cC$ on $X$ is said to be
\emph{contraadjusted}~\cite[Sections~2.5 and~4.1]{Pcosh} if
$\Ext^1_X(\cF,\cC)=0$ for all very flat quasi-coherent sheaves $\cF$
on~$X$.

 In particular, a quasi-coherent sheaf $\cC$ on an affine scheme $U$
is cotorsion if and only if the $\cO(U)$\+module $\cC(U)$ is
cotorsion.
 A quasi-coherent sheaf $\cC$ on $U$ is contraadjusted if and only
if the $\cO(U)$\+module $\cC(U)$ is contraadjusted.
 (See Section~\ref{very-flat-modules-and-ring-maps-subsecn}
for the module-theoretic background.)

 The classes of contraadjusted and cotorsion quasi-coherent sheaves
on any scheme $X$ are obviously closed under extensions and direct
summands.
 Both the classes are preserved by the direct images with respect to
quasi-compact quasi-separated morphisms of
schemes~\cite[Section~2.5]{Pcosh}.
 We denote the full subcategories of cotorsion and contraadjusted
quasi-coherent sheaves by $X\Qcoh^\cot\subset X\Qcoh^\cta\subset
X\Qcoh$.

 A more substantial theory of contraadjusted and cotorsion
quasi-coherent sheaves can be developed for quasi-compact
semi-separated schemes~$X$ \,\cite[Section~4.1]{Pcosh}.
 The key fact is that, on such a scheme $X$, any quasi-coherent
sheaf $\cM$ is a quotient sheaf of a flat quasi-coherent sheaf~$\cF$
\,\cite[Section~2.4]{M-n}, \cite[Section~3.2]{M-th}.
 Actually, one can even choose $\cF$ to be a very flat quasi-coherent
sheaf~\cite[Lemma~4.1.1]{Pcosh}.
 The explicit construction of such a flat/very flat quasi-coherent
sheaf $\cF$ for a given quasi-coherent sheaf $\cM$ suggested
in~\cite[Lemma~A.1]{EP} and further developed in~\cite[Lemmas~4.1.1
and~4.1.8]{Pcosh} plays a crucial role.

\begin{prop} \label{flat-cotorsion-pair-in-qcoh-sheaves}
 Let $X$ be a quasi-compact semi-separated scheme and\/ $\cM$ be
a quasi-coherent sheaf on~$X$.
 Then \par
\textup{(a)} there exists a short exact sequence\/ $0\rarrow\cC\rarrow
\cF\rarrow\cM\rarrow0$ in $X\Qcoh$ with a flat quasi-coherent sheaf\/
$\cF$ and a cotorsion quasi-coherent sheaf\/~$\cC$; \par
\textup{(b)} there exists a short exact sequence\/ $0\rarrow\cM\rarrow
\cC\rarrow\cF\rarrow0$ in $X\Qcoh$ with a cotorsion quasi-coherent
sheaf\/ $\cC$ and a flat quasi-coherent sheaf\/~$\cF$.
\end{prop}

\begin{proof}
 Part~(b) is ``almost'' a special case of~\cite[Corollary~4.2]{EE}
in view of the \emph{Wakamatsu lemma}~\cite[Lemma~2.1.2]{Xu}.
 The facts that every object of the Grothendieck category $X\Qcoh$
is a subobject of an injective object, and any injective
quasi-coherent sheaf is cotorsion, implying that any cotorsion
(pre)envelope is a monomorphism, are used here.
 With this approach, one also needs to use the fact that the vanishing
of $\Ext^1_X(\cF,\cC)$ for all cotorsion quasi-coherent sheaves $\cC$
implies flatness of a quasi-coherent sheaf~$\cF$.
 The latter fact is provable for quasi-compact semi-separated
schemes $X$ using existence of a flat generator in $X\Qcoh$.

 Part~(a) follows from part~(b) and the above-mentioned result
from~\cite[Section~2.4]{M-n}, \cite[Section~3.2]{M-th},
or~\cite[Lemma~A.1]{EP} by virtue of the \emph{Salce
lemma}~\cite{Sal}, \cite[Lemma~B.1.1(a)]{Pcosh}.

 Alternatively, our preferred proof can be found
in~\cite[Section~4.1]{Pcosh}; see~\cite[Corollary~4.1.11(a\+-b)]{Pcosh}.
\end{proof}

\begin{prop} \label{very-flat-cotorsion-pair-in-qcoh-sheaves}
 Let $X$ be a quasi-compact semi-separated scheme and\/ $\cM$ be
a quasi-coherent sheaf on~$X$.
 Then \par
\textup{(a)} there exists a short exact sequence\/ $0\rarrow\cC\rarrow
\cF\rarrow\cM\rarrow0$ in $X\Qcoh$ with a very flat quasi-coherent
sheaf\/ $\cF$ and a contraadjusted quasi-coherent sheaf\/~$\cC$; \par
\textup{(b)} there exists a short exact sequence\/ $0\rarrow\cM\rarrow
\cC\rarrow\cF\rarrow0$ in $X\Qcoh$ with a contraadjusted quasi-coherent
sheaf\/ $\cC$ and a very flat quasi-coherent sheaf\/~$\cF$.
\end{prop}

\begin{proof}
 This is~\cite[Corollary~4.1.4(a\+-b)]{Pcosh}.
\end{proof}

\begin{cor} \label{cotorsion-qcoh-sheaves-hereditary-cor}
 Let $X$ be a quasi-compact semi-separated scheme. \par
\textup{(a)} Let $\cF$ be a flat quasi-coherent sheaf on $X$ and\/
$\cC$ be a cotorsion quasi-coherent sheaf on~$X$.
 Then\/ $\Ext_X^n(\cF,\cC)=0$ for all $n\ge1$. \par
\textup{(b)} The full subcategory $X\Qcoh^\cot$ of cotorsion
quasi-coherent sheaves is closed under cokernels of monomorphisms
in $X\Qcoh$.
\end{cor}

\begin{proof}
 This is~\cite[Corollary~4.1.9]{Pcosh} or~\cite[Corollary~7.5.6]{Pform}.
 The assertions follow from the facts that any quasi-coherent sheaf
on $X$ is a quotient sheaf of a flat quasi-coherent sheaf and
the class of flat quasi-coherent sheaves is closed under kernels of
epimorphisms in $X\Qcoh$.
\end{proof}

\begin{cor} \label{cta-qcoh-sheaves-hereditary-cor}
 Let $X$ be a quasi-compact semi-separated scheme. \par
\textup{(a)} Let $\cF$ be a very flat quasi-coherent sheaf on $X$ and\/
$\cC$ be a contraadjusted quasi-coherent sheaf on~$X$.
 Then\/ $\Ext_X^n(\cF,\cC)=0$ for all $n\ge1$. \par
\textup{(b)} The full subcategory $X\Qcoh^\cta$ of contraadjusted
quasi-coherent sheaves is closed under cokernels of monomorphisms
in $X\Qcoh$.
\end{cor}

\begin{proof}
 This is~\cite[Corollary~4.1.2]{Pcosh} or~\cite[Corollary~7.5.3]{Pform}.
 The assertions follow from the facts that any quasi-coherent sheaf
on $X$ is a quotient sheaf of a very flat quasi-coherent sheaf (a weak
form of Proposition~\ref{very-flat-cotorsion-pair-in-qcoh-sheaves}(a))
and the class of very flat quasi-coherent sheaves is closed under
kernels of epimorphisms in $X\Qcoh$.
\end{proof}

 The following two propositions spell out \emph{antilocality}
properties of various classes of quasi-coherent sheaves,
in the sense similar to~\cite[Section~4]{Pal}.

\begin{prop} \label{cta-and-cot-quasi-coherent-sheaves-prop}
 Let $X$ be a quasi-compact semi-separated scheme with a finite affine
open covering $X=\bigcup_{\alpha=1}^N U_\alpha$.
 In this context: \par
\textup{(a)} a quasi-coherent sheaf\/ $\cC$ on $X$ is cotorsion if
and only if\/ $\cC$ is a direct summand of a finitely iterated
extension of the direct images of cotorsion quasi-coherent sheaves
from~$U_\alpha$; \par
\textup{(b)} a quasi-coherent sheaf\/ $\cC$ on $X$ is contraadjusted
if and only if\/ $\cC$ is a direct summand of a finitely iterated
extension of the direct images of contraadjusted quasi-coherent sheaves
from~$U_\alpha$.
\end{prop}

\begin{proof}
 Part~(a) is~\cite[Corollary~4.1.11(c)]{Pcosh}.
 Part~(b) is~\cite[Corollary~4.1.4(c)]{Pcosh}.
\end{proof}

\begin{prop} \label{flat-vfl-cta-cot-quasi-coherent-sheaves-prop}
 Let $X$ be a quasi-compact semi-separated scheme with a finite affine
open covering $X=\bigcup_{\alpha=1}^N U_\alpha$.
 In this context: \par
\textup{(a)} a quasi-coherent sheaf\/ $\cC$ on $X$ is flat and
cotorsion if and only if\/ $\cC$ is a direct summand of a finite
direct sum of the direct images of flat cotorsion quasi-coherent
sheaves from~$U_\alpha$; \par
\textup{(b)} a quasi-coherent sheaf\/ $\cC$ on $X$ is flat and
contraadjusted if and only if\/ $\cC$ is a direct summand of a finitely
iterated extension of the direct images of flat contraadjusted
quasi-coherent sheaves from~$U_\alpha$; \par
\textup{(c)} a quasi-coherent sheaf\/ $\cC$ on $X$ is very flat and
contraadjusted if and only if\/ $\cC$ is a direct summand of a finite
direct sum the direct images of very flat contraadjusted quasi-coherent
sheaves from~$U_\alpha$.
\end{prop}

\begin{proof}
 Part~(a) is~\cite[Lemma~4.1.12]{Pcosh}.
 Part~(b) is~\cite[Lemma~4.1.6]{Pcosh}.
 Part~(c) is~\cite[Lemma~4.1.5]{Pcosh}.
\end{proof}

\subsection{Pro-sheaves with flat contraadjusted and flat cotorsion
components}
 We start with two lemmas essentially from~\cite[Section~8.5]{Pform}.

\begin{lem} \label{closed-subscheme-flat-cta-cot-lemma}
 Let $i\:Z\rarrow X$ be a closed immersion of quasi-compact
semi-separated schemes and $i^*\:X\Qcoh\rarrow Z\Qcoh$
be the inverse image functor.  In this context: \par
\textup{(a)} the functor~$i^*$ takes flat contraadjusted quasi-coherent
sheaves on $X$ to flat contraadjusted quasi-coherent sheaves on~$Z$;
\par
\textup{(b)} if the schemes $X$ and $Z$ are Noetherian, then
the functor~$i^*$ takes flat cotorsion quasi-coherent sheaves on $X$
to flat cotorsion quasi-coherent sheaves on~$Z$.
\end{lem}

\begin{proof}
 Part~(a) is~\cite[Corollary~8.5.4(a)]{Pform}.
 Part~(b) is a particular case of~\cite[Corollary~8.5.4(b)]{Pform}.
 The arguments are based on the results of
Proposition~\ref{flat-vfl-cta-cot-quasi-coherent-sheaves-prop}(a\+-b)
above and~\cite[Lemmas~1.7.7(b) and~1.7.9(a)]{Pcosh}.
\end{proof}

 The definition of a (\emph{locally}) \emph{finitely generated}
quasi-coherent sheaf on a concentrated scheme was given
in Section~\ref{reasonable-ind-schemes-subsecn}.

\begin{lem} \label{tensor-product-flat-cta-cot-lemma}
 Let $X$ be a quasi-compact semi-separated scheme and\/ $\cM$ be
a locally finitely generated quasi-coherent sheaf on~$X$.
 In this context: \par
\textup{(a)} for any flat contraadjusted quasi-coherent sheaf\/ $\cC$
on $X$, the quasi-coherent sheaf\/ $\cM\ot_{\cO_X}\cC$ on $X$ is
contraadjusted; \par
\textup{(b)} if the scheme $X$ is Noetherian then, for any flat
cotorsion quasi-coherent sheaf\/ $\cC$ on $X$, the quasi-coherent
sheaf\/ $\cM\ot_{\cO_X}\cC$ on $X$ is cotorsion.
\end{lem}

\begin{proof}
 Part~(b) is~\cite[Corollary~8.5.3(b)]{Pform}.
 Part~(a) is a straightforward generalization
of~\cite[Corollary~8.5.3(a)]{Pform}, provable in the same way.
 The arguments are based on the results of
Proposition~\ref{flat-vfl-cta-cot-quasi-coherent-sheaves-prop}(a\+-b)
above and~\cite[Corollaries~1.7.4(a) and~1.7.5(a)]{Pcosh}.
\end{proof}

 Let $\fX$ be a (strict ind-quasi-compact) ind-semi-separated
ind-scheme.
 We will say that a flat pro-quasi-coherent pro-sheaf $\rC$ on $\fX$
is \emph{levelwise contraadjusted} if the flat quasi-coherent sheaf
$\rC^{(Z)}$ on $Z$ is contraadjusted for every closed subscheme
$Z\subset\fX$.

 Given a directed diagram of closed immersions of quasi-compact
semi-separated schemes $(X_\gamma)_{\gamma\in\Gamma}$ representing
$\fX$, it suffices to check the contraadjustedness condition from
the previous paragraph for the closed subschemes $Z=X_\gamma\subset
\fX$, \,$\gamma\in\Gamma$.
 This follows from Lemma~\ref{closed-subscheme-flat-cta-cot-lemma}(a).

 Let $\fX$ be a (strict) ind-semi-separated ind-Noetherian ind-scheme.
 We will say that a flat pro-quasi-coherent pro-sheaf $\rC$ on $\fX$
is \emph{levelwise cotorsion} if the flat quasi-coherent sheaf
$\rC^{(Z)}$ on $Z$ is cotorsion for every closed subscheme
$Z\subset\fX$.

 Given a directed diagram of closed immersions of semi-separated
Noetherian schemes $(X_\gamma)_{\gamma\in\Gamma}$ representing
$\fX$, it suffices to check the cotorsion condition from the previous
paragraph for the closed subschemes $Z=X_\gamma\subset\fX$,
\,$\gamma\in\Gamma$.
 This follows from Lemma~\ref{closed-subscheme-flat-cta-cot-lemma}(b).

 The definition of a \emph{reasonable closed immersion} of
concentrated schemes, as well as that of a \emph{reasonable ind-scheme},
was given in Section~\ref{reasonable-ind-schemes-subsecn}.

\begin{cor} \label{levelwise-cta-cot-orthogonality-cor}
\textup{(a)} Let\/ $\fX$ be a (strict ind-quasi-compact)
ind-semi-separated reasonable\/ $\aleph_0$\+ind-scheme.
 Then one has\/ $\Ext_{(\fX)}^1(\rF,\rC)=0$ for any very flat
pro-quasi-coherent pro-sheaf\/ $\rF$ and any levelwise 
contraadjusted flat pro-quasi-coherent pro-sheaf\/ $\rC$ on\/~$\fX$.
 In particular, any levelwise contraadjusted very flat
pro-quasi-coherent pro-sheaf on\/ $\fX$ is an injective object of
the exact category\/ $\fX\Proqc_\vfl$.
\par
\textup{(b)} Let\/ $\fX$ be a (strict) ind-semi-separated
ind-Noetherian\/ $\aleph_0$\+ind-scheme.
 Then one has\/ $\Ext_{(\fX)}^1(\rF,\rC)=0$ for any flat
pro-quasi-coherent pro-sheaf\/ $\rF$ and any levelwise cotorsion
flat pro-quasi-coherent pro-sheaf\/ $\rC$ on\/~$\fX$.
 In other words, any levelwise cotorsion flat pro-quasi-coherent
pro-sheaf on\/ $\fX$ is an injective object of the exact category\/
$\fX\Proqc_\fl$.
\end{cor}

\begin{proof}
 This is a generalization of~\cite[Lemma~8.5.7]{Pform} from
(adic Noetherian) formal schemes to ind-schemes.
 This is also a far-reaching partial generalization of
the contramodule-theoretic $\Ext^1$\+orthogonality results
of~\cite[Lemma~E.4.3 and Corollary~E.5.2]{Pcosh}; cf.\
Examples~\ref{ind-affine-ind-schemes-and-flat-contramodules}
and~\ref{ind-affine-ind-schemes-and-very-flat-contramodules} above.

 Part~(a): by Corollary~\ref{Ind-of-Sch-strict-qs-qc-ss-lemma}(b)
and the subsequent paragraph, together
with Lemma~\ref{strict-aleph0-ind-schemes-lemma}(b),
the ind-scheme $\fX$ can be represented by a countable directed diagram
of reasonable closed immersions of quasi-compact semi-separated schemes.
 Now the assertion follows from
Theorem~\ref{flat-proqcoh-prosheaves-Ext1-orthogonality-thm},
which is applicable in view of
Lemma~\ref{tensor-product-flat-cta-cot-lemma}(a).

 Part~(b): by Corollary~\ref{Ind-of-Sch-strict-qs-qc-ss-lemma}(b\+-c)
and the subsequent paragraph, together
with Lemma~\ref{strict-aleph0-ind-schemes-lemma}(b),
the ind-scheme $\fX$ can be represented by a countable directed diagram
of closed immersions of semi-separated Noetherian schemes.
 Now the assertion follows from
Theorem~\ref{flat-proqcoh-prosheaves-Ext1-orthogonality-thm},
which is applicable in view of
Lemma~\ref{tensor-product-flat-cta-cot-lemma}(b).
\end{proof}

\subsection{Orthogonality for complexes}
 The following theorem is a straightforward extension of
Theorem~\ref{flat-proqcoh-prosheaves-Ext1-orthogonality-thm} from
(flat pro-quasi-coherent) pro-sheaves to complexes of pro-sheaves.
 The termwise exact structure on complexes, as defined at the end of
Section~\ref{background-on-hill-lemma-subsecn}, is presumed in
the $\Ext$ notation here and below.

\begin{thm} \label{cmplxes-offlat-proqcoh-prosh-Ext1-orthogty-thm}
 Let\/ $\fX=\ilim_{\gamma\in\Gamma}X_\gamma$ be a strict\/
$\aleph_0$\+ind-scheme represented by a directed diagram of closed
immersions of concentrated schemes indexed by the ordered set of
positive integers\/ $\Gamma=\boZ_{\ge1}$.
 Let\/ $\rF^\bu$ and\/ $\rC^\bu$ be two complexes of flat
pro-quasi-coherent pro-sheaves on\/~$\fX$.
 Assume that
\begin{enumerate}
\renewcommand{\theenumi}{\roman{enumi}}
\item $\Ext^1_{\Com(X_\gamma\Qcoh)}
(\rF^{\bu\,(\gamma)},\rC^{\bu\,(\gamma)})=0$
for all\/ $\gamma\in\boZ_{\ge1}$;
\item $\Ext^1_{\Com(X_\gamma\Qcoh)}(\rF^{\bu\,(\gamma)},\>
\cI_\gamma\ot_{\cO_{X_\gamma}}\rC^{\bu\,(\gamma)})=0$
for all\/ $\gamma\in\boZ_{\ge2}$, where\/
$\cI_{\gamma}\subset\cO_{X_\gamma}$ denotes the quasi-coherent sheaf
of ideals of the closed subscheme $X_{\gamma-1}\subset X_\gamma$.
\end{enumerate}
 Then\/ $\Ext^1_{\Com(\fX\Proqc_\fl)}(\rF^\bu,\rC^\bu)=0$.
\end{thm}

\begin{proof}
 The same argument as in the proof of
Theorem~\ref{flat-proqcoh-prosheaves-Ext1-orthogonality-thm} applies.
 One just needs to work with complexes of pro-sheaves and sheaves
(short exact sequences of complexes, morphisms of complexes etc.)\
throughout, instead of single pro-sheaves and sheaves.
\end{proof}

 We need another bit of recollection of known material.

\begin{prop} \label{termwise-cta-cot-orthogty-for-qcoh-sheaves-prop}
 Let $X$ be a quasi-compact semi-separated scheme.
 Then \par
\textup{(a)} for any acyclic complex\/ $\cF^\bu$ in the exact category
of very flat quasi-coherent sheaves $X\Qcoh_\fl$ and any complex of
contraadjusted quasi-coherent sheaves\/ $\cC^\bu$ on $X$, one has\/
$\Ext^1_{\Com(X\Qcoh)}(\cF^\bu,\cC^\bu)=0$; \par
\textup{(b)} for any acyclic complex\/ $\cF^\bu$ in the exact category
of flat quasi-coherent sheaves $X\Qcoh_\fl$ and any complex of
cotorsion quasi-coherent sheaves\/ $\cC^\bu$ on $X$, one has\/
$\Ext^1_{\Com(X\Qcoh)}(\cF^\bu,\cC^\bu)=0$.
\end{prop}

\begin{proof}
 Part~(b) is a scheme generalization of the theorem of Bazzoni,
Cort\'es-Izurdiaga, and Estrada~\cite[Theorem~5.3]{BCE}.
 The formulation can be found in~\cite[Theorem~3.3]{CET}.
 (Notice that the first paragraph of the proof of~\cite[Lemma~3.2]{CET}
is erroneous, though the formulation of the lemma is true;
this mistake in the proof in~\cite{CET} was later corrected
in~\cite[Theorem~6.3(2)] {EGilO}, see also~\cite[Remark~6.7]{EGilO}.)
 See~\cite[Lemma~2.1]{Gil} or~\cite[Lemma~B.1.7]{Pcosh} for
a background exact category-theoretic lemma about extensions of
complexes.

 Another reference is~\cite[Lemma~5.1.1(b)]{Pcosh}.
 For a more general point of view, see~\cite[Proposition~2.1,
Theorem~3.3\,(2)\,$\Leftrightarrow$\,(4), and Example~8.9]{Pfltp}.

 Part~(a) is an easier version of part~(b).
 A direct reference is~\cite[Lemma~5.1.1(a)]{Pcosh}.
 For a more general point of view, see~\cite[Proposition~2.1,
Theorem~3.3\,(2)\,$\Leftrightarrow$\,(4), Corollary~6.5,
Lemma~6.7, and Example~6.12(1)]{Pfltp}.
\end{proof}

\begin{cor} \label{cmplxs-termwise-levelwise-cta-cot-orthogty-cor}
\textup{(a)} Let\/ $\fX$ be an ind-semi-separated reasonable\/
$\aleph_0$\+ind-scheme.
 Then one has\/ $\Ext_{\Com(\fX\Proqc_\fl)}^1(\rF^\bu,\rC^\bu)=0$
for any acyclic complex\/ $\rF^\bu$ in the exact category of very flat
pro-quasi-coherent pro-sheaves\/ $\fX\Proqc_\vfl$ and any complex of
levelwise contraadjusted flat pro-quasi-coherent pro-sheaves\/
$\rC^\bu$ on\/~$\fX$. \par
\textup{(b)} Let\/ $\fX$ be an ind-semi-separated ind-Noetherian\/ 
$\aleph_0$\+ind-scheme.
 Then one has\/ $\Ext_{\Com(\fX\Proqc_\fl)}^1(\rF^\bu,\rC^\bu)=0$
for any acyclic complex\/ $\rF^\bu$ in the exact category of flat
pro-quasi-coherent pro-sheaves\/ $\fX\Proqc_\fl$ and any complex of
levelwise cotorsion flat pro-quasi-coherent pro-sheaves\/ $\rC^\bu$
on\/~$\fX$.
\end{cor}

\begin{proof}
 The proof is similar to
that of Corollary~\ref{levelwise-cta-cot-orthogonality-cor}.
 The argument is based on
Theorem~\ref{cmplxes-offlat-proqcoh-prosh-Ext1-orthogty-thm},
Lemma~\ref{tensor-product-flat-cta-cot-lemma}, and
Proposition~\ref{termwise-cta-cot-orthogty-for-qcoh-sheaves-prop}.
\end{proof}

\end{document}